\documentclass[twoside,12pt]{article}
\usepackage{amssymb}

\usepackage{a4wide}
\usepackage{amsthm}
\usepackage{amsmath}
\usepackage[all]{xy}
\usepackage{enumerate}
\usepackage{fancyhdr}

\newtheorem{theorem}{Theorem}[section]
\newtheorem{proposition}[theorem]{Proposition}
\newtheorem{corollary}[theorem]{Corollary}
\newtheorem{lemma}[theorem]{Lemma}
\theoremstyle{definition}

\newtheorem*{Beweis}{Proof}
\newtheorem{definition}[theorem]{Definition}
\newtheorem{punto}[theorem]{}
\theoremstyle{remark}
\newtheorem{remark}[theorem]{Remark}
\newtheorem{ex}[theorem]{Example}

\CompileMatrices
\swapnumbers
\CompileMatrices
\input{tcilatex}

\begin{document}

\title{Rational Modules for Corings\thanks{%
MSC (2000): 16W30 \newline
Keywords: Corings; Rational Modules; Comodules; Entwining Structures;
Entwined Modules; Doi-Koppinen Structures; Doi-Koppinen Modules.}}
\author{\textbf{Jawad Y. Abuhlail}\thanks{%
Current Address: Department of Mathematical Sciences, Box \# 5046, KFUPM,
31261 Dhahran (Saudi Arabia), \textbf{email:} abuhlail@kfupm.edu.sa} \\
%EndAName
Mathematics Department, Birzeit University\\
Birzeit - Palestine}
\date{}
\maketitle

\begin{abstract}
The so called \emph{dense pairings} were studied mainly by D. Radford in his
work on coreflexive coalegbras over fields. They were generalized in a joint
paper with J. G\'{o}mez-Torrecillas and J. Lobillo to the so called \emph{%
rational pairings} over a commutative ground ring $R$ to study the interplay
between the comodules of an $R$-coalgebra $C$ and the modules of an $R$%
-algebra $A$ that admits an $R$-algebra morphism $\kappa :A\rightarrow
C^{\ast }.$ Such pairings, satisfying the so called $\alpha $\emph{-condition%
}, were called in the author's dissertation \emph{measuring }$\alpha $\emph{%
-pairings} and can be considered as the corner stone in his study of duality
theorems for Hopf algebras over commutative rings. In this paper we lay the
basis of the theory of rational modules of corings extending results on
rational modules for coalgebras to the case of arbitrary ground rings. We
apply these results mainly to categories of entwined modules (e.g.
Doi-Koppinen modules, alternative Doi-Koppinen modules) generalizing results
of Y. Doi , M. Koppinen and C. Menini et al.
\end{abstract}

\vspace*{-1cm}

\section*{Introduction}

Let $(H,A,C)$ be a right-right Doi-Koppinen structure over a commutative
ring $R,$ $\mathcal{M}(H)_{A}^{C}$ the corresponding category of
Doi-Koppinen modules and $A\#^{op}C^{\ast }$ the Koppinen opposite smash
product. If $_{R}C$ is flat, then $\mathcal{M}(H)_{A}^{C}$ is a Grothendieck
category with enough injective objects. A sufficient, however not necessary,
condition for $\mathcal{M}(H)_{A}^{C}$ to embed as a \emph{full }subcategory
of $\mathcal{M}_{A\#^{op}C^{\ast }}$ is the projectivity of $_{R}C$
\cite[Proposition 3.1]{Kop95}. A similar result for a left-right
Doi-Koppinen structure $(H,A,C)$ was obtained by Y. Doi \cite[3.1]{Doi94},
where the corresponding category of Doi-Koppinen modules $_{A}\mathcal{M}%
(H)^{C}$ was shown to be naturally isomorphic to the category of $\#$\emph{%
-rational }$\#(C,A)$\emph{-modules.} In this paper we show that these
results can be obtained under a weaker condition, that $_{R}C$ is \emph{%
locally projective}, as corollaries from the more general theory of rational
modules for \emph{corings} over a (not necessarily commutative) ring.
Moreover, we show that these categories are of type $\sigma \lbrack M],$ the
theory of which is well developed (e.g. \cite{Wis88}). This extends our
results in \cite{AG-TL2001} and \cite{Abu2001} on the category of rational
modules of an $R$-coalgebra. A fundamental tool in our work is the so called
$\alpha $\emph{-condition,} introduced in \cite{AG-TL2001} for commutative
base rings, which proved also to be very helpful in the author's study of
duality theorems for Hopf algebras \cite{Abu2001}.

The concept of a \emph{coring} over an arbitrary ground ring $R$ is due to
M. Sweedler \cite{Swe75} and is a generalization of the concept of a \emph{%
coalgebra} over a commutative ground ring. In the first section we give the
needed definitions as well as the basic properties of the category of
comodules of a coring. We introduce also the category of \emph{measuring
left }(resp. \emph{measuring right}) $R$\emph{-pairings} $\mathcal{P}_{ml}$
(resp. $\mathcal{P}_{mr}$) and the category of \emph{measuring }$R$\emph{%
-pairings }$\mathcal{P}_{m}.$ For each $(\mathcal{A},\mathcal{C})\in
\mathcal{P}_{ml}$ (resp. $(\mathcal{A},\mathcal{C})\in \mathcal{P}_{mr}$) we
consider two right (resp. left) linear topologies on $\mathcal{A}$, namely
the \emph{weak linear topology} $\mathcal{A}[\frak{T}_{ls}^{r}(\mathcal{C})]$
and the $\mathcal{C}$\emph{-adic topology} $\mathcal{T}_{-\mathcal{C}}(%
\mathcal{A})$ (resp. $\mathcal{A}[\frak{T}_{_{ls}}^{l}(\mathcal{C})]$ and $%
\mathcal{T}_{\mathcal{C}-}(\mathcal{A})$) and show that $\mathcal{A}[\frak{T}%
_{ls}^{r}(\mathcal{C})]=\mathcal{T}_{-\mathcal{C}}(\mathcal{A})$ (resp. $%
\mathcal{A}[\frak{T}_{ls}^{l}(\mathcal{C})]=\mathcal{T}_{\mathcal{C}-}(%
\mathcal{A})$).

In the second section we define the \emph{rational modules} of a measuring
left (resp. right) $R$-pairing satisfying the so called $\alpha $\emph{%
-condition.} The main result (Theorem \ref{cor-dicht}) characterizes the
measuring left $R$-pairings $(\mathcal{A},\mathcal{C})$ satisfying the $%
\alpha $-condition as those for which $_{R}\mathcal{C}$ is locally
projective and $\kappa _{P}(\mathcal{A})\subseteq $ $^{\ast }\mathcal{C}$ is
dense (equivalently, those for which $\mathcal{M}^{\mathcal{C}}=\sigma
\lbrack \mathcal{C}_{\mathcal{A}}]=\sigma \lbrack \mathcal{C}_{^{\ast }%
\mathcal{C}}]$). Theorem \ref{cd-r} provides a dual version for measuring
right $R$-pairings. For a measuring left $\alpha $-pairing $(\mathcal{A},%
\mathcal{C})$ we prove for $\mathrm{Rat}^{\mathcal{C}}(\mathcal{M}_{\mathcal{%
A}})=\sigma \lbrack \mathcal{C}_{\mathcal{A}}]$ the important \emph{%
Finiteness Theorem} (\ref{es}). The properties of the right linear topology $%
\mathcal{T}_{-\mathcal{C}}(\mathcal{A})=\mathcal{A}[\frak{T}_{ls}^{r}(%
\mathcal{C})]$ introduced in the first section will be used then to give
topological (besides the algebraic) characterizations of the rational
modules (Proposition \ref{n_rat}).

In the third section we give some applications of our results in the first
and second sections to the category of entwined modules $\mathcal{M}%
_{A}^{C}(\psi )$ corresponding to an entwining structure $(A,C,\psi )$ with $%
_{R}C$ locally projective, where $R$ is a commutative ground ring. Our
observations generalize results of Y. Doi \cite{Doi94} and M. Koppinen \cite
{Kop95} on the category of Doi-Koppinen modules $\mathcal{M}(H)_{A}^{C}$
corresponding to a Doi-Koppinen structure $(H,A,C)$ with $_{R}C$ projective
and results of C. Menini et al. (e.g. \cite{MTW2001}, \cite{MSTW2001}) on
the category of relative Hopf modules $\mathcal{M}_{A}^{H}$ with $_{R}H$
projective.

Throughout this paper $R$ denotes an associative ring with $1_{R}\neq 0_{R}.$
We consider $R$ as a right (and a left) linear topological ring with the
\emph{discrete topology}. With $\mathcal{M}_{R}$ (resp. $_{R}\mathcal{M},$ $%
_{R}\mathcal{M}_{R}$) we denote the category of right $R$-modules (resp.
left $R$-modules, $R$-bimodules). All $R$-modules are assumed to be \emph{%
unital}. For every right (resp. left) $R$-module $M$ we denote by $\vartheta
_{M}^{r}:M\otimes _{R}R\rightarrow M$ (resp. $\vartheta _{M}^{l}:R\otimes
_{R}M\rightarrow M$) the canonical isomorphisms. With $R^{op}$ we denote the
\emph{opposite ring}. For a right $R$-Module $M$ and a left $R$-module $N$
we denote with $\tau :M\otimes _{R}N\rightarrow N\otimes _{R^{op}}M$ the
canonical \emph{twist}\textbf{\ }$\mathbb{Z}$-isomorphism. For a left (resp.
a right) $R$-module $K$ we consider $K$ as a right (resp. a left) module
over its ring of endomorphisms $\mathrm{End}(_{R}K)^{op}$ (resp. $\mathrm{End%
}(K_{R})$) and a left (resp. a right) module over $\mathrm{Biend}(_{R}K):=%
\mathrm{End}(K_{\mathrm{End}(_{R}K)^{op}})$ (resp. $\mathrm{Biend}(K_{R}):=%
\mathrm{End}(_{\mathrm{End}(K_{R})}K)^{op}$), the \emph{ring of
biendomorphisms} of $K$ (e.g. \cite[6.4]{Wis88}).

For an $R$-ring $\mathcal{A}$ and an $\mathcal{A}$-module $M$ we call an $%
\mathcal{A}$-submodule $N\subset M$ $R$\emph{-cofinite}, if $M/N$ is f.g. as
an $R$-module. If $U$ is an $R$-bimodule, then for every right (resp. left) $%
R$-module $L$ we consider $\mathrm{Hom}_{-R}(U,L)$ (resp. $\mathrm{Hom}%
_{R-}(U,L)$) as a right (resp. left) $R$-module through $(fr)(u):=f(ru)$
(resp. $(rf)(u):=f(ur)$). Moreover $U^{\ast }:=\mathrm{Hom}_{-R}(U,R)$
(resp. $^{\ast }U:=\mathrm{Hom}_{R-}(U,R)$) is an $R$-bimodule through the
right (resp. left) $R$-action given above and the left (resp. right) $R$%
-action given by $(\widetilde{r}f)(u):=\widetilde{r}f(u)$ (resp. $(f%
\widetilde{r})(u):=f(u)\widetilde{r}$). With $^{\ast }U^{\ast }:=\mathrm{Hom}%
_{R-R}(U,R)$ we denote the $R$-bimodule of $R$-bilinear maps from $U$ to $R.$

\section{Preliminaries}

By{\normalsize \ }an \emph{associative }$R$\emph{-ring}\ we mean an $R$%
-bimodule with an $R$-bilinear map (\emph{multiplication}) $\mu \mathcal{_{A}%
}:\mathcal{A}\otimes _{R}\mathcal{A}\rightarrow \mathcal{A},$ such that
\begin{equation*}
\mu \mathcal{_{A}}\circ (\mu \mathcal{_{A}}\otimes id\mathcal{_{A}})=\mu
\mathcal{_{A}}\circ (id\mathcal{_{A}}\otimes \mu \mathcal{_{A}}).
\end{equation*}
If there exists an $R$-bilinear map $\eta \mathcal{_{A}}:R\rightarrow
\mathcal{A},$ such that
\begin{equation*}
\mu \mathcal{_{A}}\circ (id\mathcal{_{A}}\otimes \eta \mathcal{_{A}}%
)=\vartheta _{\mathcal{A}}^{r}\text{ and }\mu \mathcal{_{A}}\circ (\eta
\mathcal{_{A}}\otimes id\mathcal{_{A}})=\vartheta _{\mathcal{A}}^{l},
\end{equation*}
then we call $\eta \mathcal{_{A}}$ the \emph{unity} of $\mathcal{A}.$ If $%
\mathcal{A}$ and $\mathcal{B}$ are $R$-rings (with unities $\eta _{\mathcal{A%
}},$ $\eta _{\mathcal{B}}$), then an $R$-bilinear map $f:\mathcal{A}%
\rightarrow \mathcal{B}$ is called a \emph{morphism of }$R$\emph{-rings}, if
$f\circ \mu _{\mathcal{A}}=\mu _{\mathcal{B}}\circ (f\otimes f)$ (and $%
f\circ \eta _{\mathcal{A}}=\eta _{\mathcal{B}}$). The set of morphisms of $R$%
-rings form $\mathcal{A}$ to $\mathcal{B}$ is denoted by $\mathrm{Rng}_{R}(%
\mathcal{A},\mathcal{B}).$ The category of associative $R$-rings with
unities will be denoted by $\mathbf{Rng}_{R}.$

Dual to $R$-rings are $R$-corings presented by M. Sweedler \cite{Swe75}:

\begin{punto}
A \emph{coassociative }$R$\emph{-coring} is an $R$-bimodule $\mathcal{C}$
associated with an $R$-bilinear map (\emph{comultiplication}) $\Delta _{%
\mathcal{C}}:\mathcal{C}\rightarrow \mathcal{C}\otimes _{R}\mathcal{C},$
such that the following diagram is commutative:
\begin{equation*}
\xymatrix{{\mathcal {C}} \ar^(.45){\Delta_{\mathcal C}}[rr]
\ar_(.45){\Delta_{\mathcal C}}[d] & & {\mathcal {C}} \otimes_{R} {\mathcal
{C}} \ar^(.45){id \otimes \Delta_{\mathcal C}}[d]\\ {\mathcal {C}}
\otimes_{R} {\mathcal {C}} \ar_(.45){\Delta_{\mathcal C} \otimes id}[rr] & &
{\mathcal {C}} \otimes_{R} {\mathcal {C}} \otimes_{R} {\mathcal {C}} }
\end{equation*}
If there exists an $R$-bilinear map $\varepsilon _{\mathcal{C}}:\mathcal{C}%
\rightarrow R${\normalsize ,} so that the following diagram
\begin{equation*}
\xymatrix{ & {\mathcal {C}} \ar^(.45){\Delta _{\mathcal {C}}}[d] & \\ {R}
\otimes_{R} {\mathcal {C}} \ar[ur]^(.45){\vartheta _{\mathcal {C}} ^l} &
{\mathcal {C}} \otimes_{R} {\mathcal {C}} \ar^(.45){\varepsilon _{\mathcal
{C}} \otimes id}[l] \ar_(.45){id \otimes \varepsilon _{\mathcal {C}}}[r] &
{\mathcal {C}} \otimes_{R} {R} \ar[ul]_(.45){\vartheta _{\mathcal {C}} ^r} }
\end{equation*}
is commutative, then we call $\varepsilon _{\mathcal{C}}$ the \emph{counity}
of $\mathcal{C}.$ Unless the contrary is assumed, we make the convention
that an $R$-coring has a counit. If $R$ is a commutative ring, then $R$%
-corings are called $R$-coalgebras (see \cite{Swe69}).
\end{punto}

Let $(\mathcal{C},\Delta )$ be an $R$-coring. For $c\in \mathcal{C}$ we use
Sweedler-Heyneman's $\sum $-notation:
\begin{equation*}
\Delta (c)=\sum c_{1}\otimes c_{2}\in \mathcal{C}\otimes _{R}\mathcal{C}.
\end{equation*}
Moreover we define $\Delta _{n}$ inductively as $\Delta _{1}:=\Delta $ and
\begin{equation*}
\Delta _{n}:=(\Delta \otimes id^{n-1})\circ \Delta _{n-1}:\mathcal{C}%
\rightarrow \mathcal{C}^{n+1},\text{ }c\longmapsto \sum c_{1}\otimes
...\otimes c_{n+1}\text{ for }n\geq 2.
\end{equation*}

\begin{punto}
\textbf{The category of $R$-corings. }For two $R$-corings $(\mathcal{C}%
,\Delta _{\mathcal{C}}),$ $(\mathcal{D},\Delta _{\mathcal{D}})$ (with
counities $\varepsilon _{\mathcal{C}},\varepsilon _{\mathcal{D}}$) we call
an $R$-bilinear map $f:\mathcal{D}\rightarrow \mathcal{C}$\ an $R$\emph{%
-coring morphism}, if the following diagram
\begin{equation*}
\xymatrix{{\mathcal {D}} \ar^(.45){f}[rr] \ar_(.45){\Delta_{\mathcal
{D}}}[d] & & {\mathcal {C}} \ar^(.45){\Delta_{{\mathcal {C}}}}[d] \\
{\mathcal {D}} \otimes_{R} {\mathcal {D}} \ar_(.45){f \otimes f}[rr] & &
{\mathcal {C}} \otimes_{R} {\mathcal {C}}}
\end{equation*}
is commutative (and $\varepsilon _{\mathcal{C}}\circ f=\varepsilon _{%
\mathcal{D}}$). The set of $R$-coring morphisms from $\mathcal{D}$ to $%
\mathcal{C}$ is denoted by $\mathrm{Cog}_{R}(\mathcal{D},\mathcal{C}).$ The
category of coassociative $R$-corings with counities is denoted by $\mathbf{%
Corng}_{R}.$
\end{punto}

\begin{definition}
Let $M$ be a right (resp. a left) $R$-module and $N\subset M$ an $R$%
-submodule. We call $N\hookrightarrow M$ $W$\emph{-pure} for some left
(resp. right) $R$-module $W,$ if $0\rightarrow N\otimes _{R}W\rightarrow
M\otimes _{R}W$ (resp. $0\rightarrow W\otimes _{R}N\rightarrow W\otimes
_{R}M $) is exact in $_{\mathbb{Z}}\mathcal{M}.$ We call $N\hookrightarrow M$
\emph{pure} (in the sense of Cohn), if $N\hookrightarrow M$ is $W$-pure for
every left (resp. right) $R$-module $W.$ If $M$ is an $R$-bimodule and $%
N\subset M$ is an $R$-subbimodule, then we call $_{R}N_{R}\subset $ $%
_{R}M_{R}$ pure, if $N\subset M$ is pure as a right as well as a left $R$%
-submodule.
\end{definition}

\begin{punto}
Let $(\mathcal{C},\Delta _{\mathcal{C}},\varepsilon _{\mathcal{C}})$ be an $%
R $-coring. We call a \emph{pure }$R$-subbimodule $\mathcal{D}\subseteq
\mathcal{C}$ an $R$\emph{-subcoring}, if $\Delta _{\mathcal{C}}(\mathcal{D}%
)\subseteq \mathcal{D}\otimes _{R}\mathcal{D}.$ We call an $R$-subbimodule $%
K\subset \mathcal{C}$ a $\mathcal{C}$\emph{-bicoideal}\textbf{\ }(resp. a $%
\mathcal{C}$\emph{-coideal}), if $\Delta _{\mathcal{C}}(K)\subseteq \func{Im}%
(\mathcal{C}\otimes _{R}K)\cap \func{Im}(K\otimes _{R}\mathcal{C})$ (resp. $%
\Delta _{\mathcal{C}}(K)\subseteq K\wedge K:=\mathrm{Ke}(\mathcal{C}\otimes
_{R}\mathcal{C}\rightarrow \mathcal{C}/K\otimes _{R}\mathcal{C}/K)$). A
right (resp. a left) $R$-submodule $K\subset \mathcal{C}$ is called a \emph{%
right} $\mathcal{C}$\emph{-coideal}\textbf{\ }(resp. a \emph{left} $\mathcal{%
C}$\emph{-coideal}), if $\Delta _{\mathcal{C}}(K)\subset \mathrm{\func{Im}}%
(K\otimes _{R}\mathcal{C})$ (resp. $\Delta _{\mathcal{C}}(K)\subset \func{Im}%
(\mathcal{C}\otimes _{R}K)$).
\end{punto}

\begin{punto}
\textbf{The Dual rings of a coring. }(\cite{Guz85}) Let $\mathcal{C}$ be a
\emph{coassociative }$R$-coring. Then\linebreak $^{\ast }\mathcal{C}:=(%
\mathrm{Hom}_{R-}(\mathcal{C},R),\star _{l})$ is an \emph{associative }$R$%
-ring, where
\begin{equation*}
(f\star _{l}g)(c)=\sum g(c_{1}f(c_{2}))\text{ for all }f,g\in \text{ }^{\ast
}\mathcal{C}\text{ and }c\in \mathcal{C}{\normalsize ;}
\end{equation*}
$\mathcal{C}^{\ast }:=(\mathrm{Hom}_{-R}(\mathcal{C},R),\star _{r})$ is an
\emph{associative} $R$-ring, where
\begin{equation*}
(f\star _{r}g)(c)=\sum f(g(c_{1})c_{2})\text{ for all }f,g\in \mathcal{C}%
^{\ast }\text{ and }c\in \mathcal{C}{\normalsize ;}
\end{equation*}
$^{\ast }\mathcal{C}^{\ast }:=(\mathrm{Hom}_{R-R}(\mathcal{C},R),\star )$ is
an \emph{associative }$R$-ring, where
\begin{equation*}
(f\star g)(c)=\sum g(c_{1})f(c_{2})\text{ for all }f,g\in \text{ }^{\ast }%
\mathcal{C}^{\ast }\text{ and }c\in \mathcal{C}{\normalsize .}
\end{equation*}
If $\mathcal{C}$ has counity $\varepsilon _{\mathcal{C}},$ then $\varepsilon
_{\mathcal{C}}$ is a unity for $^{\ast }\mathcal{C},$ $\mathcal{C}\mathbf{%
^{\ast }}$ and $^{\ast }\mathcal{C}^{\ast }.$
\end{punto}

\begin{punto}
Let $(\mathcal{A},\mu )$ be an $R$-ring (\emph{not necessarily with unity}).
A right $\mathcal{A}$-module $M$ will be called \emph{unital} (resp. $%
\mathcal{A}$\emph{-faithful}), if $M\mathcal{A}=M$ (resp. the canonical map $%
\rho _{M}:M\rightarrow \mathrm{Hom}_{\mathbb{-}R}(\mathcal{A},M)$ is
injective). For right $\mathcal{A}$-modules $(M,\rho _{M}),(N,\rho _{N})$ an
$R$-linear map $f:M\rightarrow N$ will be called $\mathcal{A}$\emph{-linear}%
, if $\rho _{N}\circ f=(\mathcal{A},f)\circ \rho _{M}.$ The set of $\mathcal{%
A}$-linear maps from $M$ to $N$ will be denoted by $\mathrm{Hom}_{-\mathcal{A%
}}(M,N).$ With $\mathcal{M}_{\mathcal{A}}$\ (resp. $\widetilde{\mathcal{M}}_{%
\mathcal{A}}$) we denote the category of unital (resp. $\mathcal{A}$%
-faithful) right $\mathcal{A}$-modules. Analogously we define left $\mathcal{%
A}$-modules. For left $\mathcal{A}$-modules $M$ and $N,$ the set of $%
\mathcal{A}$-linear maps from $M$ to $N$ is denoted by $\mathrm{Hom}_{%
\mathcal{A}-}(M,N).$ The category of unital (resp. $\mathcal{A}$-faithful)
left{\normalsize \ }$\mathcal{A}$-modules is denoted by $_{\mathcal{A}}%
\mathcal{M}$ (resp. $\mathcal{_{A}}\widetilde{\mathcal{M}}$).

Let $\mathcal{A},\mathcal{B}$ be $R$-rings. A $(\mathcal{B},\mathcal{A})$%
-bimodule $M$ will be called \emph{unital}\textbf{\ }(resp. \textbf{$(%
\mathcal{B},\mathcal{A})$}\emph{-faithful}), if $_{\mathcal{B}}M$ and $M_{%
\mathcal{A}}$ are unital (resp. if $M$ is $\mathcal{B}$-faithful and $%
\mathcal{A}$-faithful). For $(\mathcal{B},\mathcal{A})$-bimodules $M$ and $N$
we denote the set of $\mathcal{B}$-linear $\mathcal{A}$-linear maps from $M$
to $N$ (called $(\mathcal{B},\mathcal{A})$\emph{-bilinear}) by $\mathrm{Hom}%
_{\mathcal{B}-\mathcal{A}}(M,N).$ The category of unital (resp. \emph{$(%
\mathcal{B},\mathcal{A})$}-faithful) $(\mathcal{B},\mathcal{A})$-bimodules
and $(\mathcal{B},\mathcal{A})$-bilinear maps is denoted by $_{\mathcal{B}}%
\mathcal{M}_{\mathcal{A}}$ (resp. $_{\mathcal{B}}\widetilde{\mathcal{M}}_{%
\mathcal{A}}$).
\end{punto}

\newpage Dual to modules of $R$-rings are comodules of $R$-corings:

\begin{punto}
Let $(\mathcal{C},\Delta )$ be an $R$-coring (\emph{not necessarily with
counit}). A right $\mathcal{C}$\emph{-comodule} is a right $R$-module $M$
associated with an $R$-linear map ($\mathcal{C}$\emph{-coaction})
\begin{equation*}
\varrho _{M}:M\rightarrow M\otimes _{R}\mathcal{C},\text{ }m\mapsto \sum
m_{<0>}\otimes m_{<1>},
\end{equation*}
such that the following diagram is commutative
\begin{equation*}
\xymatrix{M \ar^(.4){\varrho _M}[rr] \ar_(.45){\varrho _M}[d] & & M
\otimes_{R} {\mathcal{C}} \ar^(.45){id _M \otimes \Delta}[d] \\ M
\otimes_{R} {\mathcal{C}} \ar_(.4){\varrho _M \otimes id_{\mathcal{C}}}[rr]
& & M \otimes_{R} {\mathcal {C}} \otimes_{R} {\mathcal{C}}}
\end{equation*}
If $\varrho _{M}$ is injective, then we call $M$ \emph{counital}. For right $%
\mathcal{C}$-comodules $M,N\in \mathcal{M}^{\mathcal{C}}${\normalsize \ }we
call an $R$-linear map $f:M\rightarrow N$ a $\mathcal{C}$\emph{-comodule
morphism}{\normalsize \ }(or $\mathcal{C}$\emph{-colinear}), if the
following diagram is commutative
\begin{equation*}
\xymatrix{M \ar[rr]^{f} \ar[d]_{\varrho_M} & & N \ar[d]^{\varrho_N}\\ M
\otimes_{R} {\mathcal{C}} \ar[rr]_{f \otimes id_{\mathcal{C}}} & & N
\otimes_{R} {\mathcal{C}} }
\end{equation*}
The set of $\mathcal{C}$-colinear maps from $M$ to $N$ is denoted by $%
\mathrm{Hom}^{\mathcal{C}}(M,N).$ The category of counital right $\mathcal{C}
$-comodules and $\mathcal{C}$-colinear maps is denoted by $\mathcal{M}^{%
\mathcal{C}}.$ For a right $\mathcal{C}$-comodule $N$ we call a right $R$%
-submodule $K\subset N$ a $\mathcal{C}$\emph{-subcomodule}, if $(K,\varrho
_{K})\in \mathcal{M}^{\mathcal{C}}$ and the embedding $K\overset{\iota _{K}}{%
\hookrightarrow }N$ is $\mathcal{C}$-colinear.

Analogously we define the left $\mathcal{C}$-comodules. For two left $%
\mathcal{C}$-comodules $M,N$ we denote by $^{\mathcal{C}}\mathrm{Hom}(M,N)$
the set of $\mathcal{C}$-colinear maps from $M$ to $N.$ The catgeory of
counital\emph{\ }left $\mathcal{C}$-comodules will be denoted by $^{\mathcal{%
C}}\mathcal{M}.$
\end{punto}

\begin{lemma}
\label{inj-kounit}\emph{(}Compare \emph{\cite[Lemma 1.1.]{CC94})} Let $(%
\mathcal{C},\Delta )$ be an $R$-coring. If $\mathcal{C}$ has counity $%
\varepsilon ,$ then a right $\mathcal{C}$-comodule $(M,\varrho _{M})$ is
counital iff $\vartheta _{M}^{r}\circ (id_{M}\otimes \varepsilon )\circ
\varrho _{M}=id_{M}.$ Hence, if $M$ is counital, then $\varrho _{M}$ is a
splitting monomorphism.
\end{lemma}

\begin{definition}
Let $\frak{C}$ be a category with finite limits and finite colimits. A
functor $F:\frak{C}\rightarrow \frak{D}$ is called \emph{left-exact}\textbf{%
\ (}resp. \emph{right-exact}\textbf{)}, if $F$ preserves finite limits
(resp. finite colimits). $F$ is called \emph{exact}, if it's left-exact and
right-exact.
\end{definition}

For forthcoming reference we list some properties of the category of right
comodules of an $R$-coring $C.$ These can be found in several references
(e.g. Brzezi\'{n}ski (2002); Caenepeel et. al. (2002); Wisbauer (2002)).

\begin{proposition}
\label{propert}Let $(\mathcal{C},\Delta _{\mathcal{C}},\varepsilon _{%
\mathcal{C}})$ be an $R$-coring.

\begin{enumerate}
\item  We have a covariant functor
\begin{equation*}
-\otimes _{R}\mathcal{C}:\mathcal{M}_{R}\rightarrow \mathcal{M}^{\mathcal{C}%
},\text{ }M\mapsto (M\otimes _{R}\mathcal{C},id_{M}\otimes \Delta _{\mathcal{%
C}}).
\end{equation*}
Moreover $-\otimes _{R}\mathcal{C}$ is right adjoint to the forgetful
functor $\mathcal{F}:\mathcal{M}^{\mathcal{C}}\rightarrow \mathcal{M}_{R}$
and left adjoint to the functor $\mathrm{Hom}^{\mathcal{C}}(\mathcal{C},-):%
\mathcal{M}^{\mathcal{C}}\rightarrow \mathcal{M}_{R}.$

\item  $-\otimes _{R}\mathcal{C}:\mathcal{M}_{R}\rightarrow \mathcal{M}^{%
\mathcal{C}}$ is left-exact and $\mathcal{F}:\mathcal{M}^{\mathcal{C}%
}\rightarrow \mathcal{M}_{R}$ is right-exact \emph{(}since a right adjoint
functor preserves limits and a left adjoint functor preserves colimits \emph{%
(}e.g. \emph{\cite[Proposition 16.4.6]{Sch72})}.

\item  $\mathcal{F}$ is exact iff $_{R}\mathcal{C}$ is flat iff $-\otimes
_{R}\mathcal{C}:\mathcal{M}_{R}\rightarrow \mathcal{M}^{\mathcal{C}}$
preserves injective objects \emph{(}for the commutative case see \emph{
\cite[Proposition 8]{Wisch75})}.

\item  The category $\mathcal{M}^{\mathcal{C}}$ is cocomplete and has
cokernels. The direct limits and direct sums are formed in $\mathcal{M}_{R}.$
Moreover $-\otimes _{R}\mathcal{C}:\mathcal{M}_{R}\rightarrow \mathcal{M}^{%
\mathcal{C}}$ respects direct limits \emph{(}i.e. direct sums and cokernels%
\emph{)}.

\item  If $Q$ is a cogenerator in $\mathcal{M}_{R},$ then $Q\otimes _{R}%
\mathcal{C}$ is a cogenerator in $\mathcal{M}^{\mathcal{C}}.$ In particular $%
\mathcal{M}^{\mathcal{C}}$ has a cogenerator.

\item  If $_{R}\mathcal{C}$ is flat, then $\mathcal{M}^{\mathcal{C}}$ is a
Grothendieck catgeory with enough injective objects. In this case $-\otimes
_{R}\mathcal{C}:\mathcal{M}_{R}\rightarrow \mathcal{M}^{\mathcal{C}}$
respects inverse limits \emph{(}i.e. direct products and kernels\emph{)}.

\item  If $R_{R}$ is a cogenerator, then $\mathcal{C}$ is a cogenerator in $%
\mathcal{M}^{\mathcal{C}}.$ If $_{R}\mathcal{C}$ is flat and $R_{R}$ is
injective, then $\mathcal{C}$ is injective in $\mathcal{M}^{\mathcal{C}}.$
\end{enumerate}
\end{proposition}

\begin{remark}
If $\mathcal{C}$ is an $R$-coring such that $\mathcal{M}^{\mathcal{C}}$ is
Grothendieck, then $_{R}\mathcal{C}$ need not be flat. A counter example is
\cite[Example 1.1.]{EG-TL}. This shows that the conjecture of M. Wischnewsky
\cite[Conjecture 14]{Wisch75} is false for corings.
\end{remark}

\begin{punto}
\textbf{Bicomodules. }Let $(M,\varrho _{M}^{\mathcal{C}})$ be a right $%
\mathcal{C}$-comodule, $(M,\varrho _{M}^{\mathcal{D}})$ be a left $\mathcal{D%
}$-comodule and consider the left $\mathcal{D}$-comodule $(M\otimes _{R}%
\mathcal{C},\varrho _{M}^{\mathcal{D}}\otimes id_{\mathcal{C}})$ (resp. the
right $\mathcal{C}$-comodule $(\mathcal{D}\otimes _{R}M,id_{\mathcal{D}%
}\otimes \varrho _{M}^{\mathcal{C}})$). We call $M$ a $(\mathcal{D},\mathcal{%
C})$\emph{-bicomodule}, if $\varrho _{M}^{\mathcal{C}}:M\rightarrow M\otimes
_{R}\mathcal{C}$ is $\mathcal{D}$-colinear (equivalently, if $\varrho _{M}^{%
\mathcal{D}}:M\rightarrow \mathcal{D}\otimes _{R}M$ is $\mathcal{C}$%
-colinear). For $(\mathcal{D},\mathcal{C})$-bicomodules $M,N$ we call a $%
\mathcal{D}$-colinear $\mathcal{C}$-colinear map $f:M\rightarrow N$ a $(%
\mathcal{D},\mathcal{C})$\emph{-bicomodule morphism}\textbf{\ }(or $(%
\mathcal{D},\mathcal{C})$\emph{-bicolinear}). We say a $(\mathcal{D},%
\mathcal{C})$-bicomodule is \emph{counital}, if it's counital as a left $%
\mathcal{D}$-comodule and as a right $\mathcal{C}$-comodule. The category of
counital $(\mathcal{D},\mathcal{C})$-bicomodules and $(\mathcal{D},\mathcal{C%
})$-bicolinear maps is denoted by $^{\mathcal{D}}\mathcal{M}^{\mathcal{C}}.$
\end{punto}

\subsection*{The weak linear topology}

Next we introduce the categories of left (resp. right) $R$-pairings $%
\mathcal{P}_{l}$ (resp. $\mathcal{P}_{r}$) and the category of $R$-pairings $%
\mathcal{P}.$ For each left (resp. right) $R $-pairing $P=(V,W)$ we consider
$V$ with a right (resp. a left) linear topology, the \emph{weak linear
topology} $V[\frak{T}_{ls}^{r}(W)]$ (resp. $V[\frak{T}_{ls}^{l}(W)]$)\emph{\
}(e.g. \cite[10.3]{Kot66}, \cite{Rad73}).

\begin{punto}
\textbf{The category of }$R$\textbf{-pairings. }With a \emph{left}\textbf{\ }%
$R$\emph{-pairing} we denote a right $R$-module $V$ and a left $R$-module $W$
with an $R$-linear map $\kappa _{P}:V\rightarrow $ $^{\ast }W$ (equivalently
$\chi _{P}:W\rightarrow V^{\ast }$). For left $R$-pairings $(V,W)$ and $%
(V^{\prime },W^{\prime })$ a \emph{morphism of left }$R$\emph{-pairings}%
\textbf{\ }$(\xi ,\theta ):(V^{\prime },W^{\prime })\rightarrow (V,W)$
consists of a morphism of right $R$-modules $\xi :V\rightarrow V^{\prime }$
and a morphism of left $R$-modules $\theta :W^{\prime }\rightarrow W,$ such
that
\begin{equation}
<\xi (v),w^{\prime }>=<v,\theta (w^{\prime })>\text{ for all }v\in V\text{
and }w^{\prime }\in W^{\prime }.  \label{comp}
\end{equation}
The left $R$-pairings with the morphisms described above (and the usual
composition of pairings) build a category, which we denote by $\mathcal{P}%
_{l}.$ Analogously we define the category of right $R$-pairings $\mathcal{P}%
_{r}.$

With an $R$\emph{-pairing }we denote $R$-bimodules $V$ and $W$ with an $R$%
-bilinear map $\kappa _{P}:V\rightarrow $ $^{\ast }W^{\ast }$ (equivalently $%
\chi _{P}:W\rightarrow $ $^{\ast }V^{\ast }$). If $(V,W)$ and $(V^{\prime
},W^{\prime })\;$are $R$-pairings, then a \emph{morphism of }$R$\emph{%
-pairings}\textbf{\ }$(\xi ,\theta ):(V^{\prime },W^{\prime })\rightarrow
(V,W)$ consists of $R$-bilinear maps $\xi :V\rightarrow V^{\prime }$ and $%
\theta :W^{\prime }\rightarrow W$ with the compatibility condition (\ref
{comp}). The category of $R$-pairings is denoted by $\mathcal{P}.$
\end{punto}

\begin{punto}
\textbf{The finite topology.} Let $E$ be a right (resp. a left) $R$-module, $%
W$ a set and identify the direct product $E^{W}$ with the set of all maps
from $W$ to $E.$ If we consider $E$ with the discrete topology and the right
(resp. the left) $R$-module $E^{W}$ with the \emph{product topology}, then
the induced \emph{linear} topology on an $R$-submodule $Z\subseteq E^{W}$ is
called the \emph{finite topology} and has a neighbourhood system of $0_{Z}:$%
\begin{equation*}
\mathcal{B}{\normalsize _{f}(0_{Z})}:=\{\mathrm{An}(F)|\text{ }%
F=\{w_{1},...,w_{k}\}\subset W\text{ is a finite subset}\}.
\end{equation*}
\end{punto}

\begin{punto}
Let $P=(V,W)$ be a left (resp. a right) $R$-pairing and consider the right
(resp. the left) $R$-submodule $^{\ast }W\subset R^{W}$ (resp. $W^{\ast
}\subset R^{W}$) with the \emph{finite topology}. Then there is on $V$ a
right (resp. a left) linear topology, the \emph{weak linear topology} $V[%
\frak{T}_{ls}^{r}(W)]$ (resp. $V[\frak{T}_{ls}^{l}(W)]$), such that $\kappa
_{P}:V\rightarrow $ $^{\ast }W$ (resp. $\kappa _{P}:V\rightarrow W^{\ast }$)
is continuous. The neighbourhood system of $0_{V}$ w.r.t. this topology is
given by
\begin{equation*}
\mathcal{B}{\normalsize _{f}(0_{V})}=\{F^{\bot }:=\kappa _{P}^{-1}(\mathrm{An%
}(F))|\text{ }F=\{w_{1},...,w_{k}\}\subset W\text{ is a finite subset}\}.
\end{equation*}
The closure $\overline{X}$ of any subset $X\subseteq V$ is then given by
\begin{equation}
\overline{X}=\bigcap \{X+F^{\bot }\mid F\subset W\text{ is a finite subset}%
\}.  \label{X-closure}
\end{equation}
Hence $V[\frak{T}_{_{ls}}^{r}(W)]$ (resp. $V[\frak{T}_{_{ls}}^{l}(W)]$) is
\emph{Hausdorff} iff $V\overset{\kappa _{P}}{\hookrightarrow }$ $^{\ast }W$
(resp. $V\overset{\kappa _{P}}{\hookrightarrow }W^{\ast }$) is an embedding.
\end{punto}

\qquad Next we state some properties of the weak linear topology without
proof. For the proofs and other details the interested reader may refer to
\cite{Abu}. For the case of a commutative ground ring a reference is
\cite[Anhang]{Abu2001}.

\begin{lemma}
\label{orth-clos}Let $P=(V,W)$ be a left $R$-pairing and consider $V$ with
the weak linear topology $V[\frak{T}_{ls}^{r}(W)].$

\begin{enumerate}
\item  $\overline{X}\subseteq X^{\bot \bot }$ for any subset $X\subset V.$
Consequently every orthogonally closed right $R$-submodule of $V$ is closed.

\item  If $R_{R}$ is noetherian, then all open $R$-submodules of $V$ are $R$%
-cofinite.

\item  Let $R_{R}$ be artinian.

\begin{enumerate}
\item  Every $R$-cofinite closed right $R$-submodule $X\subset V$ is open.

\item  Let $X\subset Y\subset V$ be right $R$-submodules. If $X\subset V$ is
closed and $R$-cofinite, then $Y\subset V$ is also closed and $R$-cofinite.
\end{enumerate}
\end{enumerate}
\end{lemma}

We call the ring $R$ a \emph{QF ring}, if $R_{R}$ (equivalently $_{R}R$) is
noetherian and a cogenerator (e.g. \cite[48.15]{Wis88}).

The following result characterizes the closed and the open $R$-submodules of
$V$ w.r.t. the weak linear topology:

\begin{proposition}
\label{lrs-bet}Let $P=(V,W)$ be a left $R$-pairing and consider $V$ with the
weak linear topology $V[\frak{T}_{ls}^{r}(W)].$ Assume $R_{R}$ to be an
injective cogenerator.

\begin{enumerate}
\item  The closure of a right $R$-submodule $X\subseteq V$ is given by $%
\overline{X}=X^{\bot \bot }.$

\item  For right $R$-submodules $X\subset Y\subseteq V$ we have: $X$ is
dense in $Y$ iff $X^{\bot }=Y^{\bot }.$ If $W\hookrightarrow V^{\ast },$
then $X\subset V$ is dense iff $X^{\bot }=0.$

\item  The class of closed right $R$-submodules of $V$ is given by
\begin{equation*}
{\normalsize \{K^{\bot }|}\text{ }K\subset W\text{ is an \emph{arbitrary}
left }R\text{-submodule}\}.
\end{equation*}

\item  If $R$ is a QF-ring and $W\overset{\chi _{P}}{\hookrightarrow }%
V^{\ast }$ is an embedding, then the class of open right $R$-submodules of $%
V $ is given by
\begin{equation*}
{\normalsize \{K^{\bot }|}\text{ }K\subset W\text{ is a \emph{f.g.} left }R%
\text{-submodule}\}.
\end{equation*}
\end{enumerate}
\end{proposition}

\begin{lemma}
\label{f*-clos}Let $W,W^{\prime }$ be left $R$-modules and consider the
right $R$-modules $^{\ast }W$ and $^{\ast }W^{\prime }$ with the finite
topology. Let $\theta \in \mathrm{Hom}_{R-}(W^{\prime },W).$ If $R$ is a
QF-ring, then $\mathrm{Ke}(\theta ^{\ast }(X))=\theta ^{-1}(\mathrm{Ke}(X))$
for every right $R$-submodule $X\subset $ $^{\ast }W.$
\end{lemma}

\subsection*{The $\mathcal{C}$-adic topology}

\qquad We introduce now the category of \emph{measuring left $R$-pairings }$%
\mathcal{P}_{ml}.$ For every $(\mathcal{A},\mathcal{C})\in \mathcal{P}_{ml}$
we define the $\mathcal{C}$\emph{-adic topology }$\mathcal{T}_{-\mathcal{C}}(%
\mathcal{A})$ $($see \cite{AW97}, \cite{Ber94}), which we prove to coincide
with the linear weak topology $\mathcal{A}[\frak{T}_{_{ls}}^{r}(\mathcal{C}%
)].$

\begin{punto}
\textbf{The category of measuring }$R$-\textbf{pairings. }If $\mathcal{C}$
is an $R$-coring and $\mathcal{A}$ is an $R$-ring with a morphism of $R$%
-rings $\kappa :\mathcal{A}\rightarrow $ $^{\ast }\mathcal{C},$ $a\mapsto
\lbrack c\mapsto <a,c>],$ then we call $P:=(\mathcal{A},\mathcal{C})$ a
\emph{measuring left }$R$\emph{-pairing} (the terminology is inspired by
\cite[Definition, Page 138]{Swe69}). For measuring left $R$-pairings $(%
\mathcal{A},\mathcal{C}),$ $(\mathcal{B},\mathcal{D})$ we say a morphism of
left $R$-pairings $(\xi ,\theta ):(\mathcal{B},\mathcal{D})\rightarrow (%
\mathcal{A},\mathcal{C})$ is a \emph{morphism of measuring left }$R$\emph{%
-pairings}, if $\xi :\mathcal{A}\rightarrow \mathcal{B}$ is a morphism of $R$%
-rings and $\theta :\mathcal{D}\rightarrow \mathcal{C}$ is a morphism of $R$%
-corings. The measuring left $R$-pairings with the morphisms described above
build a subcategory $\mathcal{P}_{ml}\subset \mathcal{P}_{l}.$

A \emph{measuring right }$R$\emph{-pairing} $P=(\mathcal{A},\mathcal{C})$
consists of an $R$-ring $\mathcal{A}$ and an $R$-coring $\mathcal{C}$ with a
morphism of $R$-rings $\kappa _{P}:\mathcal{A}\rightarrow \mathcal{C}^{\ast
}.$ If $\mathcal{A}$ is an $R$-ring with a morphism of $R$-rings $\kappa
_{P}:\mathcal{A}\rightarrow $ $^{\ast }\mathcal{C}^{\ast },$ then we call $(%
\mathcal{A},\mathcal{C})$ a \emph{measuring }$R$\emph{-pairing}. The
category of measuring right $R$-pairings (resp. measuring $R$-pairings) is
denoted by $\mathcal{P}_{mr}$ (resp. $\mathcal{P}_{m}$).
\end{punto}

\begin{punto}
If $P=(\mathcal{A},\mathcal{C})$ is a measuring left (right) $R$-pairing,
then $\mathcal{C}$ becomes a right (a left) $\mathcal{A}$-module with $%
\mathcal{A}$-action given by
\begin{equation}
c\leftharpoonup a:=\sum c_{1}<a,c_{2}>\text{ (resp. }a\rightharpoonup
c:=\sum <a,c_{1}>c_{2}\text{).}  \label{A-act}
\end{equation}
If $P=(\mathcal{A},\mathcal{C})$ is a measuring $R$-pairing, then $\mathcal{C%
}$ is an $\mathcal{A}$-bimodule with the right and the left $\mathcal{A}$%
-actions in (\ref{A-act}).
\end{punto}

The following example was communicated to the author by Tomasz
Brzezi\'{n}ski:

\begin{ex}
\label{Tomasz}Let $\mathcal{C}$ be a \emph{coseparable} $R$-coring (i.e.
there exists a $\mathcal{C}$-bicolinear map $\pi :\mathcal{C}\otimes _{R}%
\mathcal{C}\rightarrow \mathcal{C}$ with $\pi \circ \Delta _{\mathcal{C}%
}=id_{\mathcal{C}}$ \cite{Guz89}, equivalently there exists a \emph{%
cointegral}\textbf{\ }$\gamma \in \mathrm{Hom}_{R-R}(\mathcal{C}\otimes _{R}%
\mathcal{C},R),$ such that $\gamma \circ \Delta _{\mathcal{C}}=\varepsilon _{%
\mathcal{C}}$ and
\begin{equation*}
\sum c_{1}\gamma (c_{2}\otimes c^{\prime })=\sum \gamma (c\otimes
c_{1}^{\prime })c_{2}^{\prime }\text{ for all }c,c^{\prime }\in \mathcal{C}.
\end{equation*}
(Brzezinski, 2002, Theorem 3.5, Corollary 3.6). Then $\mathcal{C}$ is a (non
unital) $R$-ring with multiplication
\begin{equation*}
\mu :\mathcal{C}\otimes _{R}\mathcal{C}\rightarrow \mathcal{C},\text{ }%
c\otimes \widetilde{c}\mapsto \sum c_{1}\gamma (c_{2}\otimes \widetilde{c})
\end{equation*}
and therefore $P:=(\mathcal{C},\mathcal{C})$ is a measuring left $R$-pairing
with
\begin{equation*}
\kappa _{P}:\mathcal{C}\rightarrow {}^{\ast }\mathcal{C},\text{ }c\mapsto
\lbrack c^{\prime }\mapsto \gamma (c^{\prime }\otimes c)].
\end{equation*}
\end{ex}

\begin{punto}
\textbf{Subpairings. }Let $P=(\mathcal{A},\mathcal{C})$ and assume $P\in
\mathcal{P}_{ml}$ (resp. $P\in \mathcal{P}_{mr},$ $P\in \mathcal{P}_{m}$), $%
\mathcal{J}\vartriangleleft \mathcal{A}$ an $\mathcal{A}$-ideal, $\mathcal{D}%
\overset{\iota }{\hookrightarrow }\mathcal{C}$ an $R$-subcoring with $<%
\mathcal{J},\mathcal{D}>=0$ and put $Q:=(\mathcal{A}/\mathcal{J},\mathcal{D}%
).$ Then $Q\in \mathcal{P}_{ml}$ (resp. $Q\in \mathcal{P}_{mr},$ $Q\in
\mathcal{P}_{m}$) and $(\pi ,\iota ):(\mathcal{A}/\mathcal{J},\mathcal{D}%
)\rightarrow (\mathcal{A},\mathcal{C})$ is a morphism in $\mathcal{P}_{ml}$
(resp. in $\mathcal{P}_{mr},$ in $\mathcal{P}_{m}$). We call $Q$ a\textbf{\ }%
\emph{measuring }$R$\emph{-subpairing}\textbf{\ }of $P$ and write $Q\subset
P.$ For every $R$-subcoring $\mathcal{D}\subseteq \mathcal{C}$ we have the
\emph{two-sided} $\mathcal{A}$-ideal $\mathcal{D}^{\bot }\subseteq \mathcal{A%
}$ and so $(\mathcal{A}/\mathcal{D}^{\bot },\mathcal{D})\subseteq (\mathcal{A%
},\mathcal{C})$ is a measuring $R$-subpairing. In particular, to every
measuring left (resp. right) $R$-pairing $(\mathcal{A},\mathcal{C}),$ is $(%
\mathcal{A}/\mathcal{C}^{\bot },\mathcal{C})\subseteq (^{\ast }\mathcal{C},%
\mathcal{C})$ (resp. $(\mathcal{A}/\mathcal{C}^{\bot },\mathcal{C})\subseteq
(\mathcal{C}^{\ast },\mathcal{C})$) a \emph{non-degenerate}\textbf{\ }%
measuring left (resp. right) $R$-subpairing.
\end{punto}

\begin{punto}
\textbf{Subgenerators}. Let $\mathcal{A}$ be an $R$-ring and $K$ an $%
\mathcal{A}$-module. We say that an $\mathcal{A}$-module $N$ is $K$\emph{%
-subgenerated}, if $N$ is isomorphic to a submodule of a $K$-generated $%
\mathcal{A}$-module (equivalently, if $N$ is kernel of a morphism between $K$%
-generated $\mathcal{A}$-modules). For a right $\mathcal{A}$-module $K,$ we
denote by $\sigma \lbrack K\mathcal{_{A}}]$ the \emph{full }subcategory of $%
\mathcal{M_{A}}$ whose objects are the $K$-subgenerated right $\mathcal{A}$%
-modules. For every right $\mathcal{A}$-module $M$%
\begin{equation}
\mathrm{Sp}(\sigma \lbrack K\mathcal{_{A}}],M):=\sum \{f(N)\mid \text{ }f\in
\mathrm{Hom}_{-\mathcal{A}}(N,M),\text{ }N\in \sigma \lbrack K\mathcal{_{A}}%
]\}  \label{Spur}
\end{equation}
is the largest $K$-subgenerated right $\mathcal{A}$-submodule of $M.$
Moreover $\sigma \lbrack K_{\mathcal{A}}]$ is the \emph{smallest }%
Grothendieck full subcategory of $\mathcal{M}_{{\mathcal{A}}}$ that contains
$K.$ For the well developed theory of categories of this type the reader is
referred to \cite{Wis88}.
\end{punto}

\begin{punto}
\textbf{The }$\mathcal{C}$\textbf{-adic topology}.\label{C-ad} Let $P=(%
\mathcal{A},\mathcal{C})\in \mathcal{P}_{ml}$ and consider $\mathcal{C}$
with the right $\mathcal{A}$-module structure through ($\leftharpoonup $) in
(\ref{A-act}). Then the class of right $\mathcal{A}$-ideals
\begin{equation*}
\mathcal{B}{\normalsize _{-\mathcal{C}}}(0_{\mathcal{A}}):=\{(0_{\mathcal{C}%
}:F)|\text{ }F\subset \mathcal{C}\text{ is a finite subset}\}
\end{equation*}
is a neighbourhood system of $0_{\mathcal{A}}$ for a right linear topology,
the $\mathcal{C}$\emph{-adic topology} $\mathcal{T}_{-\mathcal{C}}(\mathcal{A%
})$, and $(\mathcal{A},\mathcal{T}_{-\mathcal{C}}(\mathcal{A}))$ is a right
linear topological $R$-ring. A right $\mathcal{A}$-ideal $I\vartriangleleft
_{r}\mathcal{A}$ is open w.r.t. $\mathcal{T}_{-\mathcal{C}}(\mathcal{A})$
iff $\mathcal{A}/I$ is $\mathcal{C}$-subgenerated. If $(\mathcal{A},\frak{T}%
) $ is a right linear topological ring, then the category of $(\mathcal{A},%
\frak{T})$-discrete modules coincides with $\sigma \lbrack \mathcal{C}_{%
\mathcal{A}}]$ iff $\frak{T}=\mathcal{T}_{-\mathcal{C}}(\mathcal{A}).$ We
refer mainly to \cite{AW97} and \cite{Ber94} for detailed investigation of
this topology.
\end{punto}

\begin{lemma}
\label{C:ls=ad}Let $P=(\mathcal{A},\mathcal{C})\in \mathcal{P}_{ml}.$ The
weak linear topology $\mathcal{A}[\frak{T}_{ls}^{r}(\mathcal{C})]$ and the $%
\mathcal{C}$-adic topology $\mathcal{T}_{-\mathcal{C}}(\mathcal{A})$
coincide. In particular $(\mathcal{A},\mathcal{A}[\frak{T}_{ls}^{r}(\mathcal{%
C})])$ is a right linear topological $R$-ring and a right $\mathcal{A}$%
-module $M$ is $(\mathcal{A},\mathcal{A}[\frak{T}_{ls}^{r}(\mathcal{C})])$%
-discrete iff $M_{\mathcal{A}}$ is $\mathcal{C}$-subgenerated.
\end{lemma}

\begin{Beweis}
Let $U$ be a neighbourhood of $0_{\mathcal{A}}$ w.r.t. $\mathcal{A}[\frak{T}%
_{ls}^{r}(\mathcal{C})].$ Then there exists a f.g. left $R$-submodule $%
K\subset \mathcal{C},$ such that $K^{\bot }\subseteq U.$ But we have then
for every $a\in (0_{\mathcal{C}}:K)$ and $c\in K:$%
\begin{equation*}
<a,c>=<a,\sum \varepsilon (c_{1})c_{2}>=\varepsilon (\sum
c_{1}<a,c_{2}>)=\varepsilon (c\leftharpoonup a)=0,
\end{equation*}
and so $(0_{\mathcal{C}}:K)\subseteq K^{\bot }\subseteq U,$ i.e. $U$ is a
neighbourhood of $0_{\mathcal{A}}$ w.r.t. $\mathcal{T}_{-\mathcal{C}}(%
\mathcal{A}).$ On the other hand, let $U$ be a neighbourhood of $0_{\mathcal{%
A}}$ w.r.t. $\mathcal{T}_{-\mathcal{C}}(\mathcal{A}).$ Then there exists a
finite set $F=\{c_{1},...,c_{k}\}\subset \mathcal{C},$ such that $(0_{%
\mathcal{C}}:F)\subseteq U.$ Assume $\Delta
(c_{i})=\sum\limits_{j=1}^{n_{i}}c_{ij}\otimes \widetilde{c}_{ij}$ for $%
i=1,...,k$ and put $K:=\sum\limits_{i=1}^{k}\sum\limits_{j=1}^{n_{i}}R%
\widetilde{c}_{ij}.$ Then $K^{\bot }\subseteq (0_{\mathcal{C}}:F)\subseteq
U, $ i.e. $U$ is a neighbourhood of $0_{\mathcal{A}}$ w.r.t. $\mathcal{A}[%
\frak{T}_{ls}^{r}(\mathcal{C})].$ So $\mathcal{A}[\frak{T}_{ls}^{r}(\mathcal{%
C})]=\mathcal{T}_{-\mathcal{C}}(\mathcal{A})$ and the last statement follows
now from \ref{C-ad}.$\blacksquare $
\end{Beweis}

\qquad As a corollary of Proposition \ref{lrs-bet} (2) and Lemma \ref
{C:ls=ad} we get

\begin{corollary}
\label{dense-per}If $P=(\mathcal{A},\mathcal{C})\in \mathcal{P}_{ml}$, then
for every right $R$-submodule $Y\subset \mathcal{A}$ the following are
equivalent:

\begin{enumerate}
\item  $Y\subset \mathcal{A}$ is dense w.r.t. $\mathcal{A}[\frak{T}_{ls}^{r}(%
\mathcal{C})].$

\item  $Y\subset \mathcal{A}$ is dense w.r.t. $\mathcal{T}_{-\mathcal{C}}(%
\mathcal{A}).$

If $R_{R}$ is an injective cogenerator and $\mathcal{C}\hookrightarrow
\mathcal{A}^{\ast }$, then \emph{(1)} $\&$\ \emph{(2)} are equivalent to

\item  $Y^{\bot }:=\{c\in \mathcal{C}|<a,c>=0$ for every $a\in Y\}=0.$
\end{enumerate}
\end{corollary}

\subsection*{The $\protect\alpha $-condition}

\qquad We introduce now the category of left (resp. right) $\alpha $%
-pairings $\mathcal{P}_{l}^{\alpha }$ (resp. $\mathcal{P}_{r}^{\alpha }$)
and the category of $\alpha $-pairings $\mathcal{P}^{\alpha }.$

\begin{punto}
\label{alp-cond}We say that a left $R$-pairing $P=(V,W)\;$satisfies the (%
\emph{left})\textbf{\ }$\alpha $\emph{-condition}\textbf{, }or is \emph{left
}$\alpha $\emph{-pairing}, if for every right $R$-module $M$ the following
map is injective:
\begin{equation}
\alpha _{M}^{P}:M\otimes _{R}W\ \rightarrow \mathrm{Hom}_{-R}(V,M),\text{ }%
\sum m_{i}\otimes w_{i}\mapsto \lbrack v\mapsto \sum m_{i}<v,w_{i}>].  \notag
\end{equation}
A right $R$-pairing $P=(V,W)$ is said to satisfy the \textbf{(}\emph{right}%
\textbf{) }$\alpha $\emph{-condition}, or is a\emph{\ right }$\alpha $\emph{%
-pairing}, if for every left\emph{\ }$R$-module $M,$ the canonical map $%
\alpha _{M}^{P}:W\otimes _{R}M\rightarrow \mathrm{Hom}_{R-}(V,M)$ is
injective. With $\mathcal{P}_{l}^{\alpha }\subset \mathcal{P}_{l}$ (resp. $%
\mathcal{P}_{r}^{\alpha }\subset \mathcal{P}_{r}$) we denote the full
subcategory, whose objects satisfy the $\alpha $-condition. We call a left
(resp. a right) $R$-pairing $P=(V,W)$ \emph{dense}, if $\kappa
_{P}(V)\subseteq $ $^{\ast }W$ (resp. $\kappa _{P}(V)\subseteq W^{\ast }$)
is dense w.r.t. the finite topology.

With $\mathcal{P}_{ml}^{\alpha }\subset \mathcal{P}_{ml}$ (resp. $\mathcal{P}%
_{m}^{\alpha }\subset \mathcal{P}_{m}$) we denote the \emph{full }%
subcategory of \emph{measuring left }$\alpha $\emph{-pairings} (resp. \emph{%
measuring right }$\alpha $\emph{-pairings}). If $P\in \mathcal{P}%
_{ml}^{\alpha }$ (resp. $P\in \mathcal{P}_{mr}^{\alpha }$) and $Q\subset P$
is a measuring $R$-subpairing, then $Q$ satisfies the $\alpha $-condition as
well (see Proposition \ref{rp-rp} (1-b) below).

We say a left (resp. a right) $R$-module $W$ satisfies the $\alpha $\emph{%
-condition}, if the left $R$-pairing $(^{\ast }W,W)$ (resp. the right $R$%
-pairing $(W^{\ast },W)$) satisfies the $\alpha $-condition, i.e. if $_{R}W$
(resp. $W_{R}$) is \emph{universally torsion free} in the sense of G.
Garfinkel \cite{Gar76}.
\end{punto}

\begin{punto}
\textbf{Locally projective modules}. An $R$-module $W$ is called \emph{%
locally projective} (in the sense of B. Zimmermann-Huisgen \cite{Z-H76}), if
for every diagram
\begin{equation*}
\xymatrix{0 \ar[r] & F \ar@{.>}[dr]_{g' \circ \iota} \ar[r]^{\iota} & W
\ar[dr]^{g} \ar@{.>}[d]^{g'} & & \\ & & L \ar[r]_{\pi} & N \ar[r] & 0}
\end{equation*}
with exact rows and $F$ f.g.: for every $R$-linear map $g:W\rightarrow N,$
there exists an $R$-linear map $g^{\prime }:W\rightarrow L$, such that the
entstanding parallelogram is commutative. Note that every projective $R$%
-module is locally projective.
\end{punto}

Analog to \cite[Theorem 3.2]{Gar76} and \cite[Theorem 2.1]{Z-H76} we get the
following characterizations of the $R$-modules satisfying the $\alpha $%
-condition:

\begin{lemma}
\label{loc-p}A left \emph{(}resp. a right\emph{)} $R$-module $W$ satisfies
the $\alpha $-condition iff $_{R}W$ \emph{(}resp. $W_{R}$\emph{)} is locally
projective.
\end{lemma}

\begin{remark}
\label{flat}Let $P=(V,W)\in \mathcal{P}_{l}^{\alpha }.$ Then $W\subset
V^{\ast },$ in particular $_{R}W$ is $R$-cogenerated. If $M$ is any right $R$%
-module, then we have for every $R$-submodule $N\subset M$ the commutative
diagram
\begin{equation*}
\xymatrix{ N \otimes_{R} W \ar[rr]^{\alpha _N ^P} \ar[d]_{\iota _N \otimes
id_W} & & {\rm Hom}_{-{R}} (V,N) \ar@{^{(}->}[d] \\ M \otimes_{R} W
\ar[rr]_{\alpha_M ^P} & & {\rm Hom}_{-{R}} (V,M)}
\end{equation*}
By assumption $\alpha _{N}^{P}$ is injective and so $N\subset M$ is $W$%
-pure. Since $M$ and $N$ are arbitrary, we conclude that $_{R}W$ is flat. If
$_{R}R$ is perfect, then $_{R}W$ is projective. In particular every locally
projective left $R$-module (over a left perfect ring) is flat (projective).
So over perfect rings projectivity and local projectivity coincide.
\end{remark}

\begin{lemma}
\label{q-2}Let $P=(V,W)$ be a left $\alpha $-pairing. If $L$ is a right $R$%
-module and $K\subset L$ is an $R$-submodule, then we have for every $\sum
l_{i}\otimes w_{i}\in L\otimes _{R}W:$%
\begin{equation*}
\sum l_{i}\otimes w_{i}\in K\otimes _{R}W\Longleftrightarrow \sum
l_{i}<v,w_{i}>\in K\text{ for all }v\in V.
\end{equation*}
\end{lemma}

\begin{Beweis}
By Remark \ref{flat}\ $_{R}W$ is flat and so we get the commutative diagram
with exact rows
\begin{equation*}
\xymatrix{0 \ar[r] & K \otimes_{R} W \ar[rr]^{\iota_K \otimes id_W}
\ar@{^{(}->}[d]^{\alpha_K ^P} & & L \otimes_{R} W \ar[rr]^{\pi \otimes id_W}
\ar@{^{(}->}[d]^{\alpha_L ^P} & & L/K \otimes_{R} W \ar@{^{(}->}[d]^{\alpha
_{L/K} ^P} \ar[r] & 0\\ 0 \ar[r] & {\rm Hom} _{-{R}}(V,K)
\ar[rr]_{(V,\iota_K)} & & {\rm Hom} _{-{R}}(V,L) \ar[rr]_{(V,\pi)} & & {\rm
Hom} _{-{R}}(V,L/K) & }
\end{equation*}
Clearly $\sum l_{i}<v,w_{i}>\in K$ for every $v\in V$ iff $\sum l_{i}\otimes
w_{i}\in \mathrm{Ke}((V,\pi )\circ \alpha _{L}^{P})=\mathrm{Ke}(\alpha
_{L/K}^{P}\circ (\pi \otimes id_{W}))=\mathrm{Ke}(\pi \otimes
id_{W})=K\otimes _{R}W.\blacksquare $
\end{Beweis}

\qquad Analog to the commutative case \cite[Proposition 2.1.7]{Abu2001} we
get

\begin{proposition}
\label{rp-rp}

\begin{enumerate}
\item  Let $P=(V,W)$ be a left $R$-pairing.

\begin{enumerate}
\item  Let $W^{\prime }\subset W$ be an $R$-submodule and consider the
induced left $R$-pairing $P^{\prime }:=(V,W^{\prime }).$ If $P^{\prime }\in
\mathcal{P}_{l}^{\alpha }$, then $W^{\prime }\subset W$ is pure. If $P\in
\mathcal{P}_{l}^{\alpha },$ then $P^{\prime }\in \mathcal{P}_{l}^{\alpha }$
iff $W^{\prime }\subset W$ is pure.

\item  Let $V^{\prime }\subset V,$ $W^{\prime }\subset W$ be $R$-submodules
with $<V^{\prime },W^{\prime }>=0$ and consider the left $R$-pairing $%
Q:=(V/V^{\prime },W^{\prime }).$ If $P\in \mathcal{P}_{l}^{\alpha }$, then $%
Q\in \mathcal{P}_{l}^{\alpha }$ iff $W^{\prime }\subset W$ is pure.
\end{enumerate}

\item  Let $Q=(Y,W)$ be a left $R$-pairing, $V$ a right $R$-module, $\xi
:V\rightarrow Y$ an $R$-linear map, $P:=(V,W)$ the induced left $R$-pairing
and consider the following statements:

(i) $Q\in \mathcal{P}_{l}^{\alpha }$ and $P$ is dense;

(ii) $Q\in \mathcal{P}_{l}^{\alpha }$ and $\xi (V)\subset Y$ is dense w.r.t.
$Y[\frak{T}_{_{ls}}^{r}(W)];$

(iii) $P\in \mathcal{P}_{l}^{\alpha };$

(iv) $Q\in \mathcal{P}_{l}^{\alpha }$ and $W\overset{\chi _{P}}{%
\hookrightarrow }V^{\ast }$ is an embedding.

The following implications are always true: \emph{(i)} $\Longrightarrow $
\emph{(ii)} $\Longrightarrow $ \emph{(iii)} $\Longrightarrow $ (iv). If $%
R_{R}$ is an injective cogenerator, then \emph{(i)-(iv)}\ are equivalent.
\end{enumerate}
\end{proposition}

\qquad The proof of the following result is similar to that of
\cite[Proposition 2.5]{AG-TL2001}:

\begin{lemma}
\label{p-2}Let $V,W$ be $R$-bimodules.

\begin{enumerate}
\item  If $P=(V,W),$ $P^{\prime }=(V^{\prime },W^{\prime })$ are left $%
\alpha $-pairings, then $P\otimes _{l}P^{\prime }:=(V^{\prime }\otimes
_{R}V,W\otimes _{R}W^{\prime })$ is a left $\alpha $-pairing, where
\begin{equation*}
\kappa _{P\otimes _{l}P^{\prime }}(v^{\prime }\otimes v)(w\otimes w^{\prime
})=<v,w<v^{\prime },w^{\prime }>>=<<v^{\prime },w^{\prime }>v,w>.
\end{equation*}

\item  If $P=(V,W),$ $P^{\prime }=(V^{\prime },W^{\prime })$ are right $%
\alpha $-pairings, then $P\otimes _{r}P^{\prime }:=(V\otimes _{R}V^{\prime
},W^{\prime }\otimes _{R}W)$ is a right $\alpha $-pairing, where
\begin{equation*}
\kappa _{P^{\prime }\otimes _{r}P}(v\otimes v^{\prime })(w^{\prime }\otimes
w)=<v,<v^{\prime },w^{\prime }>w>=<v<v^{\prime },w^{\prime }>,w>.
\end{equation*}
\end{enumerate}
\end{lemma}

\section{Rational modules}

In this section we define the $\mathcal{C}$-\emph{rational }$\mathcal{A}$%
\emph{-modules} associated with a measuring left (resp. right) $R$-pairing $(%
\mathcal{A},\mathcal{C})$ satisfying the $\alpha $-condition and prove the
main result in this paper, namely Theorem \ref{cor-dicht} (resp. its dual
version Theorem \ref{cd-r}).

\begin{remark}
Let $P=(\mathcal{A},{\mathcal{C}})$ be a measuring left $R$-pairing. For
every ${}$right $R$-module $M,$ $\mathrm{Hom}_{-R}(\mathcal{A},M)$ is a
right $\mathcal{A}$-module through $(fa)(a^{\prime })=f(aa^{\prime })$ for
all $a,a^{\prime }\in \mathcal{A}$ and $\alpha _{M}^{P}:M\otimes _{R}{%
\mathcal{C}}\rightarrow \mathrm{Hom}_{-R}(\mathcal{A},M)$ is $\mathcal{A}$%
-linear. If moreover $M$ is a right $\mathcal{A}$-module, then the canonical
map $\rho _{M}:M\rightarrow \mathrm{Hom}_{-R}(\mathcal{A},M)$ is $\mathcal{A}
$-linear.
\end{remark}

\begin{punto}
\label{rat-dar}Let $\mathcal{A}$ be an $R$-ring (\emph{not necessarily with
unity}),\ $P=(\mathcal{A},\mathcal{C})$ a measuring left $\alpha $-pairing
and $M$ an $\mathcal{A}$-faithful right $\mathcal{A}$-module. Put $\mathrm{%
Rat}^{\mathcal{C}}(M_{\mathcal{A}}):=\rho _{M}^{-1}(M\otimes _{R}\mathcal{C}%
),$ i.e. $m\in \mathrm{Rat}^{\mathcal{C}}(M_{\mathcal{A}})$ iff there exists
a uniquely determined element $\sum m_{i}\otimes c_{i},$ such that $ma=\sum
m_{i}<a,c_{i}>$ for every $a\in \mathcal{A}.$ We call $M_{\mathcal{A}}$ $%
\mathcal{C}$\emph{-rational}, if $\mathrm{Rat}^{\mathcal{C}}(M_{\mathcal{A}%
})=M.$ In this case we get an $R$-linear map $\varrho _{M}:=(\alpha
_{M}^{P})^{-1}\circ \rho _{M}:M\rightarrow M\otimes _{R}\mathcal{C}.$
\end{punto}

\qquad Analog to the commutative case (e.g. \cite[Proposition 2.9]{G-T98})
we get

\begin{lemma}
\label{clos}Let $P=(\mathcal{A},\mathcal{C})$ be a measuring left $\alpha $%
-pairing \emph{(}$\mathcal{A}$ \emph{not necessarily with unity)}. For every
$(M,\rho _{M})\in $ $\widetilde{\mathcal{M}}_{\mathcal{A}}$ we have:

1. $\mathrm{Rat}^{\mathcal{C}}(M_{\mathcal{A}})\subseteq M$ is an $\mathcal{A%
}$-submodule.

2. For every $\mathcal{A}$-submodule $N\subset M$ we have $\mathrm{Rat}^{%
\mathcal{C}}(N_{\mathcal{A}})=N\cap \mathrm{Rat}^{\mathcal{C}}(M_{\mathcal{A}%
}).$

3. $\mathrm{Rat}^{\mathcal{C}}(\mathrm{Rat}^{\mathcal{C}}(M_{\mathcal{A}}))=%
\mathrm{Rat}^{\mathcal{C}}(M_{\mathcal{A}}).$

4. For every $L\in $ $\widetilde{\mathcal{M}}_{\mathcal{A}}$ and $f\in $ $%
\mathrm{Hom}_{-\mathcal{A}}(M,L)$ we have $f(\mathrm{Rat}^{\mathcal{C}}(M_{%
\mathcal{A}}))\subseteq \mathrm{Rat}^{\mathcal{C}}(L_{\mathcal{A}}).$
\end{lemma}

\qquad Analog to the commutative case \cite[Folgerung 2.2.10]{Abu2001} we
have

\begin{remark}
Let $P=(\mathcal{A},\mathcal{C})$ be a measuring left $\alpha $-pairing and
consider the embedding $\mathcal{C}\overset{\chi _{P}}{\hookrightarrow }%
\mathcal{A}^{\ast }.$ Since $\chi _{P}$ is $\mathcal{A}$-linear, we have $%
\mathcal{C}\overset{\chi _{P}}{\hookrightarrow }\mathrm{Rat}^{\mathcal{C}}((%
\mathcal{A}^{\ast })_{\mathcal{A}})$ by Lemma \ref{clos} (4). If $f\in
\mathrm{Rat}^{\mathcal{C}}((\mathcal{A}^{\ast })_{\mathcal{A}})$ with $%
\varrho (f)=\sum f_{i}\otimes c_{i},$ then we have for every $a\in \mathcal{A%
}$%
\begin{equation*}
f(a)=(fa)(1_{\mathcal{A}})=\sum f_{i}(1_{\mathcal{A}})<a,c_{i}>=\chi
_{P}(\sum f_{i}(1_{\mathcal{A}})c_{i})(a),
\end{equation*}
i.e. $f=\chi _{P}(\sum f_{i}(1_{\mathcal{A}})c_{i}).$ So $\chi _{P}$ is an
isomorphism.
\end{remark}

For a left (resp. right) measuring $\alpha $-pairing $(\mathcal{A},\mathcal{C%
})$ we denote the category of $\mathcal{C}$-rational right (resp. left) $%
\mathcal{A}$-modules and $\mathcal{A}$-linear maps by $\mathrm{Rat}^{%
\mathcal{C}}(\widetilde{\mathcal{M}}_{\mathcal{A}})$ (resp. $^{\mathcal{C}}%
\mathrm{Rat}(_{\mathcal{A}}\widetilde{\mathcal{M}}$). The subcategory of
unital $\mathcal{C}$-rational right (resp. left) $\mathcal{A}$-modules will
be denoted by $\mathrm{Rat}^{\mathcal{C}}(\mathcal{M}_{\mathcal{A}})$ (resp.
$^{\mathcal{C}}\mathrm{Rat}(_{\mathcal{A}}\mathcal{M})$).

\begin{lemma}
\label{co-rat}Let $P=(\mathcal{A},\mathcal{C})$ be a measuring left $R$%
-pairing \emph{(}$\mathcal{A}$ \emph{not necessarily with unity)}.

\begin{enumerate}
\item  If $(M,\varrho _{M})$ is a right $\mathcal{C}$-comodule, then $M$
becomes a right $\mathcal{A}$-module through
\begin{equation*}
\rho _{M}:=M\overset{\varrho _{M}}{\rightarrow }M\otimes _{R}\mathcal{C}%
\overset{\alpha _{M}^{P}}{\rightarrow }\mathrm{Hom}_{-R}(\mathcal{A},M).
\end{equation*}
If $\mathcal{A}$ has unity and $M$ is counital, then $M_{\mathcal{A}}$ is
unital \emph{(}and\ so\ $\mathcal{A}$-faithful\emph{)}.

\item  Let $(M,\varrho _{M}),(N,\varrho _{N})$ be right $\mathcal{C}$%
-comodules and consider the induced structures of right $\mathcal{A}$%
-modules $(M,\rho _{M}),$ $(N,\rho _{N}).$ If $f:M\rightarrow N$ is $%
\mathcal{C}$-colinear, then $f$ is $\mathcal{A}$-linear.

\item  Let $N$ be a right $\mathcal{C}$-comodule, $K\subset N$ a $\mathcal{C}
$-subcomodule and consider the induced right $\mathcal{A}$-module
structures\ $(N,\rho _{N}),$ $(K,\rho _{K}).$ Then $K\subset N$ is an $%
\mathcal{A}$-submodule.
\end{enumerate}
\end{lemma}

\begin{Beweis}
\begin{enumerate}
\item  Consider the left $R$-pairing $P\otimes _{l}P:=(\mathcal{A}\otimes
_{R}\mathcal{A},\mathcal{C}\otimes _{R}\mathcal{C}).$ For every right $%
\mathcal{A}$-module $(M,\rho _{M})$ we have the diagram
\begin{equation}
\xymatrix{M \ar[rrr]^(.45){\rho_M} \ar[ddd]_{\rho_M} \ar@{=}[dr] & & & {\rm
Hom} ({\mathcal {A}},M) \ar[ddd]^{(\mu,M)} \\ { } & M
\ar[d]_(.45){\varrho_M} \ar[r]^(.45){\varrho_M} & M \otimes_{R} {\mathcal C}
\ar[ur]_(.45){\alpha_M ^P} \ar[d]^(.45){id_M \otimes \Delta_{\mathcal C}} &
\\ & M \otimes_{R} {\mathcal {C}} \ar[dl]^(.45){\alpha_M ^P}
\ar[r]_(.45){\varrho_M \otimes id_{\mathcal {C}}} & M \otimes_{R} {\mathcal
{C}} \otimes_{R} {\mathcal {C}} \ar[dr]^(.45){\alpha ^{P \otimes_{l} P} _M}
& \\ {\rm Hom} ({\mathcal {A}},M) \ar[rr]_(.45){({\mathcal {A}},\rho _M)} &
& {\rm Hom} ({\mathcal {A}},{\rm Hom} (A,M)) \ar[r]_(.55){{\varsigma}^r} &
{\rm Hom} ({\mathcal {A}} \otimes_{R} {\mathcal {A}},M) }  \label{comod}
\end{equation}
(where $\varsigma ^{r}$ is the canonical isomorphism). By definition of $%
\rho _{M}$ and $\alpha _{M}^{_{P\otimes _{l}P}}$ all trapezoids are
commutative. Since $(M,\varrho _{M})$ is a right $\mathcal{C}$-comodule, the
inner rectangle is commutative and so the outer rectangle is commutative,
i.e. $(M,\rho _{M})$ is a right $\mathcal{A}$-module.

Assume $M$ to be counital. If $\mathcal{A}$ has unity, then $\kappa _{P}(1_{%
\mathcal{A}})=\varepsilon _{\mathcal{C}}$ and we have for every $m\in M$%
\begin{equation*}
m=\sum m_{<0>}\varepsilon _{\mathcal{C}}(m_{<1>})=m\varepsilon _{\mathcal{C}%
}=m1_{\mathcal{A}}\in M\mathcal{A},
\end{equation*}
i.e. $M_{\mathcal{A}}$ is unital.

\item  Consider the diagram
\begin{equation}
\xymatrix{M \ar[rrrr]^(.45){f} \ar[dd]_(.45){\varrho_M}
\ar[dr]_(.45){\rho_M} & & & & N \ar[dd]^(.45){\varrho_N}
\ar[dl]^(.45){\rho_N} \\ { } & {\rm Hom}_{-{R}}({\mathcal {A}},M)
\ar[rr]^(.45){({\mathcal {A}},f)} & & {\rm Hom}_{-{R}}({\mathcal {A}},N) &
\\ M \otimes_{R} {\mathcal {C}} \ar[ur]^(.45){\alpha_M ^P} \ar[rrrr]_(.45){f
\otimes id_{\mathcal {C}}} & & & & \ar[ul]_(.45){\alpha_N ^P} N \otimes_{R}
{\mathcal {C}} & }  \label{co-lin}
\end{equation}
The lower trapezoid is obviously commutative. The triangles are commutative
by definition of $\rho _{M},\rho _{N}.$ If the outer rectangle is
commutative ($f$ is $\mathcal{C}$-colinear), then the upper trapezoid is
commutative ($f$ is $\mathcal{A}$-linear).

\item  trivial.$\blacksquare $
\end{enumerate}
\end{Beweis}

\newpage

\begin{lemma}
\label{rat-co}Let $P=(\mathcal{A},\mathcal{C})$ be a measuring left $\alpha $%
-pairing \emph{(}$\mathcal{A}$ \emph{not necessarily with unity)}.

\begin{enumerate}
\item  If $(M,\rho _{M})\in \widetilde{\mathcal{M}}\mathcal{_{\mathcal{A}}}$
is $\mathcal{C}$-rational, then $M$ gets a structure of a \emph{counital }%
right $\mathcal{C}$-comodule through
\begin{equation*}
\varrho _{M}:M\overset{\rho _{M}}{\rightarrow }\mathrm{Hom}_{-R}(\mathcal{A}%
,M)\overset{(\alpha _{M}^{P})^{-1}}{\rightarrow }M\otimes _{R}\mathcal{C}.
\end{equation*}

\item  Let $(M,\rho _{M}),$ $(N,\rho _{N})\in \widetilde{\mathcal{M}}%
\mathcal{_{\mathcal{A}}}$ be $\mathcal{C}$-rational and consider the induced
structures of right $\mathcal{C}$-comodules $(M,\varrho _{M}),$ $(N,\varrho
_{N}).$ Then $\mathrm{Hom}^{\mathcal{C}}(M,N)=\mathrm{Hom}_{-\mathcal{A}%
}(M,N).$

\item  Let $(N,\rho _{N})\in \widetilde{\mathcal{M}}\mathcal{_{\mathcal{A}}}$
be $\mathcal{C}$-rational and consider the induced{\normalsize \ }counital
right $\mathcal{C}$-comodule $(N,\varrho _{N}).$ If $K\subset N$ is an $%
\mathcal{A}$-submodule, then $K$ gets a structure of a counital $\mathcal{C}$%
-subcomodule. Moreover $\varrho _{K}=(\varrho _{N})_{|_{K}}.$
\end{enumerate}
\end{lemma}

\begin{Beweis}
\begin{enumerate}
\item  If $(M,\rho _{M})$ is $\mathcal{C}$-rational, then by definition $%
\rho _{M}(M)\subset \alpha _{M}^{P}(M\otimes _{R}\mathcal{C}),$ i.e. $%
\varrho _{M}:=(\alpha _{M}^{P})^{-1}\circ \rho _{M}$ is well defined and we
get the commutative diagram
\begin{equation*}
\xymatrix{ {\rm Hom}_{-{R}}({\mathcal {A}},M) & \\ M
\ar@{^{(}->}[u]^{\rho_M} \ar@{.>}[r]_(.45){\varrho_M} & M \otimes_{R}
{\mathcal {C}} \ar@{^{(}->}[ul]_{\alpha_M ^P} }
\end{equation*}
By assumption $M_{\mathcal{A}}$ is $\mathcal{\mathcal{A}}$-faithful (i.e. $%
\rho _{M}$ is injective) and so $\varrho _{M}$ is injective, i.e. the
induced right $\mathcal{C}$-coaction on $M$ is counital. Consider now
diagram (\ref{comod}).\ By definition of $\varrho _{M}$ and $\alpha
_{M}^{_{P\otimes _{l}P}}$ all trapezoids are commutative. By assumption $M$
is a right $\mathcal{C}$-comodule and so the outer rectangle is commutative.
By Lemma \ref{p-2} (1) $\alpha _{M}^{P\otimes _{l}P}$ is injective and so
the inner rectangle is commutative, i.e. $(M,\varrho _{M})$ is a counital
right $\mathcal{C}$-comodule.

\item  Consider diagram (\ref{co-lin}). The lower trapezoid is obviously
commutative and the triangles are commutative by definition of $\varrho _{M}$
and $\varrho _{N}.$ Moreover $\alpha _{N}^{P}$ is injective and so the outer
rectangle is commutative ($f$ is $\mathcal{C}$-colinear) iff the upper
trapezoid is commutative ($f$ is $\mathcal{A}$-linear), i.e. $\mathrm{Hom}^{%
\mathcal{C}}(M,N)=\mathrm{Hom}_{-\mathcal{A}}(M,N).$

\item  Let $(N,\rho _{N})$ be a $\mathcal{C}$-rational right $\mathcal{A}$%
-module. If $K\subset N$ is an $\mathcal{A}$-submodule, then by Lemma \ref
{clos} (2) $\mathrm{Rat}^{\mathcal{C}}(K_{\mathcal{A}})=K\cap \mathrm{Rat}^{%
\mathcal{C}}(N_{\mathcal{A}})=K$ and so $K_{\mathcal{A}}$ is $\mathcal{C}$%
-rational, hence $K$ is a counital right $\mathcal{C}$-comodule by (1).
Moreover $K\overset{\iota _{K}}{\hookrightarrow }N$ is $\mathcal{A}$-linear
and so\ $\mathcal{C}$-colinear by (2), i.e. $K\subset N$ is a $\mathcal{C}$%
-subcomodule. Note that by assumption and Remark \ref{flat} $_{R}\mathcal{C}$
is flat, hence $\varrho _{K}=(\varrho _{N})_{|_{K}}.\blacksquare $
\end{enumerate}
\end{Beweis}

\begin{punto}
\label{Biend}For every $R$-coring $\mathcal{C}$ we have an isomorphism of $R$%
-rings $(\mathcal{C}^{\ast },\star _{r})\simeq \mathrm{End}^{\mathcal{C}}(%
\mathcal{C})$ via $f\mapsto \lbrack c\mapsto \sum f(c_{1})c_{2}]$ with
inverse $g\mapsto \varepsilon _{\mathcal{C}}\circ g$ (compare Proposition
\ref{propert} (1)). Analogously $(^{\ast }\mathcal{C},\star _{l})\simeq $ $^{%
\mathcal{C}}\mathrm{End}(\mathcal{C})^{op}$ as $R$-rings.

If $P=(\mathcal{A},\mathcal{C})\in \mathcal{P}_{ml}^{\alpha },$ then by
Lemma \ref{rat-co} (2) $\mathrm{End}^{\mathcal{C}}(\mathcal{C)}=\mathrm{End}(%
\mathcal{C}_{\mathcal{A}}\mathcal{)}$ and so
\begin{equation*}
^{\ast }\mathcal{C}\simeq \text{ }^{\mathcal{C}}\mathrm{End}(\mathcal{C}%
)^{op}\subseteq \mathrm{End}(_{\mathcal{C}^{\ast }}\mathcal{C})^{op}=\mathrm{%
End}(_{\mathrm{End}^{\mathcal{C}}(\mathcal{C)}}\mathcal{C})^{op}=\mathrm{End}%
(_{\mathrm{End}(\mathcal{C}_{\mathcal{A}}\mathcal{)}}\mathcal{C})^{op}:=%
\mathrm{Biend}(\mathcal{C}_{\mathcal{A}}),
\end{equation*}
i.e. $(^{\ast }\mathcal{C},\star _{l})$ is isomorphic to an $R$-subring of $%
\mathrm{Biend}(\mathcal{C}_{\mathcal{A}}).$ If moreover $\mathcal{C}_{R}$ is
locally projective, then $^{\mathcal{C}}\mathrm{End}(\mathcal{C})=\mathrm{End%
}(_{\mathcal{C}^{\ast }}\mathcal{C}),$ hence $^{\ast }\mathcal{C}\simeq
\mathrm{Biend}(\mathcal{C}_{\mathcal{A}}).$

On the other hand, if $P\in \mathcal{P}_{mr}^{\alpha },$ then we have
analogously $^{\mathcal{C}}\mathrm{End}(\mathcal{C)}=\mathrm{End}(_{\mathcal{%
A}}\mathcal{C)}$ and so
\begin{equation*}
\mathcal{C}^{\ast }\simeq \mathrm{End}^{\mathcal{C}}(\mathcal{C})\subseteq
\mathrm{End}(\mathcal{C}_{^{\ast }\mathcal{C}})=\mathrm{End}(\mathcal{C}_{^{%
\mathcal{C}}\mathrm{End}(\mathcal{C})^{op}})=\mathrm{End}(\mathcal{C}_{%
\mathrm{End}(_{\mathcal{A}}\mathcal{C})^{op}}):=\mathrm{Biend}(_{\mathcal{A}}%
\mathcal{C}),
\end{equation*}
i.e. $(\mathcal{C}^{\ast },\star _{r})$ is isomorphic to an $R$-subring of $%
\mathrm{Biend}(_{\mathcal{A}}\mathcal{C}).$ If moreover $_{R}\mathcal{C}$ is
locally projective, then $\mathrm{End}^{\mathcal{C}}(\mathcal{C})=\mathrm{End%
}(\mathcal{C}_{^{\ast }\mathcal{C}}),$ hence $\mathcal{C}^{\ast }\simeq
\mathrm{Biend}(_{\mathcal{A}}\mathcal{C}).$

Note that it follows from above, that in case $_{R}\mathcal{C}$ and $%
\mathcal{C}_{R}$ are locally projective we have $^{\ast }\mathcal{C}\simeq
\mathrm{Biend}(\mathcal{C}_{^{\ast }\mathcal{C}})$ and $\mathcal{C}^{\ast
}\simeq \mathrm{Biend}(_{\mathcal{C}^{\ast }}\mathcal{C})$ as $R$-rings.
\end{punto}

The following result generalizes \cite[Theorem 2.2.13]{Abu2001} from the
case of commutative base rings to the case of arbitrary rings:

\begin{proposition}
\label{equal}Let $P=({\mathcal{A}},{\mathcal{C}})$ be a measuring ${}$left $%
R $-pairing \emph{(}${\mathcal{A}}$ not necessarily with unity\emph{)}. If $%
_{R}{\mathcal{C}}$ is locally projective and $\kappa _{P}({\mathcal{A}}%
)\subseteq $ $^{\ast }\mathcal{C}$ is dense, then
\begin{equation}
\begin{tabular}{lllllll}
$\mathcal{M}^{{\mathcal{C}}}$ & $\simeq $ & $\mathrm{Rat}^{{\mathcal{C}}}(%
\widetilde{\mathcal{M}}_{\mathcal{A}})$ & $=$ & $\mathrm{Rat}^{{\mathcal{C}}%
}(\mathcal{M}_{{\mathcal{A}}})$ & $=$ & $\sigma \lbrack {\mathcal{C}}_{{%
\mathcal{A}}}]$ \\
& $\simeq $ & $\mathrm{Rat}^{{\mathcal{C}}}(\widetilde{\mathcal{M}}_{^{\ast }%
\mathcal{C}})$ & $=$ & $\mathrm{Rat}^{{\mathcal{C}}}(\mathcal{M}_{^{\ast }%
\mathcal{C}})$ & $=$ & $\sigma \lbrack {\mathcal{C}}_{^{\ast }\mathcal{C}}].$%
\end{tabular}
\label{MC=}
\end{equation}
\end{proposition}

\begin{Beweis}
\textbf{Step (1). }By Proposition \ref{rp-rp} (2) $({\mathcal{A}},{\mathcal{C%
}})$ satisfies the left $\alpha $-condition. By Lemmata \ref{co-rat} and \ref
{rat-co} we have covariant functors
\begin{equation*}
\begin{tabular}{ll}
$
\begin{tabular}{cccc}
$(-)_{{\mathcal{A}}}:$ & $\mathcal{M}^{{\mathcal{C}}}$ & $\rightarrow $ & $%
\mathrm{Rat}^{{\mathcal{C}}}(\widetilde{\mathcal{M}}_{\mathcal{A}}),$ \\
& $(M,\varrho _{M})$ & $\mapsto $ & $(M,\alpha _{M}\circ \varrho _{M}),$%
\end{tabular}
$ & $
\begin{tabular}{cccc}
$(-)^{{\mathcal{C}}}:$ & $\mathrm{Rat}^{{\mathcal{C}}}(\widetilde{\mathcal{M}%
}_{\mathcal{A}})$ & $\rightarrow $ & $\mathcal{M}^{{\mathcal{C}}},$ \\
& $(M,\rho _{M})$ & $\mapsto $ & $(M,\alpha _{M}^{-1}\circ \rho _{M}),$%
\end{tabular}
$%
\end{tabular}
\end{equation*}
that act as the identity on morphisms. Clearly we have
\begin{equation*}
(-)^{{\mathcal{C}}}\circ \text{ }(-)_{{\mathcal{A}}}=id_{\mathcal{M}^{{%
\mathcal{C}}}},\text{ }(-)_{{\mathcal{A}}}\circ (-)^{{\mathcal{C}}}=id_{%
\mathrm{Rat}^{{\mathcal{C}}}(\widetilde{\mathcal{M}}_{\mathcal{A}})},
\end{equation*}
i.e. $\mathrm{Rat}^{{\mathcal{C}}}(\widetilde{\mathcal{M}}_{\mathcal{A}%
})\simeq \mathcal{M}^{{\mathcal{C}}}.$

\textbf{Step (2). }We show now that every $\mathcal{C}$-rational right $%
\mathcal{A}$-module is unital. Let $(N,\rho _{N})\in \mathrm{Rat}^{{\mathcal{%
C}}}(\widetilde{\mathcal{M}}_{\mathcal{A}})$ and $n\in N$ with $\varrho
_{N}(n)=\sum\limits_{i=1}^{k}n_{i}\otimes c_{i}.$ By assumption $\kappa _{P}(%
{\mathcal{A}})\subseteq $ $^{\ast }\mathcal{C}$ is dense and so there exists
some $a\in {\mathcal{A}},$ such that $\kappa _{P}(a)(c_{i})=\varepsilon _{%
\mathcal{C}}(c_{i})$ for $i=1,...,k.$ So
\begin{equation*}
n=\sum\limits_{i=1}^{k}n_{i}\varepsilon _{\mathcal{C}}(c_{i})=\sum%
\limits_{i=1}^{k}n_{i}<a,c_{i}>=na\in N\mathcal{A}.
\end{equation*}

\textbf{Step (3). }Let $(N,\rho _{N})\in \mathcal{M}^{{\mathcal{C}}}.$ For
every $n\in N$ with $\varrho _{N}(n)=\sum\limits_{i=1}^{k}n_{i}\otimes c_{i}$
we have $\{c_{1},...,c_{k}\}^{\bot }\subseteq (0_{N}:n),$ hence $N_{\mathcal{%
A}}$ is $\mathcal{C}$-subgenerated. By Proposition \ref{propert} (3), Lemma
\ref{rat-co} and Step (1) $\mathcal{M}^{{\mathcal{C}}}$ is a Grothendieck
full subcategory of $\mathcal{M}_{{\mathcal{A}}}$ and so $\mathcal{M}^{{%
\mathcal{C}}}$ $=$ $\sigma \lbrack {\mathcal{C}}_{{\mathcal{A}}}]$ (since $%
\sigma \lbrack {\mathcal{C}}_{{\mathcal{A}}}]$ is the \emph{smallest }such
subcategory of $\mathcal{M}_{{\mathcal{A}}}$ containing $\mathcal{C}$).

\textbf{Step (4).}{\normalsize \ }For $\mathcal{A}=$ $^{\ast }\mathcal{C}$
we get as above
\begin{equation*}
\mathcal{M}{\normalsize ^{{\mathcal{C}}}}\simeq \mathrm{Rat}^{{\mathcal{C}}}(%
\widetilde{\mathcal{M}}_{^{\ast }\mathcal{C}})=\mathrm{Rat}^{{\mathcal{C}}}(%
\mathcal{M}_{^{\ast }\mathcal{C}})=\sigma \lbrack {\mathcal{C}}_{^{\ast }%
\mathcal{C}}].\blacksquare
\end{equation*}
\end{Beweis}

\qquad We are now ready to present the main result in this article, namely

\begin{theorem}
\label{cor-dicht}Let $P=({\mathcal{A}},{\mathcal{C}})$ be a measuring left $%
R $-pairing. Then the following are equivalent:

(i) $_{R}{\mathcal{C}}$ is locally projective and $\kappa _{P}({\mathcal{A}}%
)\subseteq $ $^{\ast }\mathcal{C}$ is dense;

(ii) $_{R}{\mathcal{C}}$ satisfies the $\alpha $-condition and $\kappa _{P}({%
\mathcal{A}})\subseteq $ $^{\ast }\mathcal{C}$ is dense;

(iii) $({\mathcal{A}},{\mathcal{C}})$ satisfies the left $\alpha $-condition;

(iv) $\sigma \lbrack {\mathcal{C}}_{^{\ast }\mathcal{C}}]\simeq \mathcal{M}^{%
{\mathcal{C}}}\simeq \sigma \lbrack {\mathcal{C}}_{\mathcal{A}}].$

If the equivalent conditions \emph{(i)-(iv)}\ are satisfied, then we have
isomorphisms of categories
\begin{equation}
\begin{tabular}{lllllll}
$\mathcal{M}^{\mathcal{C}}$ & $\simeq $ & $\mathrm{Rat}^{\mathcal{C}}(%
\widetilde{\mathcal{M}}_{\mathcal{A}})$ & $=$ & $\mathrm{Rat}^{\mathcal{C}}(%
\mathcal{M}_{\mathcal{A}})$ & $=$ & $\sigma \lbrack {\mathcal{C}}_{\mathcal{A%
}}]$ \\
& $\simeq $ & $\mathrm{Rat}^{\mathcal{C}}(\widetilde{\mathcal{M}}_{^{\ast }%
\mathcal{C}})$ & $=$ & $\mathrm{Rat}^{\mathcal{C}}(\mathcal{M}_{^{\ast }%
\mathcal{C}})$ & $=$ & $\sigma \lbrack {\mathcal{C}}_{^{\ast }\mathcal{C}}].$%
\end{tabular}
\label{sg}
\end{equation}
\end{theorem}

\begin{Beweis}
(i) $\Longleftrightarrow $ (ii)\ By Lemma \ref{loc-p} $_{R}\mathcal{C}$ is
locally projective iff it satisfies the left $\alpha $-condition.

(ii) $\Rightarrow $ (iii)\ follows from Proposition \ref{rp-rp} (2).

(iii)\ $\Rightarrow $ (iv) If $(\mathcal{A},\mathcal{C})$ satisfies the $%
\alpha $-condition, then clearly $_{R}\mathcal{C}$ satisfies the left $%
\alpha $-condition. Moreover, since (by our convention) $\mathcal{A}$ has
unity with $\kappa _{P}(1\mathcal{_{A}})=\varepsilon _{\mathcal{C}},$ we get
by an analogous argument to that in the proof of Theorem \ref{equal} the
isomorphisms of categories $\mathcal{M}^{\mathcal{C}}\simeq \sigma \lbrack
\mathcal{C}_{\mathcal{A}}]=\sigma \lbrack \mathcal{C}_{^{\ast }\mathcal{C}%
}]. $

(iv) $\Rightarrow $ (i) By Assumption $\sigma \lbrack \mathcal{C}_{^{\ast }%
\mathcal{C}}]\simeq \mathcal{M}^{\mathcal{C}}$ and it follows analog to
\cite[3.5]{Wis2002} that $_{R}\mathcal{C}$ is locally projective. Moreover,
for all $c_{1},...,c_{k}\in \mathcal{C},$ the right $\mathcal{A}$-module $%
(c_{1},...,c_{k})\mathcal{A}\subset \mathcal{C}^{k}$ is a right $^{\ast }%
\mathcal{C}$-submodule (because $\sigma \lbrack \mathcal{C}_{\mathcal{A}%
}]\simeq \sigma \lbrack \mathcal{C}_{^{\ast }\mathcal{C}}]$). Hence $%
(c_{1},...,c_{k})^{\ast }\mathcal{C}=((c_{1},...,c_{k})\mathcal{A})^{\ast }%
\mathcal{C}\subseteq (c_{1},...,c_{k})\mathcal{A},$ i.e. $\kappa _{P}(%
\mathcal{A})\subseteq $ $^{\ast }\mathcal{C}$ is dense (see \cite[15.8]
{Wis88}).

If the equivalent conditions (i)-(iv)\ are\ satisfied, then the assumptions
of Theorem \ref{equal} are satisfied and so the isomorphisms of categories (%
\ref{sg}) are evident.$\blacksquare $
\end{Beweis}

\begin{remark}
Note that if $\mathcal{A}$ has no unity, then the implication (iv) $%
\Rightarrow $ (i) in the previous theorem is still evident, if $^{\ast }%
\mathcal{C}$ is unital as a left $\mathcal{A}$-module. In this case the four
statements become equivalent, if we add ``$\kappa _{P}(\mathcal{A})\subset $
$^{\ast }\mathcal{C}$ is dense'' to statement (iii).
\end{remark}

A dual version of Theorem \ref{cor-dicht} is valid for \emph{measuring right
}$R$\emph{-pairings:}

\begin{theorem}
\label{cd-r}For a measuring right $R$-pairing $P=({\mathcal{A}},{\mathcal{C}}%
)$ the following are equivalent:

(i) ${\mathcal{C}}_{R}$ is locally projective and $\kappa _{P}({\mathcal{A}}%
)\subseteq $ $\mathcal{C}^{\ast }$ is dense;

(ii) ${\mathcal{C}}_{R}$ satisfies the $\alpha $-condition and $\kappa _{P}({%
\mathcal{A}})\subseteq $ $\mathcal{C}^{\ast }$ is dense;

(iii) $({\mathcal{A}},{\mathcal{C}})$ is a right $\alpha $-pairing;

(iv) $\sigma \lbrack _{\mathcal{C}^{\ast }}{\mathcal{C}}]\simeq $ $^{{%
\mathcal{C}}}\mathcal{M}\simeq \sigma \lbrack _{\mathcal{A}}{\mathcal{C}}].$

If the equivalent conditions \emph{(i)-(iv)}\ are satisfied, then we have
isomorphisms of categories
\begin{equation}
\begin{tabular}{lllllll}
$^{\mathcal{C}}\mathcal{M}$ & $\simeq $ & $^{\mathcal{C}}\mathrm{Rat}(_{%
\mathcal{A}}\widetilde{\mathcal{M}})$ & $=$ & $^{\mathcal{C}}\mathrm{Rat}(_{%
\mathcal{A}}\mathcal{M})$ & $=$ & $\sigma \lbrack _{\mathcal{A}}{\mathcal{C}}%
]$ \\
& $\simeq $ & $^{\mathcal{C}}\mathrm{Rat}(_{\mathcal{C}^{\ast }}\widetilde{%
\mathcal{M}})$ & $=$ & $^{\mathcal{C}}\mathrm{Rat}(_{\mathcal{C}^{\ast }}%
\mathcal{M})$ & $=$ & $\sigma \lbrack _{\mathcal{C}^{\ast }}{\mathcal{C}}].$%
\end{tabular}
\label{iso-right}
\end{equation}
\end{theorem}

\begin{remark}
We should mention here that the implications (iii) $\Rightarrow $ (iv)\ in
Theorem \ref{cor-dicht} resp. \ref{cd-r} were achieved independently in \cite
{EG-TL} Theorem 2.6' resp. Theorem 2.6 (note the interchange between \emph{%
left} and \emph{right }pairings in our notation).
\end{remark}

\newpage

As a corollary of Proposition \ref{rp-rp} and Theorem \ref{cor-dicht} we get

\begin{corollary}
\label{conseq}Let $Q=(\mathcal{B},\mathcal{C})\in \mathcal{P}_{ml},$ $\xi :%
\mathcal{A}\rightarrow \mathcal{B}$ a morphism of $R$-rings and consider the
induced measuring left pairing $P:=(\mathcal{A},\mathcal{C}).$ Then the
fallowing are equivalent:

(i) $P\in \mathcal{P}_{ml}^{\alpha };$

(ii) $Q\in \mathcal{P}_{ml}^{\alpha }$ and $\xi (\mathcal{A})\subseteq
\mathcal{B}$ is dense w.r.t. $\mathcal{B}[\frak{T}_{_{ls}}(\mathcal{C})];$

(iii) $_{R}\mathcal{C}$ is locally projective and $\kappa _{P}(\mathcal{A}%
)\subseteq $ $^{\ast }\mathcal{C}$ is dense;

(iv) $\sigma \lbrack \mathcal{C}_{^{\ast }\mathcal{C}}]\simeq \mathcal{M}^{%
\mathcal{C}}\simeq \sigma \lbrack \mathcal{C}_{\mathcal{A}}]\simeq \sigma
\lbrack \mathcal{C}_{\mathcal{B}}].$

If these equivalent conditions are satisfied, then we have isomorphisms of
categories:
\begin{equation}
\begin{tabular}{lllll}
$\mathcal{M}^{\mathcal{C}}$ & $\simeq $ & $\mathrm{Rat}^{\mathcal{C}}(%
\mathcal{M}_{\mathcal{A}})$ & $=$ & $\sigma \lbrack \mathcal{C}_{\mathcal{A}%
}]$ \\
& $\simeq $ & $\mathrm{Rat}^{\mathcal{C}}(\mathcal{M}_{^{\ast }\mathcal{C}})$
& $=$ & $\sigma \lbrack \mathcal{C}_{^{\ast }\mathcal{C}}]$ \\
& $\simeq $ & $\mathrm{Rat}^{\mathcal{C}}(\mathcal{M}_{\mathcal{B}})$ & $=$
& $\sigma \lbrack \mathcal{C}_{\mathcal{B}}].$%
\end{tabular}
\label{cat-eq}
\end{equation}
\end{corollary}

\qquad

\begin{definition}
We call a measuring left (right) $\alpha $-pairing $(\mathcal{A},\mathcal{C}%
) $ coproper, if $\mathcal{A}^{rat}:=\mathrm{R}\mathrm{at}^{\mathcal{C}}(%
\mathcal{A}_{\mathcal{A}})$ ($^{rat}\mathcal{A}:=\mathrm{Rat}^{\mathcal{C}%
}(_{\mathcal{A}}\mathcal{A})$) is dense in $\mathcal{A}.$ An $R$-coring with
$_{R}\mathcal{C}$ ($\mathcal{C}_{R}$) locally projective will be called
\emph{left coproper}\textbf{\ }(\emph{right coproper}), if $^{\Box }\mathcal{%
C}:=\mathrm{R}\mathrm{at}^{\mathcal{C}}(^{\ast }\mathcal{C}_{^{\ast }%
\mathcal{C}})\subset $ $^{\ast }\mathcal{C}$ (resp. $\mathcal{C}^{\Box }:=$ $%
^{\mathcal{C}}\mathrm{R}\mathrm{at}(\mathcal{C}_{\mathcal{C}^{\ast }}^{\ast
})\subseteq \mathcal{C}^{\ast }$) is dense. If $_{R}\mathcal{C}$ and $%
\mathcal{C}_{R}$ are locally projective, then we call $\mathcal{C}$ \emph{%
coproper}, if it is left and right coproper.
\end{definition}

\begin{proposition}
\label{coprop}\emph{(\cite{Abu})} Let $P=(\mathcal{A},\mathcal{C})$ be a
coproper left measuring $\alpha $-pairing \emph{(}i.e. $\mathcal{T}:=\mathrm{%
Rat}^{\mathcal{C}}(\mathcal{A}_{\mathcal{A}})\subset \mathcal{A}$ is dense%
\emph{)}.

1. $\mathcal{C}$ is left coproper, i.e. $^{\Box }\mathcal{C}\subset $ $%
^{\ast }\mathcal{C}$ is dense.

2. For every $f\in $ $\mathcal{T},$ there exists some $e\in $ $\mathcal{T},$
such that $fe=f.$

3. For every right $\mathcal{A}$-module $M$ we have $\mathrm{Rat}^{\mathcal{C%
}}(M_{\mathcal{A}})=M\mathcal{T}.$

4. There is an isomorphism of categories $\mathcal{M}_{^{\Box }\mathcal{C}%
}\simeq \mathcal{M}^{\mathcal{C}}\simeq \mathcal{M}_{\mathcal{T}}.$
\end{proposition}

\begin{punto}
\textbf{Birational modules.}{\normalsize \ }Let $P=(\mathcal{A},\mathcal{C})$
be a measuring left $\alpha $-pairing and $Q=(\mathcal{B},\mathcal{D})$ be a
measuring right $\alpha $-pairing ($\mathcal{A},$ $\mathcal{B}$ not
necessarily with unities). For a $(\mathcal{B},\mathcal{A})$-faithful $(%
\mathcal{B},\mathcal{A})$-bimodule $(M,\rho _{M}^{\mathcal{A}},\rho _{M}^{%
\mathcal{B}})$ it's obvious that $^{\mathcal{D}}\mathrm{Rat}(_{\mathcal{B}%
}M) $ is a right $\mathcal{A}$-module, $\mathrm{Rat}^{\mathcal{C}}(M_{%
\mathcal{A}})$ is a left $\mathcal{B}$-module, and
\begin{equation}
\mathrm{Rat}^{\mathcal{C}}((^{{\normalsize \mathcal{D}}}\mathrm{Rat}(_{%
{\normalsize \mathcal{B}}}M))_{\mathcal{A}})=\text{ }^{{\normalsize \mathcal{%
D}}}\mathrm{Rat}(_{{\normalsize \mathcal{B}}}M){\normalsize \cap }\mathrm{Rat%
}^{\mathcal{C}}(M_{\mathcal{A}})=\text{ }^{{\normalsize \mathcal{D}}}\mathrm{%
Rat}(_{\mathcal{B}}(\mathrm{Rat}^{\mathcal{C}}(M_{\mathcal{A}})))
\label{ratC-D}
\end{equation}
is a $(\mathcal{B},\mathcal{A})$-subbimodule of $M,$ which we call the $(%
\mathcal{D},\mathcal{C})$\emph{-birational}\textbf{\ }$(\mathcal{B},\mathcal{%
A})$\emph{-subbimodule}\textbf{\ }of $M.$ If $M=\mathrm{Rat}^{\mathcal{C}%
}((^{\mathcal{D}}\mathrm{Rat}(_{\mathcal{B}}M))_{\mathcal{A}}),$ then we
call $_{\mathcal{B}}M_{\mathcal{A}}$ $(\mathcal{D},\mathcal{C})$\emph{%
-birational}.

With $^{\mathcal{D}}\mathrm{Rat}^{\mathcal{C}}(_{\mathcal{B}}\widetilde{%
\mathcal{M}}_{\mathcal{A}})\subset $ $_{\mathcal{B}}\widetilde{\mathcal{M}}_{%
\mathcal{A}}$ we denote the full subcategory of $(\mathcal{D},\mathcal{C})$%
-birational $(\mathcal{B},\mathcal{A})$-bimodules. The subcategory of unital
$(\mathcal{D},\mathcal{C})$-birational $(\mathcal{B},\mathcal{A})$-bimodules
is denoted with $^{\mathcal{D}}\mathrm{Rat}^{\mathcal{C}}(_{\mathcal{B}}%
\mathcal{M}_{\mathcal{A}}).$
\end{punto}

\qquad As a generalization of the corresponding result for coalgebras over
base fields (e.g. \cite[Theorem 2.3.3]{DNR2001}) resp. over commutative
rings (\cite[Folgerung 2.2.19]{Abu2001}) we get

\begin{theorem}
\label{bicom}Let $P=(\mathcal{A},\mathcal{C})$ be a measuring left $R$%
-pairing and $Q=(\mathcal{B},\mathcal{D})$ be a measuring right $R$-pairing
\emph{(}$\mathcal{A},$ $\mathcal{B}$ not necessarily with unities\emph{)}.
If $\mathcal{C},\mathcal{D}$ are locally projective, $\kappa _{P}(\mathcal{A}%
)\subseteq $ $^{\ast }\mathcal{C}$ and $\kappa _{Q}(\mathcal{B})\subseteq
\mathcal{D}^{\ast }$ are dense, then there are isomorphisms of categories
\begin{equation*}
\begin{tabular}{llllll}
$^{\mathcal{D}}\mathcal{M}^{\mathcal{C}}$ & $\simeq $ & $^{\mathcal{D}}%
\mathrm{Rat}^{\mathcal{C}}(_{\mathcal{B}}\widetilde{\mathcal{M}}_{\mathcal{A}%
})$ & $=$ & $^{\mathcal{D}}\mathrm{Rat}^{\mathcal{C}}(_{\mathcal{B}}\mathcal{%
M}_{\mathcal{A}})$ &  \\
\multicolumn{1}{c}{} & $\simeq $ & $^{\mathcal{D}}\mathrm{Rat}^{\mathcal{C}%
}(_{\mathcal{D}^{\ast }}\widetilde{\mathcal{M}}_{^{\ast }\mathcal{C}})$ & $=$
& $^{\mathcal{D}}\mathrm{Rat}^{\mathcal{C}}(_{\mathcal{D}^{\ast }}\mathcal{M}%
_{^{\ast }\mathcal{C}})$ &
\end{tabular}
\end{equation*}
\end{theorem}

\begin{Beweis}
Let $M$ be an arbitrary $R$-bimodule. In view of the previous results in
this section it's enough to show that $M$ is a counital $(\mathcal{D},%
\mathcal{C})$-bicomodule iff it's a $(\mathcal{D},\mathcal{C})$-birational $(%
\mathcal{B},\mathcal{A})$-bimodule. If $M$ is a counital $(\mathcal{D},%
\mathcal{C})$-bicomodule, then $M$ is by Lemma \ref{co-rat} (1) a $\mathcal{C%
}$-rational right $\mathcal{A}$-module and analogously a $\mathcal{D}$%
-rational left $\mathcal{B}$-module. Moreover $\alpha _{M}^{Q}$ is obviously
$\mathcal{A}$-linear, $\varrho _{M}^{\mathcal{D}}$ is by assumption $%
\mathcal{C}$-colinear, hence $\mathcal{A}$-linear by Lemma \ref{co-rat} (2).
Consequently $\rho _{M}^{\mathcal{B}}=\alpha _{M}^{Q}\circ \varrho _{M}^{%
\mathcal{D}}$ is $\mathcal{A}$-linear, i.e. $M$ is a $(\mathcal{D},\mathcal{C%
})$-birational $(\mathcal{B},\mathcal{A})$-bimodule.

On the other hand, let $M$ be a $(\mathcal{D},\mathcal{C})$-birational $(%
\mathcal{B},\mathcal{A})$-bimodule. By Lemma \ref{rat-co} $M$ is a counital
right $\mathcal{C}$-comodule and analogously a counital left $\mathcal{D}$%
-comodule. Since $M$ is a $(\mathcal{B},\mathcal{A})$-bimodule, $\rho _{M}^{%
\mathcal{B}}$ is $\mathcal{A}$-linear and so we have for all $a\in \mathcal{A%
}$ and $m\in M:$%
\begin{equation*}
\alpha _{M}^{Q}(\varrho _{M}^{\mathcal{D}}(ma))=\rho _{M}^{\mathcal{B}%
}(ma)=\rho _{M}^{\mathcal{B}}(m)a=(\alpha _{M}^{Q}(\varrho _{M}^{\mathcal{D}%
}(m))a=\alpha _{M}^{Q}(\varrho _{M}^{\mathcal{D}}(m)a),
\end{equation*}
hence $\varrho _{M}^{\mathcal{D}}$ is $\mathcal{A}$-linear by the
injectivity of $\alpha _{M}^{Q}.$ By Lemma \ref{rat-co} (2), $\varrho _{M}^{%
\mathcal{D}}$ is $\mathcal{C}$-colinear, i.e. $M$ is a counital $(\mathcal{D}%
,\mathcal{C})$-bicomodule.$\blacksquare $
\end{Beweis}

\qquad As a consequence of Theorems \ref{cor-dicht} and \ref{cd-r} we get

\begin{proposition}
\label{id-coid}

\begin{enumerate}
\item  Let $P=(\mathcal{A},\mathcal{C})\in \mathcal{P}_{ml}.$ If $K\subset
\mathcal{C}$ is a right $\mathcal{C}$-coideal \emph{(}resp. a left $\mathcal{%
C}$-coideal, a $\mathcal{C}$-bicoideal\emph{)}, then $K^{\bot }$ is a left $%
\mathcal{A}$-ideal \emph{(}resp. a right $\mathcal{A}$-ideal, an $\mathcal{A}
$-ideal\emph{)}. If $K$ is a $\mathcal{C}$-coideal, then $K^{\bot }\subset
\mathcal{A}$ is an $R$-subring with unity $1_{\mathcal{A}}.$ If $P\in
\mathcal{P}_{ml}^{\alpha }$ and $I\subset \mathcal{A}$ is a left $\mathcal{A}
$-ideal, then $I^{\bot }\subset \mathcal{C}$ is a right $\mathcal{C}$%
-coideal.

\item  Let $P=(\mathcal{A},\mathcal{C})\in \mathcal{P}_{mr}.$ If $K\subset
\mathcal{C}$ is a left $\mathcal{C}$-coideal \emph{(}resp. a right $\mathcal{%
C}$-coideal, a $\mathcal{C}$-bicoideal\emph{)}, then $K^{\bot }$ is a right $%
\mathcal{A}$-ideal \emph{(}resp. a left $\mathcal{A}$-ideal, an $\mathcal{A}$%
-ideal\emph{)}. If $K$ is a $\mathcal{C}$-coideal, then $K^{\bot }\subset
\mathcal{A}$ is an $R$-subring with unity $1_{\mathcal{A}}.$ If $P\in
\mathcal{P}_{mr}^{\alpha }$ and $I\subset \mathcal{A}$ is a right $\mathcal{A%
}$-ideal, then $I^{\bot }\subset \mathcal{C}$ is a left $\mathcal{C}$%
-coideal.
\end{enumerate}
\end{proposition}

\qquad

\begin{lemma}
\label{P-rs}Let $X$ be set and $XR$ the free $R$-module with basis $X.$ If $%
R_{R}$ is noetherian, then for every right $R$-module $M,$ the following $R$%
-linear map is injective
\begin{equation}
\beta _{M}:M\otimes _{R}R^{X}\rightarrow M^{X},\text{ }m\otimes f\longmapsto
\lbrack x\longmapsto mf(x)].  \label{bet_M}
\end{equation}
Hence $\widetilde{P}:=(XR,R^{X})$ is a left $\alpha $-pairing.
\end{lemma}

\begin{Beweis}
Let $M$ be an arbitrary right $R$-module and write $M$ as a direct limit of
its f.g. $R$-submodules $M={\underrightarrow{lim}}_{\Lambda }M_{\lambda }$ (
\cite[24.7]{Wis88}). For every $\lambda \in \Lambda ,$ $M_{\lambda }$ is
f.p. in $\mathcal{M}_{R}\ $and we have by (\cite[25.4]{Wis88}) the
isomorphisms
\begin{equation*}
\beta _{M_{\lambda }}:M_{\lambda }\otimes _{R}R^{X}\rightarrow M_{\lambda
}^{X},\text{ }m\otimes f\mapsto \lbrack x\mapsto mf(x)]
\end{equation*}
Moreover for each $\lambda \in \Lambda $ the restriction of $\beta _{M}$ on $%
M_{\lambda }$ is equal to $\beta _{M_{\lambda }}$ and so the following map
is injective:
\begin{equation*}
\beta _{M}=\underrightarrow{lim}\beta _{M_{\lambda }}:\underrightarrow{lim}%
M_{\lambda }\otimes _{R}R^{X}\rightarrow \underrightarrow{lim}M_{\lambda
}^{X}\subset M^{X}.
\end{equation*}
Obviously $\widetilde{P}\in \mathcal{P}_{l}^{\alpha }$ iff $\beta _{M}$ is
injective for every $M\in $ $\mathcal{M}_{R}.\blacksquare $
\end{Beweis}

\begin{corollary}
\label{uno}Let $W,W^{\prime }$ be $R$-bimodules, $X\subset $ $^{\ast
}W,X^{\prime }\subset $ $^{\ast }W^{\prime }$ be $R$-subbimodules and
consider the canonical $R$-linear maps
\begin{equation*}
\kappa :X^{\prime }\otimes _{R}X\rightarrow \text{ }^{\ast }(W\otimes
_{R}W^{\prime })\text{ and }\chi :W\otimes _{R}W^{\prime }\rightarrow
(X^{\prime }\otimes _{R}X)^{\ast }.
\end{equation*}
If $R_{R}$ is noetherian, $W_{R}$ is flat and $\mathrm{Ke}(X)_{R}\subset $ $%
W_{R}$ is pure, then
\begin{equation}
\mathrm{Ke}(\kappa (X^{\prime }\otimes _{R}X))\simeq \mathrm{Ke}(X)\otimes
_{R}W^{\prime }+W\otimes _{R}\mathrm{Ke}({X^{\prime }}).  \label{KeX-Y}
\end{equation}
\newline
\end{corollary}

\begin{Beweis}
Consider the embeddings $E:=W/\mathrm{Ke}(X)\hookrightarrow X^{\ast },$ $%
E^{\prime }:=W^{\prime }/\mathrm{Ke}({X^{\prime })}\hookrightarrow
R^{X^{\prime }}$ and the commutative diagram
\begin{equation*}
\xymatrix{ W \otimes_{R} W' \ar[rrr]^{\chi } \ar[dd]_{\pi \otimes \pi '} & &
& (X' \otimes_{R} X)^* \ar@{^{(}->}[dd]^{\iota} \\ & W/{\rm
{Ke}(X)}\otimes_{R} {R}^{X'} \ar@{^{(}->}[r] & X^* \otimes_{R} {R}^{X'}
\ar@{^{(}->}[rd]^{\beta_{X^*}} & \\ W/{\rm {Ke}(X)}\otimes_{R}
W'/{\rm{Ke}(X')} \ar@{^{(}->}[ur] \ar@{.>}[rrr]_(.6){\delta} & & &
({X^*})^{X'} }
\end{equation*}
It follows by assumptions that $W/\mathrm{Ke}(X)$ is flat in $\mathcal{M}%
_{R} $ and $_{R}R^{X^{\prime }}$ is flat (e.g. \cite[36.5, 26.6]{Wis88}).
Moreover $\beta _{X^{\ast }}$ is injective by Lemma \ref{P-rs}, hence $%
\delta $ is injective. It follows then by \cite[II-3.6]{Bou74} that
\begin{equation*}
\begin{tabular}{lllll}
\textrm{Ke}$(\kappa (X\otimes _{R}X^{\prime }))$ & $:=$ & $\mathrm{Ke}(\chi
) $ & $=$ & \textrm{Ke}$(\delta \circ (\pi \otimes \pi ^{\prime }))$ \\
& $=$ & $\mathrm{Ke}(\pi _{X}\otimes \pi _{X^{\prime }})$ & $=$ & \textrm{Ke}%
$(X)\otimes _{R}W^{\prime }+W\otimes _{R}\mathrm{Ke}({X^{\prime }}%
).\blacksquare $%
\end{tabular}
\end{equation*}
\end{Beweis}

\begin{proposition}
\label{An(B)-coid}Let $R$ be a QF ring and $\mathcal{C}$ an $R$-coring. If $%
\mathcal{A}\subseteq $ $^{\ast }\mathcal{C}$ is an $R$-subring \emph{(}with $%
\varepsilon _{\mathcal{C}}\in \mathcal{A}$\emph{)}, $\mathcal{C}_{R}$ is
flat and $\mathrm{Ke}(\mathcal{A})_{R}\subset \mathcal{C}_{R}$ is pure, then
$\Delta _{\mathcal{C}}(\mathrm{Ke}(\mathcal{A}))\subseteq \mathrm{Ke}(%
\mathcal{A})\otimes _{R}\mathcal{C}+\mathcal{C}\otimes _{R}\mathrm{Ke}(%
\mathcal{A})$ \emph{(}$\mathrm{Ke}(\mathcal{A})\subset \mathcal{C}$ is a $%
\mathcal{C}$-coideal\emph{)}.
\end{proposition}

\begin{Beweis}
Let $\mathcal{A}\subseteq $ $^{\ast }\mathcal{C}$ be an $R$-subring and
consider the $R$-linear map
\begin{equation*}
\kappa :\mathcal{A}\otimes _{R}\mathcal{A}\rightarrow \text{ }^{\ast }(%
\mathcal{C}\otimes _{R}\mathcal{C}),\text{ }a\otimes b\mapsto \lbrack
c\otimes d\mapsto <b,c<a,d>>].
\end{equation*}
If $\mathcal{C}_{R}$ is flat and $\mathrm{Ke}(\mathcal{A})_{R}\subset $ $%
\mathcal{C}_{R}$ is pure, then we have by Corollary \ref{uno} and Lemma \ref
{f*-clos}:
\begin{equation}
\mathrm{Ke}(\mathcal{A})\subseteq \mathrm{Ke}(\Delta _{\mathcal{C}}^{\ast
}(\kappa (\mathcal{A}\otimes _{R}\mathcal{A})))=\Delta _{\mathcal{C}}^{-1}(%
\mathrm{Ke}(\mathcal{A})\otimes _{R}\mathcal{C}+\mathcal{C}\otimes _{R}%
\mathrm{Ke}(\mathcal{A})),  \label{Bij-cood-al}
\end{equation}
i.e. $\Delta _{\mathcal{C}}(\mathrm{Ke}(\mathcal{A}))\subseteq \mathrm{K}%
\mathrm{e}(\mathcal{A})\otimes _{R}\mathcal{C}+\mathcal{C}\otimes _{R}%
\mathrm{Ke}(\mathcal{A}).$ If $\varepsilon _{\mathcal{C}}\in \mathcal{A},$
then $\varepsilon _{\mathcal{C}}(\mathrm{Ke}(\mathcal{A}))=0,$ hence $%
\mathrm{Ke}(\mathcal{A})\subset \mathcal{C}$ is a $\mathcal{C}$-coideal.$%
\blacksquare $
\end{Beweis}

\begin{corollary}
\label{CinMD}Let $\mathcal{C}$ be an $R$-coring and assume that $_{R}%
\mathcal{C}$ is locally projective. For every $R$-coring $\mathcal{D}$ with
an injective morphism of $R$-corings $\iota _{\mathcal{D}}:\mathcal{D}%
\hookrightarrow \mathcal{C}$ we have:

\begin{enumerate}
\item  $_{R}\mathcal{D}\subseteq $ $_{R}\mathcal{C}$ is pure iff $Q:=(^{\ast
}\mathcal{C},\mathcal{D})\in \mathcal{P}_{ml}^{\alpha }.$ In this case $_{R}%
\mathcal{D}$ is locally projective, $\iota _{\mathcal{D}}^{\ast }(^{\ast }%
\mathcal{C})\subseteq $ $^{\ast }\mathcal{D}$ is dense and there are
isomorphisms of categories
\begin{equation*}
\begin{tabular}{lllll}
$\mathcal{M}^{\mathcal{D}}$ & $\simeq $ & $\mathrm{Rat}^{\mathcal{D}}(%
\mathcal{M}_{^{\ast }\mathcal{D}})$ & $=$ & $\sigma \lbrack \mathcal{D}%
_{^{\ast }\mathcal{D}}]$ \\
& $\simeq $ & $\mathrm{Rat}^{\mathcal{D}}(\mathcal{M}_{^{\ast }\mathcal{C}})$
& $=$ & $\sigma \lbrack \mathcal{D}_{^{\ast }\mathcal{C}}].$%
\end{tabular}
\end{equation*}

\item  Let $_{R}R$ be $\mathcal{C}$-injective. Then $_{R}\mathcal{D}%
\subseteq $ $_{R}\mathcal{C}$ is pure iff $_{R}\mathcal{D}$ is locally
projective.

\item  If $_{R}\mathcal{D}\subseteq $ $_{R}\mathcal{C}$ is pure, then
\begin{equation*}
\mathcal{D}=\mathcal{C}\Longleftrightarrow \mathcal{M}{\normalsize ^{%
\mathcal{D}}}=\mathcal{M}{\normalsize ^{\mathcal{C}}}\Longleftrightarrow
\mathcal{C}_{^{\ast }\mathcal{C}}\text{ is }\mathcal{D}\text{-rational.}
\end{equation*}
\end{enumerate}
\end{corollary}

\begin{Beweis}
\begin{enumerate}
\item  Since $\iota _{\mathcal{D}}$ is a morphism of $R$-corings, it follows
that $\iota _{\mathcal{D}}^{\ast }:$ $^{\ast }\mathcal{C}\rightarrow $ $%
^{\ast }\mathcal{D}$ is a morphism of $R$-rings, i.e. $Q$ is a measuring
left $R$-pairing. The result follows now by Theorem \ref{cor-dicht} and the
commutativity of the following diagram for every right $R$-module $M$
\begin{equation}
\xymatrix{ M \otimes_{R} {\mathcal D} \ar@{.>}[drr]^{\alpha_M ^{\mathcal Q}}
\ar[rr]^(.45){id_M \otimes \iota_{\mathcal D}} \ar[d]_{\alpha_M ^{\mathcal
D}} & & M \otimes_{R} {\mathcal C} \ar[d]^{\alpha_M ^{\mathcal C}} \\ {\rm
Hom} _{-{R}} (^{*} {\mathcal D},M) \ar[rr]_{(\iota_{\mathcal D} ^*,M)} & &
{\rm Hom} _{-{R}} (^{*}{\mathcal C},M) }  \label{subco}
\end{equation}

\item  If $_{R}R$ is $\mathcal{C}$-injective, then $\iota _{\mathcal{D}%
}^{\ast }:$ $^{\ast }\mathcal{C}\rightarrow $ $^{\ast }\mathcal{D}$ is
surjective. Hence, for every right $R$-module $M$ the map $(\iota _{\mathcal{%
D}}^{\ast },M)$ in diagram (\ref{subco}) is injective and the result follows.

\item  It is enough to prove: $\mathcal{C}\in \mathrm{Rat}^{\mathcal{D}}(%
\mathcal{M}_{^{\ast }\mathcal{C}})\Longrightarrow \mathcal{D}=\mathcal{C}.$

Assume $\mathcal{C}\in \mathrm{Rat}^{\mathcal{D}}(\mathcal{M}_{^{\ast }%
\mathcal{C}})\simeq \mathcal{M}^{\mathcal{D}}.$ Then there exists a right $R$%
-linear map
\begin{equation*}
\varrho :\mathcal{C}\rightarrow \mathcal{C}\otimes _{R}\mathcal{D},\text{ }%
c\longmapsto \sum\limits_{i=1}^{k_{c}}c_{i}\otimes d_{i},
\end{equation*}
such that $c\leftharpoonup f=\sum\limits_{i=1}^{k_{c}}c_{i}\iota _{\mathcal{D%
}}^{\ast }(f)(d_{i})$ for every $f\in $ $^{\ast }\mathcal{C}.$ Consider the
following diagram
\begin{equation}
\xymatrix{ {\mathcal C} \ar[rr]^(.4){\varrho} \ar[drr]_{\Delta_{\mathcal C}}
& & {\mathcal C} \otimes_{R} {\mathcal D} \ar@{^{(}->}[d]^{id_{\mathcal C}
\otimes \iota_{\mathcal D}} \\ & & {\mathcal C} \otimes_{R} {\mathcal C}}
\label{CincomD}
\end{equation}
and the left $\alpha $-pairing $P\otimes _{l}P:=(^{\ast }\mathcal{C}\otimes
_{R}$ $^{\ast }\mathcal{C},\mathcal{C}\otimes _{R}\mathcal{C})$ (see Lemma
\ref{p-2} (1)). Then we have for all $c\in \mathcal{C}$ and $f,g\in $ $%
^{\ast }\mathcal{C}:$%
\begin{equation*}
\begin{tabular}{lllllll}
$\chi _{P\otimes _{l}P}(\sum c_{1}\otimes c_{2})(f\otimes g)$ & $=$ & $\sum
g(c_{1}f(c_{2}))$ & $=$ & $g(c\leftharpoonup f)$ &  &  \\
& $=$ & $g(\sum\limits_{i=1}^{k_{c}}c_{i}\iota _{\mathcal{D}}^{\ast
}(f)(d_{i}))$ & $=$ & $\chi _{P\otimes
_{l}P}(\sum\limits_{i=1}^{k_{c}}c_{i}\otimes \iota (d_{i}))(f\otimes g),$ &
&
\end{tabular}
\end{equation*}
and so $\sum c_{1}\otimes c_{2}=\sum\limits_{i=1}^{k_{c}}c_{i}\otimes \iota
(d_{i}),$ i.e. diagram (\ref{CincomD}) is commutative. Hence for every $c\in
\mathcal{C}$ we have
\begin{equation*}
c=\sum \varepsilon _{\mathcal{C}}(c_{1})c_{2}=\sum\limits_{i=1}^{k_{c}}%
\varepsilon _{\mathcal{C}}(c_{i})\iota (d_{i})\in \iota (\mathcal{D}),
\end{equation*}
i.e. $\mathcal{C}=\mathcal{D}.\blacksquare $
\end{enumerate}
\end{Beweis}

\begin{remark}
Even if $\mathcal{C}$ is an $R$-coring and $\mathcal{D}\subset \mathcal{C}$
is an $R$-subbimodule with $\Delta _{\mathcal{C}}(\mathcal{D})\subseteq
\func{Im}(\mathcal{D}\otimes _{R}\mathcal{D}),$ $\mathcal{D}$ may have no $R$%
-coring structure such that the natural embedding $\iota _{\mathcal{D}}:%
\mathcal{D}\hookrightarrow \mathcal{C}$ is a morphism of $R$-corings. For $%
\mathcal{D}$ to be an $R$-subcoring of $\mathcal{C}$ we need $_{R}\mathcal{D}%
_{R}\subset $ $_{R}\mathcal{C}_{R}$ to be \emph{pure} (in the sense of
Cohn). A counterexample for coalgebras over commutative rings can be found
in \cite[Page 56]{NS82}.
\end{remark}

An important role by studying the category of rational representations of a
left measuring pairing $P\in \mathcal{P}_{ml}^{\alpha }$ is played by the

\begin{punto}
\textbf{Finiteness Theorem}\label{es}. Let $P=(\mathcal{A},\mathcal{C})\in
\mathcal{P}_{ml}^{\alpha }.$

\begin{enumerate}
\item  If $M\in \mathrm{Rat}^{\mathcal{C}}(\mathcal{M}_{\mathcal{A}})$, then
there exists for every finite subset $\{m_{1},...,m_{k}\}\subset M$ some $%
N\in \mathrm{Rat}^{\mathcal{C}}(\mathcal{M}_{\mathcal{A}})$, such that $%
N\subset M$ and $N_{R}$ is finitely generated.

\item  Every finite subset of $\mathcal{C}$ is contained in a right $%
\mathcal{C}$-coideal, which is f.g. in $\mathcal{M}_{R}.$
\end{enumerate}
\end{punto}

\begin{Beweis}
\begin{enumerate}
\item  Let $M\in \mathrm{Rat}^{\mathcal{C}}(\mathcal{M}_{\mathcal{A}})$ and $%
\{m_{1},...,m_{k}\}\subset M.$ Then $m_{i}\mathcal{A}\subset M$ is an $%
\mathcal{A}$-submodule, hence a $\mathcal{C}$-subcomodule. Moreover $%
m_{i}\in m_{i}\mathcal{A}$ and consequently there exists a subset $%
\{(m_{ij},c_{ij})\}_{j=1}^{n_{i}}\subset m_{i}\mathcal{A}\times \mathcal{C},$
so that $\varrho _{M}(m_{i})=\sum\limits_{j=1}^{n_{i}}m_{ij}\otimes c_{ij}$
for $i=1,...,k.$ Obviously $N:=\sum\limits_{i=1}^{k}m_{i}\mathcal{A}%
=\sum\limits_{i=1}^{k}\sum\limits_{j=1}^{n_{i}}m_{ij}R\subset M$ is a $%
\mathcal{C}$-subcomodule and contains $\{m_{1},...,m_{k}\}.$

\item  This is a special case of (1).$\blacksquare $
\end{enumerate}
\end{Beweis}

\qquad For every $(\mathcal{A},\mathcal{C})\in \mathcal{P}_{ml}^{\alpha }$
we get from the isomorphism of categories $\mathrm{Rat}^{\mathcal{C}}(%
\mathcal{M}_{\mathcal{A}})\simeq \sigma \lbrack \mathcal{C}_{\mathcal{A}}]$
and \cite[2.9]{Wis99}:

\begin{corollary}
\label{sg-loc}Let $(\mathcal{A},\mathcal{C})\in \mathcal{P}_{ml}^{\alpha }.$

\begin{enumerate}
\item  $\mathrm{Rat}^{\mathcal{C}}(\mathcal{M}_{\mathcal{A}})$ is $(\mathcal{%
A},R)$-finite \emph{(}i.e. a $\mathcal{C}$-rational right $\mathcal{A}$%
-module is f.g. in $\mathcal{M}_{\mathcal{A}}$ iff it's f.g. in $\mathcal{M}%
_{R}$\emph{)}.

\item  If $_{R}R$ is perfect, then every $\mathcal{C}$-rational right $%
\mathcal{A}$-module satisfies the descending chain condition w.r.t. the f.g.
$\mathcal{A}$-submodules.

\item  If $R_{R}$ is noetherian, then every $\mathcal{C}$-rational right $%
\mathcal{A}$-module is locally noetherian.

\item  If $R_{R}$ is artinian, then every f.g. $\mathcal{C}$-rational right $%
\mathcal{A}$-module has finite length.
\end{enumerate}
\end{corollary}

\begin{proposition}
\label{MC=AM}For every dense measuring left $R$-pairing $P=(\mathcal{A},%
\mathcal{C})$ the following are equivalent:

\begin{enumerate}
\item  $\mathcal{M}^{\mathcal{C}}\simeq \sigma \lbrack \mathcal{C}_{\mathcal{%
A}}]=\mathcal{M}_{\mathcal{A}/\mathrm{An}_{\mathcal{A}}(\mathcal{C)}};$

\item  The functor $-\otimes _{R}\mathcal{C}:\mathcal{M}_{R}\rightarrow
\mathcal{M}_{\mathcal{A}/\mathrm{An}_{\mathcal{A}}(\mathcal{C)}}$ has a left
adjoint;

\item  $_{R}\mathcal{C}$ is f.g. and projective;

\item  $_{\mathcal{C}^{\ast }}\mathcal{C}$ is f.g. and $_{R}\mathcal{C}$ is
locally projective;

\item  $\mathcal{A}/\mathrm{An}_{\mathcal{A}}(\mathcal{C})$ is f.g. in $%
\mathcal{M}_{R}$ and $_{R}\mathcal{C}$ is locally projective.
\end{enumerate}
\end{proposition}

\begin{Beweis}
With the help of Proposition \ref{propert} and Theorem \ref{cor-dicht}, the
equivalence of the first four statements can be established as in \cite[3.6]
{Wis2002}.

If we assume (1) or (5), then we conclude that $\mathcal{M}^{\mathcal{C}%
}\simeq \sigma \lbrack \mathcal{C}_{\mathcal{A}}]\simeq \sigma \lbrack
\mathcal{C}_{^{\ast }\mathcal{C}}]$ by Proposition \ref{equal} and is $(%
\mathcal{A},R)$-finite by Corollary \ref{sg-loc}. The result follows then by
the fact that in this case $\sigma \lbrack \mathcal{C}_{\mathcal{A}}]=%
\mathcal{M}_{\mathcal{A}/\mathrm{An}_{\mathcal{A}}(\mathcal{C})}$ iff $%
\mathcal{A}/\mathrm{An}_{\mathcal{A}}(\mathcal{C})$ is f.g. in $\mathcal{M}%
_{R}$ \cite[2.9 (3)]{Wis99}.$\blacksquare $
\end{Beweis}

\begin{remark}
It follows from Proposition \ref{MC=AM} that for an $R$-coring $\mathcal{C}$
with $_{R}\mathcal{C}$ locally projective and every dense $R$-subring $%
\mathcal{A}\subseteq $ $^{\ast }\mathcal{C}$ $:$
\begin{equation*}
\mathcal{A}_{R}\text{ is f.g.}\Longleftrightarrow \text{ }_{R}\mathcal{C}%
\text{ is f.g.}\Longleftrightarrow \text{ }_{R}\mathcal{C}\text{ is f.g. and
projective.}
\end{equation*}
In particular a f.g. locally projective $R$-coring is projective.
\end{remark}

\qquad The following result gives topological characterizations of the $%
\mathcal{C}$-rational right $\mathcal{A}$-modules. Here we generalize some
of those characterizations given by D. Radford in \cite[2.2]{Rad73} in the
case of coalgebras over base fields to the case of corings over arbitrary
(artinian) ground rings. See also \cite[Proposition 2.2.26]{Abu2001} for the
case of coalgebras over commutative rings.

\begin{proposition}
\label{n_rat}Let $P=(\mathcal{A},\mathcal{C})\in \mathcal{P}_{ml}^{\alpha }$
and consider $\mathcal{A}$ with the $\mathcal{C}$-adic topology $\mathcal{T}%
_{-\mathcal{C}}(\mathcal{A})=\mathcal{A}[\frak{T}_{ls}^{r}(\mathcal{C})].$
If $M$ is a unital right $\mathcal{A}$-module, then for every $m\in M$ the
following are equivalent:

\begin{enumerate}
\item  there exists a finite subset $F=\{c_{1},...,c_{k}\}\subset \mathcal{C}
$, such that $(0_{\mathcal{C}}:F)\subseteq (0_{M}:m);$

\item  $m\mathcal{A}$ is $\mathcal{C}$-subgenerated;

\item  $m\in \mathrm{Rat}^{\mathcal{C}}(M\mathcal{_{A}});$

\item  there exists a f.g. left $R$-submodule $K\subset \mathcal{C},$ such
that $K^{\bot }\subseteq (0_{M}:m).$

If $R_{R}$ is artinian, then \emph{(1)-(4)} are moreover equivalent to:

\item  $(0_{M}:m)$ contains an $R$-cofinite closed $R$-submodule;

\item  $(0_{M}:m)$ is an $R$-cofinite closed right $\mathcal{A}$-ideal.
\end{enumerate}
\end{proposition}

\begin{Beweis}
(1) $\Longrightarrow $ (2) By assumption and \ref{C-ad} $m\in N:=\mathrm{Sp}%
(\sigma \lbrack \mathcal{C}_{\mathcal{A}}],M).$ Moreover $m\mathcal{A}%
\subset $ $N$ is a right $\mathcal{A}$-submodule and is consequently $%
\mathcal{C}$-subgenerated.

(2) $\Longrightarrow $ (3) By assumption and Theorem \ref{cor-dicht} $m\in m%
\mathcal{A}\subset \mathrm{Rat}^{\mathcal{C}}(M\mathcal{_{A}}).$

(3)\ $\Longrightarrow $ (4) Let $\varrho
(m)=\sum\limits_{i=1}^{k}m_{i}\otimes c_{i}$ and put $K:=\sum%
\limits_{i=1}^{k}Rc_{i}\subset \mathcal{C}.$ Then clearly $K^{\bot
}\subseteq (0_{M}:m).$

(4) $\Longrightarrow $ (1) For every left $R$-submodule $K\subseteq \mathcal{%
C}$ we have $(0_{\mathcal{C}}:K)\subseteq K^{\bot }.$

Let $R_{R}$ be \emph{artinian}.

(4) $\Longrightarrow $ (5). Assume $K^{\bot }\subseteq (0_{M}:m)$ for some $%
K=\sum\limits_{i=1}^{k}Rc_{i}\subset \mathcal{C}.$ Since $\mathcal{A}%
/K^{\bot }\hookrightarrow $ $^{\ast }K$ and $R_{R}$ is noetherian, we
conclude that $K^{\bot }\subset \mathcal{A}$ is $R$-cofinite. Moreover $%
K^{\bot }$ is by Lemma \ref{orth-clos} (1) closed.

The implications (5) $\Rightarrow $ (6) $\Rightarrow $ (1) follow from Lemma
\ref{orth-clos} (3).$\blacksquare $
\end{Beweis}

\section{Applications}

In what follows $R$ is a \emph{commutative }ring and $\mathcal{M}_{R}$ is
the category of $R$-(bi)modules. For an $R$-algebra $(A,\mu _{A},\eta _{A})$
and an $R$-coalgebra $(C,\Delta _{C},\varepsilon _{C})$ we consider $(%
\mathrm{Hom}_{R}(C,A),\star ):=\mathrm{Hom}_{R}(C,A)$ as an $R$-algebra with
the so called \emph{convolution product }$(f\star g)(c):=\sum
f(c_{1})g(c_{2})$ and unity $\eta _{A}\circ \varepsilon _{C}.$ With this
definition $C$ becomes a $C^{\ast }$-bimodule through the left and the right
$C^{\ast }$-action $f\rightharpoonup c=\sum c_{1}f(c_{2})$ and $%
c\leftharpoonup f=\sum f(c_{1})c_{2}.$

\subsection*{Entwined Modules}

Next we apply our results in the previous sections to the category of
entwined modules corresponding to a right-right entwining structure $%
(A,C,\psi ).$ These were introduced by T. Brzezi\'{n}ski and S. Majid \cite
{BM98} as a generalization of the Doi-Koppinen modules corresponding to a
right-right Doi-Koppinen structure (see \ref{DK}).

\begin{punto}
\label{ent-str}\ A \emph{right-right entwining structure}\textbf{\ }$%
(A,C,\psi )$ over $R$ consists of an $R$-algebra $A,$ an $R$-coalgebra $C$
and an $R$-linear map
\begin{equation*}
\psi :C\otimes _{R}A\rightarrow A\otimes _{R}C,\text{ }c\otimes a\mapsto
\sum a_{\psi }\otimes c^{\psi },
\end{equation*}
such that
\begin{equation}
\begin{tabular}{llllll}
$\sum (a\widetilde{a})_{\psi }\otimes c^{\psi }$ & $=$ & $\sum a_{\psi }%
\widetilde{a}_{\Psi }\otimes c^{\psi \Psi },$ & $\sum (1_{A})_{\psi }\otimes
c^{\psi }$ & $=$ & $1_{A}\otimes c,$ \\
$\sum a_{\psi }\otimes \Delta _{C}(c^{\psi })$ & $=$ & $\sum a_{\psi \Psi
}\otimes c_{1}^{\Psi }\otimes c_{2}^{\psi },$ & $\sum a_{\psi }\varepsilon
_{C}(c^{\psi })$ & $=$ & $\varepsilon _{C}(c)a.$%
\end{tabular}
\label{rr-ent}
\end{equation}
\end{punto}

\begin{punto}
Let $(A,C,\psi )\;$be a right-right entwining structure. An \emph{entwined
module}\textbf{\ }corresponding to $(A,C,\psi )$ is a right $A$-module $M,$
which is also a right $C$-comodule through $\varrho _{M},$ such that
\begin{equation*}
\varrho _{M}(ma)=\sum m_{<0>}a_{\psi }\otimes m_{<1>}^{\psi }\text{ for all }%
m\in M\text{ and }a\in A.
\end{equation*}
The category of right-right entwined modules and $A$-linear $C$-colinear
morphisms is denoted by $\mathcal{M}_{A}^{C}(\psi ).$ For $M,N\in \mathcal{M}%
_{A}^{C}(\psi )$ we denote by $\mathrm{Hom}_{A}^{C}(M,N)$ the set of $A$%
-linear $C$-colinear morphisms from $M$ to $N.$ By a remark of M. Takeuchi
(e.g. \cite[Proposition 2.2]{Brz}) $\mathcal{C}:=A\otimes _{R}C$ is an $A$%
-coring with $A$-bimodule structure given by
\begin{equation}
a(\widetilde{a}\otimes c):=a\widetilde{a}\otimes c,\text{ }(\widetilde{a}%
\otimes c)a:=\sum \widetilde{a}a_{\psi }\otimes c^{\psi },  \label{AotC-mod}
\end{equation}
comultiplication
\begin{equation*}
\Delta _{\mathcal{C}}:A\otimes _{R}C\rightarrow (A\otimes _{R}C)\otimes
_{A}(A\otimes _{R}C),\text{ }a\otimes c\mapsto \sum (a\otimes c_{1})\otimes
_{A}(1_{A}\otimes c_{2})
\end{equation*}
and counity $\varepsilon _{\mathcal{C}}:=\vartheta _{A}^{r}\circ
(id_{A}\otimes \varepsilon _{C}).$ Moreover $\mathcal{M}_{A}^{C}(\psi
)\simeq \mathcal{M}^{\mathcal{C}}.$
\end{punto}

\begin{lemma}
\label{lem-entw}\emph{(See \cite[4.2]{Wis2001})} Let $(A,C,\psi )$ be a
right-right entwining structure over $R$ and consider the corresponding $A$%
-coring $\mathcal{C}:=A\otimes _{R}C.$

\begin{enumerate}
\item  $\#_{\psi }^{op}(C,A):=\mathrm{Hom}_{R}(C,A)$ is an $A$-ring with $A$%
-bimodule structure given by $(af)(c):=\sum a_{\psi }f(c^{\psi }),$ $%
(fa)(c):=f(c)a,$ multiplication
\begin{equation}
(f\cdot g)(c)=\sum f(c_{2})_{\psi }g(c_{1}^{\psi }),  \label{re-m}
\end{equation}
and unity $\eta _{A}\circ \varepsilon _{C}.$

\item  $\#_{\psi }^{op}(C,A)\simeq $ $^{\ast }\mathcal{C}$ as $A$-rings via
\begin{equation}
\varphi :\mathrm{Hom}_{R}(C,A)\rightarrow \mathrm{Hom}_{A-}(A\otimes
_{R}C,A),\text{ }f\mapsto \lbrack a\otimes c\mapsto af(c)]  \label{vphi-iso}
\end{equation}
with inverse $h\mapsto \lbrack c\mapsto h(1_{A}\otimes c)].$
\end{enumerate}
\end{lemma}

\begin{punto}
A \emph{left-right entwining structure} is a triple $(A,C,\psi ),$ where $A$
is an $R$-algebra, $C$ is an $R$-coalgebra and
\begin{equation*}
\psi :A\otimes _{R}C\rightarrow A\otimes _{R}C,\text{ }a\otimes c\mapsto
\sum a_{\psi }\otimes c^{\psi },
\end{equation*}
is an $R$-linear map such that the conditions in (\ref{rr-ent}) are
satisfied with the first of them replaced by
\begin{equation*}
\sum (a\widetilde{a})_{\psi }\otimes c^{\psi }=\sum a_{\psi }\widetilde{a}%
_{\Psi }\otimes c^{\Psi \psi }\text{ for all }a,\widetilde{a}\in A,c\in C.
\end{equation*}
\end{punto}

\begin{punto}
Let $(A,C,\psi )$ be a left-right entwining structure. With an entwined
module corresponding to $(A,C,\psi )$ we mean a left $A$-module $M,$ which
is also a right $C$-comodule through $\varrho _{M},$ s.t.
\begin{equation*}
\varrho _{M}(am)=\sum a_{\psi }m_{<0>}\otimes m_{<1>}^{\psi }\text{ for all }%
a\in A\text{ and }m\in M.
\end{equation*}
The category of left-right entwined modules and $A$-linear $C$-colinear
morphisms is denoted by $_{A}\mathcal{M}^{C}(\psi ).$ For $M,N\in $ $_{A}%
\mathcal{M}^{C}(\psi )$ we denote by $_{A}\mathrm{Hom}^{C}(M,N)$ the set of
all $A$-linear $C$-colinear morphisms from $M$ to $N.$ It's easy to see that
$(A^{op},C,\psi \circ \tau )$ is a right-right entwining structure, hence $%
\mathcal{D}:=A^{op}\otimes _{R}C$ is an $A^{op}$-coring and $_{A}\mathcal{M}%
^{C}(\psi )\simeq \mathcal{M}_{A^{op}}^{C}(\psi \circ \tau )\simeq \mathcal{M%
}^{A^{op}\otimes _{R}C}.$ Moreover $\#_{\psi \circ \tau
}^{op}(C,A^{op})\simeq $ $^{\ast }\mathcal{D}$ as $A^{op}$-rings and $%
\#_{\psi }(C,A):=(\#_{\psi \circ \tau }^{op}(C,A^{op}))^{op}$ is an $A$-ring
with multiplication
\begin{equation}
(f\cdot g)(c)=\sum f(c_{1}^{\psi })g(c_{2})_{\psi }\text{ for all }f,g\in
\mathrm{Hom}_{R}(C,A)\text{ and }c\in C  \label{lr-m}
\end{equation}
and unity $\eta _{A}\circ \varepsilon _{C}.$
\end{punto}

\begin{punto}
\label{s-rat}Let $(A,C,\psi )$ be a right-right entwining structure over $R$
and consider the corresponding $A$-coring $\mathcal{C}:=A\otimes _{R}C.$ We
say that $(A,C,\psi )$ satisfies the $\alpha $\emph{-condition}, if for
every right $A$-module $M$ the following map is injective
\begin{equation*}
\alpha _{M}^{\psi }:M\otimes _{R}C\rightarrow \mathrm{Hom}_{R}(\#_{\psi
}^{op}(C,A),M),\text{ }m\otimes c\mapsto \lbrack f\mapsto mf(c)]
\end{equation*}
(equivalently, if $_{A}\mathcal{C}$ is locally projective).
\end{punto}

\qquad Inspired by \cite[3.1]{Doi94} we present

\begin{definition}
\begin{enumerate}
\item  Let $(A,C,\psi )$ be a right-right entwining structure satisfying the
$\alpha $-condition. Let $M\in \mathcal{M}_{\#_{\psi }^{op}(C,A)},$ $\rho
_{M}:M\rightarrow \mathrm{Hom}_{-A}(\#_{\psi }^{op}(C,A),M)$ the canonical
map and $\mathrm{Rat}^{C}(M):=\rho _{M}^{-1}(M\otimes _{R}C).$ If $\mathrm{%
Rat}^{C}(M_{\#_{\psi }^{op}(C,A)})=M,$ then we call $M$ $\#$\emph{-rational}
and set $\varrho _{M}:=(\alpha _{M}^{\psi })^{-1}\circ \rho
_{M}:M\rightarrow M\otimes _{R}C.$ The class of $\#$\emph{-rational} right $%
\#_{\psi }^{op}(C,A)$-modules build a full subcatgery of $\mathcal{M}%
_{\#_{\psi }^{op}(C,A)},$ which we denote with $\mathrm{Rat}^{C}(\mathcal{M}%
_{\#_{\psi }^{op}(C,A)}).$ For a left-right entwining structure $(A,C,\psi
), $ the $\alpha $-condition and the category of $\#$\emph{-rational left}%
\textbf{\ }$\#_{\psi }(C,A)$\emph{-modules} are analogously defined.
\end{enumerate}
\end{definition}

\begin{lemma}
\label{ent-alph}Let $(A,C,\psi )$ be a right-right entwining structure over $%
R$ and consider the corresponding $A$-coring $\mathcal{C}:=A\otimes _{R}C.$

\begin{enumerate}
\item  If $_{R}C$ is flat \emph{(}resp. projective, f.g.\emph{)}, then $_{A}%
\mathcal{C}$ is flat \emph{(}resp. projective, f.g.\emph{)}.

\item  If $_{R}C$ is locally projective, then $_{A}\mathcal{C}$ is locally
projective \emph{(}i.e. $(A,C,\psi )$ satisfies the left $\alpha $-condition%
\emph{)}.
\end{enumerate}
\end{lemma}

\begin{Beweis}
\begin{enumerate}
\item  Clear.

\item  For every right $A$-module $M$ we have the commutative diagram
\begin{equation*}
\xymatrix{ M \otimes_A (A \otimes_R C) \ar[rr]^(.45){\alpha_M ^{\mathcal C}}
\ar@{=}[d] & & {\rm Hom}_{-A} ({}^{\ast }{\mathcal C},M) \ar[d]^{(\gamma
,M)} \\ M \otimes_R C \ar[rr]_{\alpha_M ^C} & & {\rm Hom}_{R} (C^*,M) }
\end{equation*}
where $\gamma :C^{\ast }\rightarrow $ $^{\ast }\mathcal{C},$ $f\mapsto
\lbrack a\otimes c\mapsto af(c)].$ The result follows then by Lemma \ref
{loc-p}.$\blacksquare $
\end{enumerate}
\end{Beweis}

\begin{lemma}
\label{rat-ent}Let $(A,C,\psi )$ be a right-right entwining structure.

\begin{enumerate}
\item  If $M\in \mathcal{M}_{A}^{C}(\psi )$ and $N\subset M$ is a $C$%
-subcomodule, then $NA$ is a subobject of $M$ in $\mathcal{M}_{A}^{C}(\psi
). $

\item  Assume $_{R}C$ to be locally projective.

\begin{enumerate}
\item  For every right $\#_{\psi }^{op}(C,A)$ module $M$ we have $\mathrm{Rat%
}^{C}(M_{\#_{\psi }^{op}(C,A)})\in \mathcal{M}_{A}^{C}(\psi ).$

\item  If $M\in \mathcal{M}_{A}^{C}(\psi ),$ then $M$ becomes a $\#$%
-rational right $\#_{\psi }^{op}(C,A)$-module through
\begin{equation*}
mf=\sum m_{<0>}f(m_{<1>})\text{ for all }m\in M\text{ and }f\in \#_{\psi
}^{op}(C,A).
\end{equation*}

\item  Assume $A\in \mathcal{M}_{A}^{C}(\psi ),$ so that $\sum
1_{<0>}\otimes 1_{<1>}\in \mathcal{C}$ is a group-like element. Let $M\in
\mathcal{M}_{\#_{\psi }^{op}(C,A)}$ be $C$-rational and put $M^{co\mathcal{C}%
}:=\{m\in M\mid $ $\varrho _{M}(m)=\sum m1_{<0>}\otimes 1_{<1>}\}.$ If $\Psi
_{M}:M^{co\mathcal{C}}\otimes _{B}A\rightarrow M,$ $m\otimes a\mapsto ma$ is
surjective, then $M$ is $\#$-rational.
\end{enumerate}
\end{enumerate}
\end{lemma}

\begin{Beweis}
\begin{enumerate}
\item  For every $n\in N$ and $a\in A$ we have
\begin{equation*}
\varrho _{M}(na)=\sum n_{<0>}a_{\psi }\otimes n_{<1>}^{\psi }\in NA\otimes
_{R}C.
\end{equation*}
Consequently $NA\subset M$ is a $C$-subcomodule with structure map $(\varrho
_{M})_{\mid _{NA}},$ hence $NA\subset M$ is a subobject in $\mathcal{M}%
_{A}^{C}(\psi ).$

\item  Assume $_{R}C$ to be locally projective.

\begin{enumerate}
\item  Every right $\#_{\psi }^{op}(C,A)$-module $M$ becomes a right $A$%
-module through the canonical algebra morphism $\iota _{A}:A\rightarrow
\#_{\psi }^{op}(C,A)$ and a left $C^{\ast }$-module through the algebra
anti-morphism $\iota _{C^{\ast }}:C^{\ast }\rightarrow \#_{\psi }^{op}(C,A).$
Clearly $\mathrm{Rat}^{C}(M_{\#_{\psi }^{op}(C,A)})\subset M$ is a $C$%
-rational left $C^{\ast }$-module. Moreover, for all $m\in \mathrm{Rat}%
^{C}(M_{\#_{\psi }^{op}(C,A)}),$ $a\in A$ and $g\in \#_{\psi }^{op}(C,A)$ we
have
\begin{equation*}
\begin{tabular}{lllll}
$\lbrack ma]g$ & $=$ & $(m\iota _{A}(a))g$ &  &  \\
& $=$ & $m(\iota _{A}(a)\cdot g)$ &  &  \\
& $=$ & $\sum m_{<0>}(\iota _{A}(a)\cdot g))(m_{<1>})$ &  &  \\
& $=$ & $\sum m_{<0>}\iota _{A}(a)(m_{<1>2})_{\psi }g(m_{<1>1}^{\psi })$ &
&  \\
& $=$ & $\sum m_{<0>}(a\varepsilon _{C}(m_{<1>2}))_{\psi }g(m_{<1>1}^{\psi
}) $ &  &  \\
& $=$ & $\sum m_{<0>}a_{\psi }g(m_{<1>}^{\psi }),$ &  &
\end{tabular}
\end{equation*}
i.e. $ma\in \mathrm{Rat}^{C}(M_{\#_{\psi }^{op}(C,A)})$ with $\varrho
(ma)=\sum m_{<0>}a_{\psi }\otimes m_{<1>}^{\psi }$ and the result follows.

\item  Let $M\in \mathcal{M}_{A}^{C}(\psi ).$ Then for every $m\in M$ and $%
f,g\in \#_{\psi }^{op}(C,A)\;$we have
\begin{equation*}
\begin{tabular}{lllll}
$m(f\cdot g)$ & $=$ & $\sum m_{<0>}(f\cdot g)(m_{<1>})$ &  &  \\
& $=$ & $\sum m_{<0>}f(m_{<1>2})_{\psi }g(m_{<1>1}^{\psi })$ &  &  \\
& $=$ & $\sum m_{<0><0>}f(m_{<1>})_{\psi }g(m_{<0><1>}^{\psi })$ &  &  \\
& $=$ & $(\sum m_{<0>}f(m_{<1>}))g$ &  &  \\
& $=$ & $(mf)g.$ &  &
\end{tabular}
\end{equation*}

\item  Let $m\in M$ be arbitrary. By assumption $m=\Psi
_{M}(\sum\limits_{i=1}^{k}n_{i}\otimes a_{i})$ for some $\sum%
\limits_{i=1}^{k}n_{i}\otimes a_{i}\in M^{co\mathcal{C}}\otimes _{R}A,$
hence we have for all $f\in \#_{\psi }^{op}(C,A):$%
\begin{equation*}
\begin{tabular}{lllll}
$mf$ & $=$ & $\sum\limits_{i=1}^{k}n_{i}(a_{i}f)$ & $=$ & $%
\sum\limits_{i=1}^{k}n_{i}1_{<0>}a_{i_{_{\psi }}}f(1_{<1>}^{\psi })$ \\
& $=$ & $\sum\limits_{i=1}^{k}(n_{i}a_{i})_{<0>}f((n_{i}a_{i})_{<1>})$ & $=$
& $\sum\limits_{i=1}^{k}m_{<0>}f(m_{<1>}).\blacksquare $%
\end{tabular}
\end{equation*}
\end{enumerate}
\end{enumerate}
\end{Beweis}

\qquad The main result in this section is

\begin{theorem}
\label{ent-sg}Let $(A,C,\psi )$ be a right-right entwining structure and
consider the corresponding $A$-coring $\mathcal{C}:=A\otimes _{R}C.$

\begin{enumerate}
\item  If $_{R}C$ is flat, then $\mathcal{M}_{A}^{C}(\psi )$ is a
Grothendieck category with enough injective objects.

\item  If $_{R}C$ is locally projective \emph{(}resp. f.g. and projective%
\emph{)}, then
\begin{equation}
\mathcal{M}_{A}^{C}(\psi )\simeq \mathrm{Rat}^{C}(\mathcal{M}_{\#_{\psi
}^{op}(C,A)})\simeq \sigma \lbrack (A\otimes _{R}C)_{\#_{\psi }^{op}(C,A)}]%
\text{\ \emph{(}resp. }\mathcal{M}_{A}^{C}(\psi )\simeq \mathcal{M}%
_{\#_{\psi }^{op}(C,A)}\text{\emph{)}.}  \label{iso-sm}
\end{equation}
\end{enumerate}
\end{theorem}

\begin{Beweis}
\begin{enumerate}
\item  If $_{R}C$ is flat, then $_{A}\mathcal{C}$ is flat and the result
follows by the isomorphism $\mathcal{M}_{A}^{C}(\psi )\simeq \mathcal{M}^{%
\mathcal{C}}$ and Proposition \ref{propert} (3) (this is a generalization of
\cite[Section 2.8, Corollary 4]{CMZ2002}, where $(A,C,\psi )$ is \emph{%
monoidal }and $R$ is a base field).

\item  If $_{R}C$ is locally projective, then $_{A}\mathcal{C}$ satisfies
the $\alpha $-condition by Lemma \ref{ent-alph} (2), i.e. $(^{\ast }\mathcal{%
C},\mathcal{C})\in \mathcal{P}_{ml}^{\alpha }.$ The result follows now by
Theorem \ref{cor-dicht} and Lemma \ref{rat-ent}. If $_{R}C$ is f.g. and
projective, then $_{A}\mathcal{C}$ is also f.g. and projective and the
result follows by Proposition \ref{MC=AM}.$\blacksquare $
\end{enumerate}
\end{Beweis}

\begin{corollary}
\label{lr-ent-sg}Let $(A,C,\psi )$ be a left-right entwining structure and
consider the corresponding $A^{op}$-coring $\mathcal{D}:=A^{op}\otimes
_{R}C. $ If $_{R}C$ is flat, then by Theorem \ref{ent-sg} $_{A}\mathcal{M}%
^{C}(\psi )\simeq \mathcal{M}_{A^{op}}^{C}(\psi \circ \tau )$ is a
Grothendieck category with enough injective objects. If moreover $_{R}C$ is
locally projective \emph{(}resp. f.g. and projective\emph{)}, then
\begin{equation*}
_{A}\mathcal{M}^{C}(\psi )\simeq \mathrm{Rat}^{C}(_{\#_{\psi }(C,A)}\mathcal{%
M})=\sigma \lbrack _{\#_{\psi }(C,A)}\mathcal{D}]\;\text{\emph{(}resp. }_{A}%
\mathcal{M}^{C}(\psi )\simeq \text{ }_{\#_{\psi }(C,A)}\mathcal{M}\text{%
\emph{)}.}
\end{equation*}
\end{corollary}

\begin{punto}
\label{func}Let $(A,C,\psi )$ be a right-right entwining structure.

\begin{enumerate}
\item  By \cite[Corollaries 3.4, 3.7]{Brz99} $-\otimes _{R}^{c}A:\mathcal{M}%
^{C}\rightarrow \mathcal{M}_{A}^{C}(\psi )$ is a functor, where for every $%
N\in \mathcal{M}^{C}$ we consider the canonical right $A$-module $N\otimes
_{R}^{c}A:=N\otimes _{R}A$ with the $C$-coaction $n\otimes a\mapsto \lbrack
\sum n_{<0>}\otimes a_{\psi }\otimes n_{<1>}^{\psi }].$ Moreover $-\otimes
_{R}A$ is left adjoint to the forgetful functor $\mathcal{F}_{A}:\mathcal{M}%
_{A}^{C}(\psi )\rightarrow \mathcal{M}^{C},$ where for every $N\in \mathcal{M%
}^{C}$ and $M\in \mathcal{M}_{A}^{C}(\psi )$
\begin{equation}
\mathrm{Hom}_{A}^{C}(N\otimes _{R}^{c}A,M)\rightarrow \mathrm{Hom}^{C}(N,M),%
\text{ }g\mapsto g(-\otimes 1_{A})  \label{ot-A}
\end{equation}
is a functorial isomorphism with inverse $f\mapsto \lbrack n\otimes a\mapsto
f(n)a].$

\item  By the isomorphism $\mathcal{M}_{A}^{C}(\psi )\simeq \mathcal{M}^{%
\mathcal{C}}$ and Proposition \ref{propert} (1) $-\otimes _{R}C\simeq
-\otimes _{A}\mathcal{C}:\mathcal{M}_{A}\rightarrow \mathcal{M}_{A}^{C}(\psi
)$ is a functor, where for every $N\in \mathcal{M}_{A}$ we consider the
canonical right $C$-comodule $N\otimes _{R}C$ with the $A$-action $(n\otimes
c)a\mapsto \sum na_{\psi }\otimes c^{\psi }.$ Moreover $-\otimes _{R}C$ is
right adjoint to the forgetful functor $\mathcal{F}^{C}:\mathcal{M}%
_{A}^{C}(\psi )\rightarrow \mathcal{M}_{A}$ and left adjoint to $\mathrm{Hom}%
_{A}^{C}(\mathcal{C},-):\mathcal{M}_{A}^{C}(\psi )\rightarrow \mathcal{M}%
_{A}.$
\end{enumerate}
\end{punto}

\begin{definition}
Let $C$ be an $R$-coalgebra.

\begin{enumerate}
\item  $C$ is said to be \emph{left}\textbf{\ }(\emph{right})\textbf{\ }%
\emph{Quasi-co-Frobenius}, if $C$ is cogenerated by $C^{\ast }$ as a left (a
right) $C^{\ast }$-module (i.e. $C$ is a torsionless $C^{\ast }$-module \cite
{G-TN95}).

\item  Assume $_{R}C$ to be locally projective. After \cite{AT78} we call $C$
\emph{left coproper} (resp. \emph{right coproper}), if $C^{\Box }:=\mathrm{%
Rat}^{C}(_{C^{\ast }}C^{\ast })$ (resp. $^{\Box }C:=\mathrm{Rat}%
^{C}(C_{C^{\ast }}^{\ast })$) is dense in $C^{\ast }.$ We call $C$ \emph{%
coproper}, if $C$ is left and right coproper.
\end{enumerate}
\end{definition}

\begin{corollary}
\label{cogen}Let $(A,C,\psi )$ be a right-right entwining structure and
consider the corresponding $A$-coring $\mathcal{C}:=A\otimes _{R}C.$

\begin{enumerate}
\item  Let $C$ be projective in $\mathcal{M}^{C}$ \emph{(}e.g. $R$ is a QF
ring and $C$ is left Quasi-co-Frobenius\emph{)}. Then $C\otimes _{R}^{c}A$
is projective in $\mathcal{M}_{A}^{C}(\psi ).$ If moreover $\psi $ is
bijective, then $A\otimes _{R}C$ is also projective in $\mathcal{M}%
_{A}^{C}(\psi ).$

\item  If $_{R}C$ is locally projective and left coproper, then $C^{\Box
}\otimes _{R}^{c}A$ is a generator in $\mathcal{M}_{A}^{C}(\psi ).$

\item  If $A$ is a cogenerator in $\mathcal{M}_{A},$ then $A\otimes _{R}C$
is a cogenerator in $\mathcal{M}_{A}^{C}(\psi ).$ If $_{R}C$ is flat and $%
A_{A}$ is injective, then $A\otimes _{R}C$ is injective in $\mathcal{M}%
_{A}^{C}(\psi ).$
\end{enumerate}
\end{corollary}

\begin{Beweis}
\begin{enumerate}
\item  This follows from the functorial isomorphism (\ref{ot-A}): $\mathrm{%
Hom}_{A}^{C}(C\otimes _{R}^{c}A,M)\simeq \mathrm{Hom}^{C}(C,M)$ for every $%
M\in \mathcal{M}_{A}^{C}(\psi ).$ If $R$ is a QF ring and $C$ is left
Quasi-co-Frobenius, then $C$ is projective in $\mathcal{M}^{C}$ by \cite
{MTW2001}\emph{. }Note that $\psi $ is a morphism in $\mathcal{M}%
_{A}^{C}(\psi ),$ hence $A\otimes _{R}C\simeq C\otimes _{R}^{c}A$ in $%
\mathcal{M}_{A}^{C}(\psi ),$ if $\psi $ is bijective.

\item  If $_{R}C$ is locally projective and left coproper, then $C^{\Box }$
is a generator in $\sigma \lbrack _{C^{\ast }}C]\simeq \mathcal{M}^{C}$ by
\cite[2.6]{Wis99}. The result follows then by the functorial isomorphism (%
\ref{ot-A}): $\mathrm{Hom}_{A}^{C}(C^{\Box }\otimes _{R}^{c}A,M)\simeq
\mathrm{Hom}^{C}(C^{\Box },M)$ for every $M\in \mathcal{M}_{A}^{C}(\psi ).$

\item  Consider the corresponding $A$-coring $\mathcal{C}:=A\otimes _{R}C.$
If $A$ is a cogenerator in $\mathcal{M}_{A},$ then by Lemma \ref{propert}
(2) $\mathcal{C}$ is a cogenerator in $\mathcal{M}^{\mathcal{C}}\simeq
\mathcal{M}_{A}^{C}(\psi ).$ If $_{R}C$ is flat, then $_{A}\mathcal{C}$ is
flat and the second statement follows by Proposition \ref{propert} (7).$%
\blacksquare $
\end{enumerate}
\end{Beweis}

\subsection*{Doi-Koppinen Modules}

In what follows we consider a fundamental class of entwined modules, namely
the class of Doi-Koppinen modules introduced independently by Y. Doi \cite
{Doi92} and M. Koppinen \cite{Kop95}.

\begin{punto}
\label{DK}A \emph{right-right Doi-Koppinen structure} over $R$ is a triple $%
(H,A,C)$ consisting of an $R$-bialgebra $H,$ a right $H$-comodule algebra $A$
and a right $H$-module coalgebra $C.$ A \emph{right-right Doi-Koppinen module%
} for $(H,A,C)$ is a right $A$-module $M,$ which is also a right $C$%
-comodule through $\varrho _{M},$ such that
\begin{equation*}
\varrho _{M}(ma)=\sum m_{<0>}a_{<0>}\otimes m_{<1>}a_{<1>}\text{ for all }%
m\in M\text{ and }a\in A.
\end{equation*}
With $\mathcal{M}(H)_{A}^{C}$ we denote the category of right-right
Doi-Koppinen modules and $A$-linear $C$-colinear morphisms. By \cite[Page
295]{Brz99} $(A,C,\psi )$ is a right-right entwining structure and $\mathcal{%
M}(H)_{A}^{C}\simeq \mathcal{M}_{A}^{C}(\psi ),$ where
\begin{equation}
\psi :C\otimes _{R}A\rightarrow A\otimes _{R}C,\text{ }c\otimes a\mapsto
\sum a_{<0>}\otimes ca_{<1>}.  \label{psi-DK}
\end{equation}
By Lemma \ref{lem-entw} $\#^{op}(C,A):=\mathrm{Hom}_{R}(C,A)$ is an $A$-ring
with $A$-bimodule structure
\begin{equation*}
(af)(c):=\sum a_{<0>}f(ca_{<1>})\text{ and }(fa)(c):=f(c)a,
\end{equation*}
multiplication
\begin{equation}
(f\cdot g)(c)=\sum f(c_{2})_{<0>}g(c_{1}f\left( c_{2}\right) _{<1>})
\label{DK-mult}
\end{equation}
and unity $\eta _{A}\circ \varepsilon _{C}.$ Moreover we have, with $%
\mathcal{C}:=A\otimes _{R}C$ the corresponding $A$-coring, an isomorphism of
$A$-rings $\#^{op}(C,A)\simeq $ $^{\ast }\mathcal{C}.$ The $R$-algebra $%
\#^{op}(C,A)$ was introduced by M. Koppinen \cite[2.2]{Kop95}.
\end{punto}

\begin{punto}
Let $H$ be an $R$-bialgebra. Since $H$ itself is a right $H$-module
coalgebra with structure map $\mu _{H},$ it turns out that, for every right $%
H$-comodule algebra $A,$ the triple $(H,A,H)$ is a right-right Doi-Koppinen
structure and $\mathcal{M}(H)_{A}^{H}=\mathcal{M}_{A}^{H},$ the category of
\emph{relative Hopf modules}\textbf{\ }investigated in \cite{Doi83}.\textbf{%
\ }Note also that $H$ is a right $H$-comodule algebra with structure map $%
\Delta _{H}$ and it turns out that, for every right $H$-module coalgebra $C,$
the triple $(H,H,C)$ is a right-right Doi-Koppinen structure and $\mathcal{M}%
(H)_{H}^{C}=\mathcal{M}_{[C,H]},$ the category of \emph{Doi's }$[C,H]$\emph{%
-modules}\textbf{\ }introduced in \cite{Doi83}. Finally $(H,H,H)$ is a
right-right Doi-Koppinen structure and $\mathcal{M}(H)_{H}^{H}=\mathcal{M}%
_{H}^{H},$ the category of \emph{Hopf modules} studied by M. Sweedler
\cite[4.1]{Swe69}.
\end{punto}

\qquad The following result is easy to prove.

\begin{lemma}
\label{rr-ring}Let $(H,A,C)\;$be a right-right Doi-Koppinen structure over $%
R,$ $\mathcal{C}:=A\otimes _{R}C$ the corresponding $A$-coring and $%
T\subseteq C^{\ast }$ a left $H$-module subalgebra.

\begin{enumerate}
\item  $A\#^{op}T:=A\otimes _{R}T$ is an $A$-ring with $A$-bimodule
structure
\begin{equation}
\widetilde{a}(a\#f):=\sum \widetilde{a}_{<0>}a\#\widetilde{a}_{<1>}f\text{
and }(a\#f)\widetilde{a}:=a\widetilde{a}\#f  \label{AC*-bimod}
\end{equation}
and multiplication
\begin{equation}
(a\#f)\cdot (b\#g):=\sum a_{<0>}b\#(a_{<1>}g)\star f.  \label{op-smash}
\end{equation}
If $\varepsilon _{C}\in T,$ then $1_{A}\#\varepsilon _{C}$ is a unity for $%
A\#^{op}T$ and $A\rightarrow A\#^{op}T,$ $a\mapsto a\#\varepsilon _{C}$ is a
morphism of $A$-rings.

\item  We have a morphism of $A$-rings
\begin{equation*}
\beta :A\#^{op}T\rightarrow \#^{op}(C,A),\text{ }a\#f\mapsto \lbrack
c\mapsto af(c)].
\end{equation*}
Hence $Q:=(A\#^{op}T,\mathcal{C})$ is a measuring left $A$-pairing with
\begin{equation*}
\kappa _{Q}:=\varphi \circ \beta :A\#^{op}T\rightarrow \text{ }^{\ast }%
\mathcal{C},\text{ }a\#f\mapsto \lbrack \widetilde{a}\otimes c\mapsto
\widetilde{a}af(c)].
\end{equation*}
\end{enumerate}
\end{lemma}

\begin{theorem}
\label{HAC-iso}Let $(H,A,C)$ be a right-right Doi-Koppinen structure and
consider the corresponding $A$-coring $\mathcal{C}:=A\otimes _{R}C.$

\begin{enumerate}
\item  If $_{R}C$ is flat, then $\mathcal{M}(H)_{A}^{C}$ is a Grothendieck
catgeory with enough injective objects.

\item  Let $T\subseteq C^{\ast }$ be an $A$-pure left $H$-module subalgebra
and $Q:=(A\#^{op}T,\mathcal{C}).$ If $_{R}C$ is locally projective \emph{(}%
resp. f.g. and projective\emph{) }and $T\subseteq C^{\ast }$ is dense, then $%
\beta (A\#^{op}T)\subseteq \#^{op}(C,A)$ is dense, $Q\in \mathcal{P}%
_{ml}^{\alpha }$ and we have isomorphisms of categories
\begin{equation*}
\mathcal{M}(H)_{A}^{C}\simeq \sigma \lbrack \mathcal{C}_{\#^{op}(C,A)}]=%
\sigma \lbrack \mathcal{C}_{A\#^{op}T}]\;\text{\emph{(}resp. }\mathcal{M}%
(H)_{A}^{C}\simeq \mathcal{M}_{\#^{op}(C,A)}\simeq \mathcal{M}%
_{A\#^{op}C^{\ast }}\text{\emph{)}.}
\end{equation*}
\end{enumerate}
\end{theorem}

\begin{Beweis}
\begin{enumerate}
\item  Since $\mathcal{M}(H)_{A}^{C}\simeq \mathcal{M}_{A}^{C}(\psi ),$
where $\psi $ is defined in (\ref{psi-DK}), the result follows by
Proposition \ref{ent-sg} (1).

\item  Consider the left measuring $A$-pairing $P:=(A\#^{op}C^{\ast },%
\mathcal{C})$ and let $\phi :C^{\ast }\rightarrow A\otimes _{R}C^{\ast },$ $%
f\mapsto 1_{A}\otimes f.$ Then we have for every right $A$-module $M$ the
following commutative diagram
\begin{equation*}
\xymatrix{ M \otimes_A (A \otimes_R C) \ar[rr]^{\alpha_M ^{P}} \ar@{=}[d] &
& {\rm Hom}_{-A} (A \otimes_R C^*,M) \ar[d]^{(\phi,M)} \\ M \otimes_R C
\ar[rr]_{\alpha_M ^{C}} & & {\rm Hom}_{R} (C^*,M) }
\end{equation*}
Let $_{R}C$ be locally projective. Then $\alpha _{M}^{C}$ is injective and
so $\alpha _{M}^{P}$ is injective. Since $M$ is an arbitrary right $A$%
-module, $P$ satisfies the $\alpha $-condition and we get by Theorem \ref
{ent-sg} the category isomorphisms $\mathcal{M}(H)_{A}^{C}\simeq \sigma
\lbrack \mathcal{C}_{\#^{op}(C,A)}]\simeq \sigma \lbrack \mathcal{C}%
_{A\#^{op}C^{\ast }}].$ It follows then by Theorem \ref{cor-dicht} that $%
\kappa _{P}(A\#^{op}C^{\ast })\subseteq $ $^{\ast }\mathcal{C}$ is dense. If
$T\subseteq C^{\ast }$ is an $A$-pure dense left $H$-module subalgebra, then
obviously $A\#^{op}T\subseteq A\#^{op}C^{\ast }$ is dense, hence $\kappa
_{Q}(A\#^{op}T)\subseteq $ $^{\ast }\mathcal{C}$ is dense. Since $^{\ast }%
\mathcal{C}\simeq \#^{op}(C,A)$ it follows then that $\beta
(A\#^{op}T)\subseteq \#^{op}(C,A)$ is dense.

If $_{R}C$ is f.g. and projective, then $\mathcal{M}(H)_{A}^{C}\simeq
\mathcal{M}_{\#^{op}(C,A)}$ by Theorem \ref{ent-sg} (2). Note that in this
case $A\#^{op}C^{\ast }\simeq \#^{op}(C,A)$ and the result follows.$%
\blacksquare $
\end{enumerate}
\end{Beweis}

\begin{punto}
\label{lr-DK}A left-right Doi-Koppinen structure is a triple $(H,A,C),\;$%
where $H$ is an $R$-bialgebra, $A$ is a right $H$-comodule algebra and $C$
is a left $H$-module coalgebra. A Doi-Koppinen module corresponding to $%
(H,A,C)$ is a left $A$-module $M,$ which is also a right $C$-comodule
through $\varrho _{M},$ s.t. $\varrho _{M}(am)=\sum a_{<0>}m_{<0>}\otimes
a_{<1>}m_{<1>}.$ The category of left-right Doi-Koppinen modules and $A$%
-linear $C$-colinear morphisms is denoted by $_{A}\mathcal{M}^{C}(H).$ It
turns out that $(H^{op},A^{op},C)$ is a right-right Doi-Koppinen structure,
hence $\#(C,A):=(\#^{op}(C,A^{op}))^{op}$ is an $A$-ring with multiplication
\begin{equation}
(f\cdot g)(c)=\sum f(g\left( c_{2}\right) _{<1>}c_{1})g(c_{2})_{<0>}.
\label{(C,A)}
\end{equation}
and unity $\eta _{A}\circ \varepsilon _{C}.$ For every right $H$-module
subalgebra $T\subseteq C^{\ast }$ (with $\varepsilon _{C}\in T$) the smash
product $A\#T:=(A^{op}\#^{op}T)^{op}$ is an $A$-ring with multiplication
\begin{equation}
(a\#f)\cdot (b\#g):=\sum ab_{<0>}\#(fb_{<1>})\star g  \label{smash}
\end{equation}
(and unity $1_{A}\#\varepsilon _{C}$). In fact the $R$-algebra $\#(C,A)$
(resp. $A\#T$) was introduced in \cite[2.1]{Kop95} (resp. in \cite[Page 375]
{Doi92}).
\end{punto}

\begin{corollary}
Let $(H,A,C)$ be a left-right Doi-Koppinen structure, $\mathcal{D}%
:=A^{op}\otimes _{R}C$ the corresponding $A^{op}$-coring and $\beta
:A\#C^{\ast }\rightarrow \#(C,A)\;$the canonical morphism. If $_{R}C$ is
flat, then $_{A}\mathcal{M}(H)^{C}\simeq \mathcal{M}(H^{op})_{A^{op}}^{C}$
is a Grothendieck category with enough injective objects. If $T\subseteq
C^{\ast }$ is an $A$-pure dense right $H$-module subalgebra and $_{R}C$ is
locally projective \emph{(}resp. f.g. and projective\emph{)}, then $\beta
(A\#T)\subseteq \#(C,A)$ is dense and
\begin{equation*}
_{A}\mathcal{M}(H)^{C}\simeq \sigma \lbrack _{\#(C,A)}\mathcal{D}]\simeq
\sigma \lbrack _{A\#T}\mathcal{D}]\;\text{\emph{(}resp. }_{A}\mathcal{M}%
(H)^{C}\simeq \text{ }_{\#(C,A)}\mathcal{M}\simeq \text{ }_{A\#C^{\ast }}%
\mathcal{M}\text{\emph{)}.}
\end{equation*}
\end{corollary}

Next we extend some results of \cite{MSTW2001} on relative Hopf modules to
the general case of right-right Doi-Koppinen modules.

\begin{proposition}
\label{rat-dk}Let $(H,A,C)$ be a right-right Doi-Koppinen structure such
that $_{R}C$ is locally projective, $\mathcal{C}:=A\otimes _{R}C$ the
corresponding $A$-coring and $T\subseteq C^{\ast }$ an $A$-pure dense left $%
H $-module subalgebra.

\begin{enumerate}
\item  Let $M$ be a right $A$-module and a left $T$-module. If for all $f\in
T,$ $a\in A$ and $m\in M$ we have
\begin{equation*}
f[ma]=\sum ((a_{<1>}f)m)a_{<0>},
\end{equation*}
then $\mathrm{Rat}^{C}(_{T}M)\in \mathcal{M}(H)_{A}^{C}.$ Consequently $M\in
\mathcal{M}(H)_{A}^{C}$ iff $M=\mathrm{Rat}^{C}(_{T}M).$

\item  If $\varepsilon _{C}\in T,$ then for every right $A\#^{op}T$-module $%
M $ we have: $\mathrm{Rat}^{C}(_{T}M)=\mathrm{Rat}^{\mathcal{C}%
}(M_{A\#^{op}T}) $ and $M\in \mathcal{M}(H)_{A}^{C}$ iff $M=\mathrm{Rat}%
^{C}(_{T}M).$
\end{enumerate}
\end{proposition}

\begin{Beweis}
\begin{enumerate}
\item  Since $_{R}C$ is locally projective and $T\subseteq C^{\ast }$ is
dense, it follows by \cite[Satz 2.2.13]{Abu2001} that $\mathrm{Rat}^{C}(_{T}%
\mathcal{M})\simeq \mathcal{M}^{C}.$ Moreover we have for all $m\in \mathrm{%
Rat}^{C}(_{T}M),$ $f\in T$ and $a\in A:$%
\begin{equation*}
\begin{tabular}{lllll}
$f[ma]$ & $=$ & $\sum ((a_{<1>}f)m)a_{<0>}$ & $=$ & $\sum
(m_{<0>}(a_{<1>}f)(m_{<1>}))a_{<0>}$ \\
& $=$ & $\sum f(m_{<1>}a_{<1>})m_{<0>}a_{<0>},$ &  &
\end{tabular}
\end{equation*}
i.e. $ma\in \mathrm{Rat}^{C}(_{T}M)$ with $\varrho _{M}(ma)=\sum
m_{<0>}a_{<0>}\otimes m_{<1>}a_{<1>},$ hence $\mathrm{Rat}^{C}(_{T}M)\in
\mathcal{M}(H)_{A}^{C}.$ On the other hand, if $M\in \mathcal{M}(H)_{A}^{C},$
then $M$ is in particular a right $C$-comodule and so $M=\mathrm{Rat}%
^{C}(_{T}M).$

\item  Clearly $\mathrm{Rat}^{\mathcal{C}}(M_{A\#^{op}T})\subseteq \mathrm{%
Rat}^{C}(_{T}M).$ On the other hand we have for all $f\in T,$ $a\in A$ and $%
m\in \mathrm{Rat}^{C}(_{T}M):$%
\begin{equation*}
\begin{tabular}{lllll}
$m(a\#f)$ & $=$ & $m((1_{A}\#f)\cdot (a\#\varepsilon _{C}))$ & $=$ & $%
(m(1_{A}\#f))(a\#\varepsilon _{C})$ \\
& $=$ & $(fm)(a\#\varepsilon _{C})$ & $=$ & $(\sum
m_{<0>}f(m_{<1>}))(a\#\varepsilon _{C}))$ \\
& $=$ & $\sum m_{<0>}af(m_{<1>}),$ &  &
\end{tabular}
\end{equation*}
i.e. $m\in \mathrm{Rat}^{\mathcal{C}}(M_{A\#^{op}T})$ with $\varrho
_{M}(m)=\sum m_{<0>}\otimes _{A}(1_{A}\otimes m_{<1>}).$ Note that for all $%
f\in T,$ $a\in A$ and $m\in \mathrm{Rat}^{C}(_{T}M)$ we have by a similar
argument that $f[ma]=\sum ((a_{<1>}f)m)a_{<0>}$ and the result follows by
(1).$\blacksquare $\newpage
\end{enumerate}
\end{Beweis}

\qquad As a direct consequence of Proposition \ref{rat-dk} we get

\begin{corollary}
Let $(H,A,C)\;$be a right-right \emph{(}resp. a left-right\emph{)}
Doi-Koppinen structure and assume $_{R}C$ to be locally projective. If $%
\mathcal{M}^{C}$ is closed under extensions in $_{C^{\ast }}\mathcal{M},$
then $\mathcal{M}(H)_{A}^{C}$ \emph{(}resp. $_{A}\mathcal{M}(H)^{C}$\emph{)}%
\ is closed under extensions in $\mathcal{M}_{A\#^{op}C^{\ast }}$ \emph{(}%
resp. in $_{A\#C^{\ast }}\mathcal{M}$\emph{)}.
\end{corollary}

The proof of the following result is with slight modifications along the
lines of \cite[1.9]{MSTW2001}.

\begin{corollary}
\label{sub-d}Let $(H,A,C)$ be a right-right \emph{(}resp. a left-right\emph{)%
} Doi-Koppinen structure with $_{R}C$ locally projective and consider the
corresponding $A$-coring $\mathcal{C}:=A\otimes _{R}C$ \emph{(}resp. the $%
A^{op}$-coring $\mathcal{D}:=A^{op}\otimes _{R}C$\emph{)}. If $C$ is left
coproper and $C^{\Box }\subseteq C^{\ast }$ is $A$-pure, then
\begin{equation*}
\mathcal{M}(H)_{A}^{C}=\mathcal{M}_{A\#^{op}C^{\Box }}=\mathcal{M}_{^{\Box }%
\mathcal{C}}\;\text{\emph{(}resp. }_{A}\mathcal{M}(H)^{C}\simeq \text{ }%
_{A\#C^{\Box }}\mathcal{M}\simeq \text{ }_{(^{\Box }\mathcal{D})^{op}}%
\mathcal{M}\text{\emph{)}.}
\end{equation*}
\end{corollary}

\begin{Beweis}
Let $(H,A,C)$ be a right-right Doi-Koppinen structure, $\mathcal{A}%
:=A\#^{op}C^{\ast },$ $P:=(\mathcal{A},\mathcal{C})\in \mathcal{P}_{ml}$ and
assume $_{R}C$ to be locally projective. Since $C$ is left coproper, it
follows by the isomorphism of categories $\mathrm{Rat}^{C}(_{C^{\ast }}%
\mathcal{M})=\sigma \lbrack _{C^{\ast }}C]$ and analog to \cite[2.6]{Wis99}
that $\mathrm{Rat}^{C}(_{C^{\ast }}N)=C^{\Box }N$ for every left $C^{\ast }$%
-module $N.$ By Lemma \ref{rat-dk} (2) we have then
\begin{equation*}
\mathcal{T}:=\mathrm{Rat}^{\mathcal{C}}(\mathcal{A}_{\mathcal{A}})=\mathrm{%
Rat}^{C}(_{C^{\ast }}\mathcal{A})=C^{\Box }\mathcal{A}=\mathcal{A}%
(1_{A}\#^{op}C^{\Box })=A\#^{op}C^{\Box }.
\end{equation*}
Since $C^{\Box }\subset C^{\ast }$ is dense, $\mathcal{T}=A\otimes
_{R}C^{\Box }\subset A\otimes _{R}C^{\ast }$ is dense by Theorem \ref
{HAC-iso} and the result follows by Proposition \ref{coprop}. The
corresponding result for left-right Doi-Koppinen structures follows by
symmetry.$\blacksquare $
\end{Beweis}

\qquad Next we consider three examples of right-right Doi-Koppinen
structures (see \cite{CMZ2002}).

\begin{punto}
\textbf{Yetter-Drinfel'd modules. }Let $(H,K,A,C)$ be a \emph{%
Yetter-Drinfel'd datum,} $(K^{op}\otimes _{R}H,A,C)$ the corresponding
right-right Doi-Koppinen structure and consider the category of
Yetter-Drinfel'd modules $\mathcal{YD}(K,H)_{A}^{C}\simeq \mathcal{M}%
(K^{op}\otimes _{R}H)_{A}^{C}.$ If $_{R}C$ is flat, then $\mathcal{YD}%
(K,H)_{A}^{C}$ is a Grothendieck category with enough injective objects
(this generalizes \cite[Section 4.4., Corollary 31]{CMZ2002}, where a base
field is assumed). If $_{R}C$ is locally projective (resp. f.g. and
projective), then we have with $\mathcal{C}:=A\otimes _{R}C$ the
corresponding $A$-coring
\begin{equation*}
\mathcal{YD}(K,H)_{A}^{C}\simeq \sigma \lbrack \mathcal{C}%
_{\#^{op}(C,A)}]\simeq \sigma \lbrack \mathcal{C}_{A\#^{op}C^{\ast })}]\text{
(resp. }\mathcal{YD}(K,H)_{A}^{C}\simeq \mathcal{M}_{\#^{op}(C,A)}\simeq
\mathcal{M}_{A\#^{op}C^{\ast }}\text{).}
\end{equation*}
\end{punto}

\begin{punto}
\textbf{Long dimodules.} Let $A$ be an $R$-algebra, $C$ an $R$-coalgebra, $%
(R,A,C)$ the trivial right-right Doi-Koppinen structure and consider the
category of \emph{Long dimodules} $\mathcal{L}_{A}^{C}\simeq \mathcal{M}%
(R)_{A}^{C}.$ If $_{R}C$ is flat, then $\mathcal{L}_{A}^{C}$ is a
Grothendieck category with enough injective objects. If moreover $_{R}C$ is
locally projective (resp. f.g. and projective), then we have with $\mathcal{C%
}:=A\otimes _{R}C$ the corresponding $A$-coring
\begin{equation*}
\mathcal{L}_{A}^{C}\simeq \sigma \lbrack \mathcal{C}_{\#^{op}(C,A)}]\simeq
\sigma \lbrack \mathcal{C}_{A\#^{op}C^{\ast }}]\text{ (resp. }\mathcal{L}%
_{A}^{C}\simeq \mathcal{M}_{\#^{op}(C,A)}=\mathcal{M}_{A\#^{op}C^{\ast }}%
\text{).}
\end{equation*}
\end{punto}

\begin{punto}
\textbf{Modules graded by }$G$\textbf{-sets}. Let $G$ be a group, $A$ a $G$%
-graded $R$-algebra, $X$ a right $G$-set (e.g. $X=G$), $(RG,A,RX)$ the
corresponding right-right Doi-Koppinen structure and denote by $gr$-$%
(G,A,X)\simeq \mathcal{M}(RG)_{A}^{RX}$ the category of $RX$-graded right $A$%
-modules. Since the free $R$-module $RX$ is in particular locally
projective, we get by Theorem \ref{HAC-iso} (2) isomorphisms of categories
\begin{equation*}
gr\text{-}(G,A,X)\simeq \sigma \lbrack (A\otimes
_{R}RX)_{\#^{op}(RX,A)}]\simeq \sigma \lbrack (A\otimes
_{R}RX)_{A\#^{op}(RX)^{\ast }}].
\end{equation*}
If moreover $X$ is finite, then $RX$ is in particular f.g. and projective,
hence
\begin{equation*}
gr\text{-}(G,A,X)\simeq \mathcal{M}_{\#^{op}(RX,A)}\simeq \mathcal{M}%
_{A\#^{op}(RX)^{\ast }}.
\end{equation*}
\end{punto}

\subsection*{Alternative Doi-Koppinen modules}

It turns out from work of D. Tambara \cite{Tam90} that every entwining
structure $(A,C,\psi ),$ for which $_{R}A$ is f.g. and projective, can be
obtained from a Doi-Koppinen structure with a suitable auxiliary $R$%
-bialgebra giving rise to the entwining map $\psi .$ P. Schauenburg has
shown in \cite{Sch2000} that this is not the case in general. However, if $%
_{R}C$ is f.g. and projective, then he remarks that $(A,C,\psi )$ can be
derived form what he calls an \emph{alternative Doi-Koppinen structure.}

\begin{punto}
Let $H$ be an $R$-bialgebra, $A$ a right $H$-module algebra and $C$ a right $%
H$-comodule coalgebra. Then $(H,A,C)$ is called a \emph{right-right
alternative Doi-Koppinen structure}. It turns out, that $(A,C,\psi )$ is a
right-right entwining structure, where
\begin{equation*}
\psi :C\otimes _{R}A\rightarrow A\otimes _{R}C,\text{ }c\otimes a\mapsto
\sum ac_{<1>}\otimes c_{<0>}.
\end{equation*}
We denote the corresponding category of entwined modules (called \emph{%
alternative right-right Doi-Koppinen modules}) by $a\mathcal{M}(H)_{A}^{C}.$
As for other categories of entwined moudles, if $_{R}C$ is flat, then $a%
\mathcal{M}(H)_{A}^{C}$ is a Grothendieck category with enough injective
objects. If moreover $_{R}C$ is locally projective (resp. f.g. and
projective), then
\begin{equation*}
a\mathcal{M}(H)_{A}^{C}\simeq \sigma \lbrack (A\otimes _{R}C)_{\#^{op}(C,A)}]%
\text{ (resp. }a\mathcal{M}(H)_{A}^{C}\simeq \mathcal{M}_{\#^{op}(C,A)}).
\end{equation*}
\end{punto}

\textbf{Acknowledgments.} Most of the results in this paper are
generalizations of results in my dissertation at the Heinrich-Heine
Universit\"{a}t (D\"{u}sseldorf - Deutschland). I am so grateful to my
supervisor Prof. Robert Wisbauer for the continuous support and
encouragement. I also thank Tomasz Brzezi\'{n}ski for drawing my attention
to the theory of corings and entwined modules during my visit to him in
Swansea and for example (\ref{Tomasz}). Many thanks go to Jos\'{e}
G\'{o}mez-Torrecillas for his inspiring ideas and for the useful preprints
on the subject he sent me.


\begin{thebibliography}{Abu05}
\bibitem{Abu}  J.Y. Abuhlail, \emph{On the linear weak topology and
dual pairings over rings}, Topology and its Applications \textbf{149},
161-175 (2005).

\bibitem{Abu2001}  J.Y. Abuhlail, \emph{Dualit\"{a}tstheoreme for
Hopf-Algebren \"{u}ber Ringen}, \textbf{Ph.D. Dissertation,} Heinrich-Heine
Universit\"{a}t, D\"{u}sseldorf - Germany (2001).
http://www.ulb.uni-duesseldorf.de/diss/mathnat/2001/abuhlail.html

\bibitem{AG-TL2001}  J.Y. Abuhlail, J. G\'{o}mez-Torrecillas and F. Lobillo,
\emph{Duality and rational modules for Hopf algebras over commutative rings}%
, J. Algebra \textbf{240}, 165-184 (2001).

\bibitem{AT78}  H. Allen and D. Trushin, \emph{Coproper coalgebras}, J.
Algebra \textbf{54}, 203-215 (1978).

\bibitem{AW97}  T. Albu and R. Wisbauer, $M$\emph{-density, }$M$\emph{-adic
completion and }$M$\emph{-subgeneration}, Rend. Semin. Mat. Univ. Padova
\textbf{98}, 141-159 (1997).

\bibitem{Ber94}  J. Berning, \emph{Beziehungen zwischen links-linearen
Toplogien und Modulkategorien}, \textbf{Dissertation}, Heinrich-Heine
Universit\"{a}t, D\"{u}sseldorf - Germany (1994).

\bibitem{BM98}  T. Brzezi\'{n}ski and S. Majid, \emph{Coalgebra bundles},
Comm. Math. Phys. \textbf{191}, 467-492 (1998).

\bibitem{Brz}  T. Brzezi\'{n}ski, \emph{The structure of corings. Induction
functors, Maschke-type theorem, and Frobenius and Galois-type properties},
Algebr. Represent. Theory \textbf{5, }389-410 (2002).

\bibitem{Brz99}  T. Brzezi\'{n}ski, \emph{On modules associated to coalgebra
Galois extensions}, J. Algebra \textbf{215} (1999), 290-317.

\bibitem{Bou74}  N. Bourbaki, \emph{Elements of Mathematics, Algebra I,
Chapters }\textbf{1-3}, Hermann (1974).

\bibitem{CC94}  C. Cai and H. Chen, \emph{Coactions, smash products and Hopf
modules}, J. Algebra \textbf{167}, 85-89 (1994).

\bibitem{CMZ2002}  S.~Caenepeel, G.~Militaru, and S. Zhu, \emph{Frobenius
and Separable Functors for Generalized Module Categories and Nonlinear
Equations}{\normalsize , }Lect. Not. Math.{\normalsize \ }\textbf{1787},
Springer-Verlag, Berlin (2002).

\bibitem{DNR2001}  S. D\v{a}sc\v{a}lescu, C. N\v{a}st\v{a}sescu and \c{S}.
Raianu, \emph{Hopf Algebras: an Introduction}, Pure and Applied Mathematics
\textbf{235}, Marcel Dekker, New York (2001).

\bibitem{Doi94}  Y. Doi, \emph{Generalized smash products and Morita
contexts for arbitrary Hopf algebras}, J. Bergen (ed.) et al., Advances in
Hopf Algebras, Lect. Notes Pure Appl. Math. \textbf{158}, Marcle Dekker, New
York, 39-53 (1994).

\bibitem{Doi92}  Y. Doi, \emph{Unifying Hopf modules}, J. Algebra \textbf{153%
}, 373-385 (1992).

\bibitem{Doi83}  Y. Doi, \emph{On the structure of relative Hopf modules},
Comm. Algebra \textbf{11}, 243-255 (1983).

\bibitem{EG-TL}  L. El Kaoutit, J. G\'{o}mez-Torrecillas and F.J. Lobillo,
\emph{Semisimple corings}, preprint, to appear in Algebra Colloquium.

\bibitem{Gar76}  G. Garfinkel, \emph{Universally torsionless and trace
modules}, J. Amer. Math. Soc. \textbf{215}, 119-144 (1976).

\bibitem{G-T98}  J. G\'{o}mez-Torrecillas, \emph{Coalgebras and comodules
over a commutative ring}, Rom. J. Pure Appl. Math. \textbf{43}, 591-603
(1998).

\bibitem{G-TN95}  J. G\'{o}mez-Torrecillas and C. N\u{a}st\u{a}sescu, \emph{%
Quasi-co-Frobenius coalgebras}, J. Algebra \textbf{174}, 909-923 (1995).

\bibitem{Guz89}  F. Guzman, \emph{Cointegrations, relative cohomology for
comodules and coseparable corings}, J. Algebra \textbf{126}, 211-224 (1989).

\bibitem{Guz85}  F. Guzman, $\emph{Cointegratin}$ \emph{and relative
cohomology for comodules}, \textbf{Ph.D. Dissertation}, Syracuse University
USA (1985).

\bibitem{Kot66}  G. K\"{o}the, \emph{Topologische lineare R\"{a}ume }$I$,
Die Grundlagen der mathematischen Wissenschaften \textbf{107}, Berlin:
Springer-Verlag (1966).

\bibitem{Kop95}  M. Koppinen, \emph{Variations on the smash product with
applications to group-graded rings}, J. Pure Appl. Algebra \textbf{104},
61-80 (1995).

\bibitem{Kop92}  M. Koppinen, \emph{A duality theorem for crossed products
of Hopf algebras, }J. Algebra \textbf{146}, 153-174 (1992).

\bibitem{MSTW2001}  C. Menini, A. Seidel, B. Torrecillas and R. Wisbauer, $A$%
-$H$\emph{-bimodules and equivalences}, Comm. Algebra \textbf{29}, 4619-4640
(2001).

\bibitem{MTW2001}  C. Menini, B. Torrecillas and R. Wisbauer, \emph{Strongly
rational comodules and semiperfect Hopf algebras over QF rings}, J. Pure
Appl. Algebra \textbf{155}, 237-255\ (2001).

\bibitem{NS82}  W. Nichols and M. Sweedler, \emph{Hopf algebras and
combinatorics}, Contemp. Math. \textbf{6}, 49-84 (1982).

\bibitem{MZ97}  C. Menini and M. Zuccoli, \emph{Equivalence theorems and
Hopf-Galois extensions}, J. Algebra \textbf{194}, 245-274 (1997).

\bibitem{Rad73}  D. Radford, \emph{Coreflexive coalgebras, }J. Algebra
\textbf{26}, 512-535 (1973).

\bibitem{Sch2000}  P. Schauenburg, \emph{Doi-Koppinen Hopf modules versus
entwined modules}, New York J. Math. \textbf{6}, 325-329 (2000).

\bibitem{Sch72}  H. Schubert, \emph{Categories}, Springer-Verlag (1972).

\bibitem{Swe75}  M. Sweedler, \emph{The predual theorem to the
Jacobson-Bourbaki theorem}, Trans. Amer. Math. Soc. \textbf{213}, 391-406
(1975).

\bibitem{Swe69}  M. Sweedler, \emph{Hopf Algebras}, New York: Benjamin,
(1969).

\bibitem{Tam90}  D. Tambara, \emph{The coendomorphism bialgebra of an algebra%
}, J. Fac. Sci., Univ. Tokyo, Sect. I A \textbf{37}, 425-456 (1990).

\bibitem{Wis2002}  R. Wisbauer, \emph{On the category of comodules for
corings}, Proc. 3rd. Int. Pal. Conf.: Math. $\&$ Math. Edu., Bethlehem
(Palestine), S. Elaydi et al. (ed.), World Scientific, New Jersy, ISBN
981-02-4720-6, 325-336 (2002).

\bibitem{Wis2001}  R. Wisbauer, \emph{Weak corings}, J. Algebra \textbf{245}%
, 123-160 (2001).

\bibitem{Wis99}  R. Wisbauer, \emph{Semiperfect coalgebras over rings},
Algebra and Combinatorics. Papers from the International Congress ICAC'97
Hongkong, K.-P. Shum et al. (ed.), Singapore:\ Springer-Verlag, 487-512
(1999).

\bibitem{Wis88}  R.\thinspace Wisbauer, \emph{Grundlagen der Modul- und
Ringtheorie}, M\"{u}nchen: Verlag Reinhard Fischer (1988); \emph{Foundations
of Module and Ring Theory}, Gordon and Breach, Reading (1991).

\bibitem{Wisch75}  M. Wischnewsky, \emph{On linear representations of affine
groups. I}, Pac. J. Math. \textbf{61}, 551-572 (1975).

\bibitem{Z-H76}  B. Zimmermann-Huisgen, \emph{Pure submodules of direct
products of free modules}, Math. Ann. \textbf{224}, 233-245 (1976).
\end{thebibliography}
\end{document}